\documentclass[12pt,a4paper]{article}

\usepackage{amsfonts,amssymb,amsthm,amsmath,latexsym}
\usepackage{eqsection,indent}
\usepackage{subeqnarray, eepicemu, url,cite,bm}
\usepackage{multirow}  
\usepackage{tensor}  
\usepackage{graphicx,graphbox}
\usepackage{float}  
\usepackage{mycaption}

\date{July 4, 2019 \\[2mm] revised August 13, 2020}   

\begin{document}

\title{\vspace*{-2cm}
       Lattice paths and branched continued fractions  \\[5mm]
       II.~Multivariate Lah polynomials \\ and Lah symmetric functions
      }

\author{ \\
      \hspace*{-1cm}
      {\large Mathias P\'etr\'eolle${}^1$ and Alan D.~Sokal${}^{1,2}$}
   \\[5mm]
     \hspace*{-1.3cm}
      \normalsize
           ${}^1$Department of Mathematics, University College London,
                    London WC1E 6BT, UK   \\[1mm]
     \hspace*{-2.9cm}
      \normalsize
           ${}^2$Department of Physics, New York University,
                    New York, NY 10003, USA
       \\[5mm]
     \hspace*{-0.5cm}
     {\tt mathias.petreolle@gmail.com}, {\tt sokal@nyu.edu}  \\[1cm]
}

\maketitle
\thispagestyle{empty}   

\begin{abstract}
We introduce the generic Lah polynomials $L_{n,k}(\bm{\phi})$,
which enumerate unordered forests of increasing ordered trees
with a weight $\phi_i$ for each vertex with $i$ children.
We show that, if the weight sequence $\bm{\phi}$ is Toeplitz-totally positive,
then the triangular array of generic Lah polynomials is totally positive
and the sequence of row-generating polynomials $L_n(\bm{\phi},y)$
is coefficientwise Hankel-totally positive.
Upon specialization we obtain results for the
Lah symmetric functions and multivariate Lah polynomials
of positive and negative type.
The multivariate Lah polynomials of positive type
are also given by a branched continued fraction.
Our proofs use mainly the method of production matrices;
the production matrix is obtained by a bijection
from ordered forests of increasing ordered trees
to labeled partial \L{}ukasiewicz paths.
We also give a second proof of the continued fraction
using the Euler--Gauss recurrence method.
\end{abstract}

\bigskip
\noindent
{\bf Key Words:}
Lah polynomial, Bell polynomial, Eulerian polynomial, symmetric function,
increasing tree, forest, \L{}ukasiewicz path, Stirling permutation,
continued fraction, branched continued fraction,
production matrix, Toeplitz matrix, Hankel matrix,
total positivity, Toeplitz-total positivity, Hankel-total positivity.

\bigskip
\bigskip
\noindent
{\bf Mathematics Subject Classification (MSC 2010) codes:}
05A15 (Primary);
05A05, 05A18, 05A19, 05A20, 05C30, 05E05, 11B37, 11B73,
15B48, 30B70 (Secondary).

\vspace*{1cm}

\newtheorem{theorem}{Theorem}[section]
\newtheorem{proposition}[theorem]{Proposition}
\newtheorem{lemma}[theorem]{Lemma}
\newtheorem{corollary}[theorem]{Corollary}
\newtheorem{definition}[theorem]{Definition}
\newtheorem{conjecture}[theorem]{Conjecture}
\newtheorem{question}[theorem]{Question}
\newtheorem{problem}[theorem]{Problem}
\newtheorem{openproblem}[theorem]{Open Problem}
\newtheorem{example}[theorem]{Example}

\renewcommand{\theenumi}{\alph{enumi}}
\renewcommand{\labelenumi}{(\theenumi)}
\def\eop{\hbox{\kern1pt\vrule height6pt width4pt
depth1pt\kern1pt}\medskip}
\def\prf{\par\noindent{\bf Proof.\enspace}\rm}
\def\rmk{\par\medskip\noindent{\bf Remark\enspace}\rm}

\newcommand{\textbfit}[1]{\textbf{\textit{#1}}}

\newcommand{\bigdash}{%
\smallskip\begin{center} \rule{5cm}{0.1mm} \end{center}\smallskip}

\newcommand{\safepar}{ {\protect\hfill\protect\break\hspace*{5mm}} }

\newcommand{\be}{\begin{equation}}
\newcommand{\ee}{\end{equation}}
\newcommand{\<}{\langle}
\renewcommand{\>}{\rangle}
\newcommand{\widebar}{\overline}
\def\reff#1{(\protect\ref{#1})}
\def\spose#1{\hbox to 0pt{#1\hss}}
\def\ltapprox{\mathrel{\spose{\lower 3pt\hbox{$\mathchar"218$}}
    \raise 2.0pt\hbox{$\mathchar"13C$}}}
\def\gtapprox{\mathrel{\spose{\lower 3pt\hbox{$\mathchar"218$}}
    \raise 2.0pt\hbox{$\mathchar"13E$}}}
\def\textprime{${}^\prime$}
\def\proof{\par\medskip\noindent{\sc Proof.\ }}
\def\firstproof{\par\medskip\noindent{\sc First Proof.\ }}
\def\secondproof{\par\medskip\noindent{\sc Second Proof.\ }}
\def\alternateproof{\par\medskip\noindent{\sc Alternate Proof.\ }}
\def\algebraicproof{\par\medskip\noindent{\sc Algebraic Proof.\ }}
\def\combinatorialproof{\par\medskip\noindent{\sc Combinatorial Proof.\ }}
\def\proofof#1{\bigskip\noindent{\sc Proof of #1.\ }}
\def\firstproofof#1{\bigskip\noindent{\sc First Proof of #1.\ }}
\def\secondproofof#1{\bigskip\noindent{\sc Second Proof of #1.\ }}
\def\thirdproofof#1{\bigskip\noindent{\sc Third Proof of #1.\ }}
\def\alternateproofof#1{\bigskip\noindent{\sc Alternate Proof of #1.\ }}
\def\algebraicproofof#1{\bigskip\noindent{\sc Algebraic Proof of #1.\ }}
\def\combinatorialproofof#1{\bigskip\noindent{\sc Combinatorial Proof of #1.\ }}
\def\sketchofproof{\par\medskip\noindent{\sc Sketch of proof.\ }}
\renewcommand{\qed}{ $\square$ \bigskip}
\newcommand{\myendremark}{ $\blacksquare$ \bigskip}
\def\half{ {1 \over 2} }
\def\third{ {1 \over 3} }
\def\twothird{ {2 \over 3} }
\def\smfrac#1#2{{\textstyle{#1\over #2}}}
\def\smhalf{ {\smfrac{1}{2}} }
\newcommand{\real}{\mathop{\rm Re}\nolimits}
\renewcommand{\Re}{\mathop{\rm Re}\nolimits}
\newcommand{\imag}{\mathop{\rm Im}\nolimits}
\renewcommand{\Im}{\mathop{\rm Im}\nolimits}
\newcommand{\sgn}{\mathop{\rm sgn}\nolimits}
\newcommand{\tr}{\mathop{\rm tr}\nolimits}
\newcommand{\supp}{\mathop{\rm supp}\nolimits}
\newcommand{\disc}{\mathop{\rm disc}\nolimits}
\newcommand{\diag}{\mathop{\rm diag}\nolimits}
\newcommand{\tridiag}{\mathop{\rm tridiag}\nolimits}
\newcommand{\AZ}{\mathop{\rm AZ}\nolimits}
\newcommand{\NC}{\mathop{\rm NC}\nolimits}
\newcommand{\PF}{{\rm PF}}
\newcommand{\rk}{\mathop{\rm rk}\nolimits}
\newcommand{\perm}{\mathop{\rm perm}\nolimits}
\def\hboxscript#1{ {\hbox{\scriptsize\em #1}} }
\renewcommand{\emptyset}{\varnothing}
\newcommand{\eqdef}{\stackrel{\rm def}{=}}

\newcommand{\restrict}{\upharpoonright}

\newcommand{\compinv}{{\langle -1 \rangle}}   

\newcommand{\scra}{{\mathcal{A}}}
\newcommand{\scrb}{{\mathcal{B}}}
\newcommand{\scrc}{{\mathcal{C}}}
\newcommand{\scrd}{{\mathcal{D}}}
\newcommand{\scrdtilde}{{\widetilde{\mathcal{D}}}}
\newcommand{\scre}{{\mathcal{E}}}
\newcommand{\scrf}{{\mathcal{F}}}
\newcommand{\scrg}{{\mathcal{G}}}
\newcommand{\scrh}{{\mathcal{H}}}
\newcommand{\scri}{{\mathcal{I}}}
\newcommand{\scrj}{{\mathcal{J}}}
\newcommand{\scrk}{{\mathcal{K}}}
\newcommand{\scrl}{{\mathcal{L}}}
\newcommand{\scrm}{{\mathcal{M}}}
\newcommand{\scrn}{{\mathcal{N}}}
\newcommand{\scro}{{\mathcal{O}}}
\newcommand\scroo{
  \mathchoice
    {{\scriptstyle\mathcal{O}}}
    {{\scriptstyle\mathcal{O}}}
    {{\scriptscriptstyle\mathcal{O}}}
    {\scalebox{0.6}{$\scriptscriptstyle\mathcal{O}$}}
  }
\newcommand{\scrp}{{\mathcal{P}}}
\newcommand{\scrq}{{\mathcal{Q}}}
\newcommand{\scrr}{{\mathcal{R}}}
\newcommand{\scrs}{{\mathcal{S}}}
\newcommand{\scrt}{{\mathcal{T}}}
\newcommand{\scrv}{{\mathcal{V}}}
\newcommand{\scrw}{{\mathcal{W}}}
\newcommand{\scrz}{{\mathcal{Z}}}
\newcommand{\SP}{{\mathcal{SP}}}
\newcommand{\ST}{{\mathcal{ST}}}

\newcommand{\bfa}{{\mathbf{a}}}
\newcommand{\bfb}{{\mathbf{b}}}
\newcommand{\bfc}{{\mathbf{c}}}
\newcommand{\bfd}{{\mathbf{d}}}
\newcommand{\bfe}{{\mathbf{e}}}
\newcommand{\bfh}{{\mathbf{h}}}
\newcommand{\bfj}{{\mathbf{j}}}
\newcommand{\bfi}{{\mathbf{i}}}
\newcommand{\bfk}{{\mathbf{k}}}
\newcommand{\bfl}{{\mathbf{l}}}
\newcommand{\bfL}{{\mathbf{L}}}
\newcommand{\bfm}{{\mathbf{m}}}
\newcommand{\bfn}{{\mathbf{n}}}
\newcommand{\bfp}{{\mathbf{p}}}
\newcommand{\bfr}{{\mathbf{r}}}
\newcommand{\bfu}{{\mathbf{u}}}
\newcommand{\bfv}{{\mathbf{v}}}
\newcommand{\bfw}{{\mathbf{w}}}
\newcommand{\bfx}{{\mathbf{x}}}
\newcommand{\bfX}{{\mathbf{X}}}
\newcommand{\bfy}{{\mathbf{y}}}
\newcommand{\bfz}{{\mathbf{z}}}
\renewcommand{\k}{{\mathbf{k}}}
\newcommand{\n}{{\mathbf{n}}}
\newcommand{\vv}{{\mathbf{v}}}
\newcommand{\bv}{{\mathbf{v}}}
\newcommand{\w}{{\mathbf{w}}}
\newcommand{\x}{{\mathbf{x}}}
\newcommand{\y}{{\mathbf{y}}}
\newcommand{\cc}{{\mathbf{c}}}
\newcommand{\zero}{{\mathbf{0}}}
\newcommand{\one}{{\mathbf{1}}}
\newcommand{\bmm}{{\mathbf{m}}}

\newcommand{\ahat}{{\widehat{a}}}
\newcommand{\Zhat}{{\widehat{Z}}}

\newcommand{\C}{{\mathbb C}}
\newcommand{\D}{{\mathbb D}}
\newcommand{\Z}{{\mathbb Z}}
\newcommand{\N}{{\mathbb N}}
\newcommand{\Q}{{\mathbb Q}}
\newcommand{\PP}{{\mathbb P}}
\newcommand{\R}{{\mathbb R}}
\newcommand{\RR}{{\mathbb R}}
\newcommand{\E}{{\mathbb E}}

\newcommand{\Sym}{{\mathfrak{S}}}
\newcommand{\SymB}{{\mathfrak{B}}}
\newcommand{\Alt}{{\mathrm{Alt}}}

\newcommand{\germanA}{{\mathfrak{A}}}
\newcommand{\germanB}{{\mathfrak{B}}}
\newcommand{\germanQ}{{\mathfrak{Q}}}
\newcommand{\germanh}{{\mathfrak{h}}}

\newcommand{\myle}{\preceq}
\newcommand{\myge}{\succeq}
\newcommand{\mygt}{\succ}

\newcommand{\B}{{\sf B}}
\newcommand{\OB}{B^{\rm ord}}
\newcommand{\OS}{{\sf OS}}
\newcommand{\OO}{{\sf O}}
\newcommand{\OSP}{{\sf OSP}}
\newcommand{\Eu}{{\sf Eu}}
\newcommand{\ERR}{{\sf ERR}}
\newcommand{\sfB}{{\sf B}}
\newcommand{\sfD}{{\sf D}}
\newcommand{\sfE}{{\sf E}}
\newcommand{\sfG}{{\sf G}}
\newcommand{\sfJ}{{\sf J}}
\newcommand{\sfL}{{\sf L}}
\newcommand{\sfLhat}{{\widehat{{\sf L}}}}
\newcommand{\sfLtilde}{{\widetilde{{\sf L}}}}
\newcommand{\sfP}{{\sf P}}
\newcommand{\sfQ}{{\sf Q}}
\newcommand{\sfS}{{\sf S}}
\newcommand{\sfT}{{\sf T}}
\newcommand{\sfW}{{\sf W}}
\newcommand{\sfMV}{{\sf MV}}
\newcommand{\AMV}{{\sf AMV}}
\newcommand{\BM}{{\sf BM}}
\newcommand{\emIB}{B^{\rm irr}}
\newcommand{\emIP}{P^{\rm irr}}
\newcommand{\emOB}{B^{\rm ord}}
\newcommand{\emCB}{B^{\rm cyc}}
\newcommand{\emSC}{P^{\rm cyc}}

\newcommand{\lev}{{\rm lev}}
\newcommand{\stat}{{\rm stat}}
\newcommand{\cyc}{{\rm cyc}}
\newcommand{\mysteryone}{{\rm mys1}}
\newcommand{\mysterytwo}{{\rm mys2}}
\newcommand{\Asc}{{\rm Asc}}
\newcommand{\asc}{{\rm asc}}
\newcommand{\Des}{{\rm Des}}
\newcommand{\des}{{\rm des}}
\newcommand{\Exc}{{\rm Exc}}
\newcommand{\exc}{{\rm exc}}
\newcommand{\Wex}{{\rm Wex}}
\newcommand{\wex}{{\rm wex}}
\newcommand{\Fix}{{\rm Fix}}
\newcommand{\fix}{{\rm fix}}
\newcommand{\lrmax}{{\rm lrmax}}
\newcommand{\rlmax}{{\rm rlmax}}
\newcommand{\Rec}{{\rm Rec}}
\newcommand{\rec}{{\rm rec}}
\newcommand{\Arec}{{\rm Arec}}
\newcommand{\arec}{{\rm arec}}
\newcommand{\ERec}{{\rm ERec}}
\newcommand{\erec}{{\rm erec}}
\newcommand{\EArec}{{\rm EArec}}
\newcommand{\earec}{{\rm earec}}
\newcommand{\recarec}{{\rm recarec}}
\newcommand{\nonrec}{{\rm nonrec}}
\newcommand{\Cpeak}{{\rm Cpeak}}
\newcommand{\cpeak}{{\rm cpeak}}
\newcommand{\Cval}{{\rm Cval}}
\newcommand{\cval}{{\rm cval}}
\newcommand{\Cdasc}{{\rm Cdasc}}
\newcommand{\cdasc}{{\rm cdasc}}
\newcommand{\Cddes}{{\rm Cddes}}
\newcommand{\cddes}{{\rm cddes}}
\newcommand{\cdrise}{{\rm cdrise}}
\newcommand{\cdfall}{{\rm cdfall}}
\newcommand{\Peak}{{\rm Peak}}
\newcommand{\peak}{{\rm peak}}
\newcommand{\Val}{{\rm Val}}
\newcommand{\val}{{\rm val}}
\newcommand{\Dasc}{{\rm Dasc}}
\newcommand{\dasc}{{\rm dasc}}
\newcommand{\Ddes}{{\rm Ddes}}
\newcommand{\ddes}{{\rm ddes}}
\newcommand{\inv}{{\rm inv}}
\newcommand{\maj}{{\rm maj}}
\newcommand{\rs}{{\rm rs}}
\newcommand{\cross}{{\rm cr}}
\newcommand{\crosshat}{{\widehat{\rm cr}}}
\newcommand{\nest}{{\rm ne}}
\newcommand{\rodd}{{\rm rodd}}
\newcommand{\reven}{{\rm reven}}
\newcommand{\lodd}{{\rm lodd}}
\newcommand{\leven}{{\rm leven}}
\newcommand{\sg}{{\rm sg}}
\newcommand{\bl}{{\rm bl}}
\newcommand{\tran}{{\rm tr}}
\newcommand{\area}{{\rm area}}
\newcommand{\ret}{{\rm ret}}
\newcommand{\peaks}{{\rm peaks}}
\newcommand{\hl}{{\rm hl}}
\newcommand{\sll}{{\rm sl}}
\newcommand{\negg}{{\rm neg}}
\newcommand{\imp}{{\rm imp}}
\newcommand{\osg}{{\rm osg}}
\newcommand{\ons}{{\rm ons}}
\newcommand{\isg}{{\rm isg}}
\newcommand{\ins}{{\rm ins}}
\newcommand{\LL}{{\rm LL}}
\newcommand{\height}{{\rm ht}}
\newcommand{\as}{{\rm as}}

\newcommand{\ba}{{\bm{a}}}
\newcommand{\bahat}{{\widehat{\bm{a}}}}
\newcommand{\sfa}{{{\sf a}}}
\newcommand{\bb}{{\bm{b}}}
\newcommand{\bc}{{\bm{c}}}
\newcommand{\bchat}{{\widehat{\bm{c}}}}
\newcommand{\bd}{{\bm{d}}}
\newcommand{\bee}{{\bm{e}}}
\newcommand{\beh}{{\bm{eh}}}
\newcommand{\bff}{{\bm{f}}}
\newcommand{\bg}{{\bm{g}}}
\newcommand{\bh}{{\bm{h}}}
\newcommand{\bll}{{\bm{\ell}}}
\newcommand{\bp}{{\bm{p}}}
\newcommand{\br}{{\bm{r}}}
\newcommand{\bs}{{\bm{s}}}
\newcommand{\bu}{{\bm{u}}}
\newcommand{\bw}{{\bm{w}}}
\newcommand{\bx}{{\bm{x}}}
\newcommand{\by}{{\bm{y}}}
\newcommand{\bz}{{\bm{z}}}
\newcommand{\bA}{{\bm{A}}}
\newcommand{\bB}{{\bm{B}}}
\newcommand{\bC}{{\bm{C}}}
\newcommand{\bE}{{\bm{E}}}
\newcommand{\bF}{{\bm{F}}}
\newcommand{\bG}{{\bm{G}}}
\newcommand{\bH}{{\bm{H}}}
\newcommand{\bI}{{\bm{I}}}
\newcommand{\bJ}{{\bm{J}}}
\newcommand{\bL}{{\bm{L}}}
\newcommand{\bLhat}{{\widehat{\bm{L}}}}
\newcommand{\bM}{{\bm{M}}}
\newcommand{\bN}{{\bm{N}}}
\newcommand{\bP}{{\bm{P}}}
\newcommand{\bQ}{{\bm{Q}}}
\newcommand{\bR}{{\bm{R}}}
\newcommand{\bS}{{\bm{S}}}
\newcommand{\bT}{{\bm{T}}}
\newcommand{\bW}{{\bm{W}}}
\newcommand{\bX}{{\bm{X}}}
\newcommand{\bY}{{\bm{Y}}}
\newcommand{\bIB}{{\bm{B}^{\rm irr}}}
\newcommand{\bOB}{{\bm{B}^{\rm ord}}}
\newcommand{\bOS}{{\bm{OS}}}
\newcommand{\bERR}{{\bm{ERR}}}
\newcommand{\bSP}{{\bm{SP}}}
\newcommand{\bMV}{{\bm{MV}}}
\newcommand{\bBM}{{\bm{BM}}}
\newcommand{\balpha}{{\bm{\alpha}}}
\newcommand{\balphapre}{{\bm{\alpha}^{\rm pre}}}
\newcommand{\bbeta}{{\bm{\beta}}}
\newcommand{\bgamma}{{\bm{\gamma}}}
\newcommand{\bdelta}{{\bm{\delta}}}
\newcommand{\bkappa}{{\bm{\kappa}}}
\newcommand{\bmu}{{\bm{\mu}}}
\newcommand{\bomega}{{\bm{\omega}}}
\newcommand{\bsigma}{{\bm{\sigma}}}
\newcommand{\btau}{{\bm{\tau}}}
\newcommand{\bphi}{{\bm{\phi}}}
\newcommand{\bphihat}{{\skew{3}\widehat{\vphantom{t}\protect\smash{\bm{\phi}}}}}
\newcommand{\bpsi}{{\bm{\psi}}}
\newcommand{\bxi}{{\bm{\xi}}}
\newcommand{\bzeta}{{\bm{\zeta}}}
\newcommand{\bone}{{\bm{1}}}
\newcommand{\bzero}{{\bm{0}}}

\newcommand{\Cbar}{{\overline{C}}}
\newcommand{\Dbar}{{\overline{D}}}
\newcommand{\dbar}{{\overline{d}}}
\def\Btilde{{\widetilde{B}}}
\def\Ctilde{{\widetilde{C}}}
\def\Ftilde{{\widetilde{F}}}
\def\Gtilde{{\widetilde{G}}}
\def\Htilde{{\widetilde{H}}}
\def\Lhat{{\widehat{L}}}
\def\Ltilde{{\widetilde{L}}}
\def\Ptilde{{\widetilde{P}}}
\def\ptilde{{\widetilde{p}}}
\def\Chat{{\widehat{C}}}
\def\ctilde{{\widetilde{c}}}
\def\zbar{{\overline{Z}}}
\def\pitilde{{\widetilde{\pi}}}
\def\omegahat{{\widehat{\omega}}}

\newcommand{\sech}{{\rm sech}}

%
%
\newcommand{\sn}{{\rm sn}}
\newcommand{\cn}{{\rm cn}}
\newcommand{\dn}{{\rm dn}}
\newcommand{\sm}{{\rm sm}}
\newcommand{\cm}{{\rm cm}}

%
%
\newcommand{\zfz}{ {{}_0 \! F_0} }
\newcommand{\zfo}{ {{}_0  F_1} }
\newcommand{\ofz}{ {{}_1 \! F_0} }
\newcommand{\ofo}{ {{}_1 \! F_1} }
\newcommand{\oft}{ {{}_1 \! F_2} }

%
%
\newcommand{\FHyper}[2]{ {\tensor[_{#1 \!}]{F}{_{#2}}\!} }
\newcommand{\FHYPER}[5]{ {\FHyper{#1}{#2} \!\biggl(
   \!\!\begin{array}{c} #3 \\[1mm] #4 \end{array}\! \bigg|\, #5 \! \biggr)} }
\newcommand{\tfo}{ {\FHyper{2}{1}} }
\newcommand{\tfz}{ {\FHyper{2}{0}} }
\newcommand{\threefz}{ {\FHyper{3}{0}} }
\newcommand{\FHYPERbottomzero}[3]{ {\FHyper{#1}{0} \hspace*{-0mm}\biggl(
   \!\!\begin{array}{c} #2 \\[1mm] \hbox{---} \end{array}\! \bigg|\, #3 \! \biggr)} }
\newcommand{\FHYPERtopzero}[3]{ {\FHyper{0}{#1} \hspace*{-0mm}\biggl(
   \!\!\begin{array}{c} \hbox{---} \\[1mm] #2 \end{array}\! \bigg|\, #3 \! \biggr)} }

\newcommand{\phiHyper}[2]{ {\tensor[_{#1}]{\phi}{_{#2}}} }
\newcommand{\psiHyper}[2]{ {\tensor[_{#1}]{\psi}{_{#2}}} }
\newcommand{\PhiHyper}[2]{ {\tensor[_{#1}]{\Phi}{_{#2}}} }
\newcommand{\PsiHyper}[2]{ {\tensor[_{#1}]{\Psi}{_{#2}}} }
\newcommand{\phiHYPER}[6]{ {\phiHyper{#1}{#2} \!\left(
   \!\!\begin{array}{c} #3 \\ #4 \end{array}\! ;\, #5, \, #6 \! \right)\!} }
\newcommand{\psiHYPER}[6]{ {\psiHyper{#1}{#2} \!\left(
   \!\!\begin{array}{c} #3 \\ #4 \end{array}\! ;\, #5, \, #6 \! \right)} }
\newcommand{\PhiHYPER}[5]{ {\PhiHyper{#1}{#2} \!\left(
   \!\!\begin{array}{c} #3 \\ #4 \end{array}\! ;\, #5 \! \right)\!} }
\newcommand{\PsiHYPER}[5]{ {\PsiHyper{#1}{#2} \!\left(
   \!\!\begin{array}{c} #3 \\ #4 \end{array}\! ;\, #5 \! \right)\!} }
\newcommand{\zerophizero}{ {\phiHyper{0}{0}} }
\newcommand{\ophizero}{ {\phiHyper{1}{0}} }
\newcommand{\zphio}{ {\phiHyper{0}{1}} }
\newcommand{\ophio}{ {\phiHyper{1}{1}} }
\newcommand{\tphio}{ {\phiHyper{2}{1}} }
\newcommand{\tphiz}{ {\phiHyper{2}{0}} }
\newcommand{\tPhio}{ {\PhiHyper{2}{1}} }
\newcommand{\opsio}{ {\psiHyper{1}{1}} }

%
%
\newcommand{\stirlingsubset}[2]{\genfrac{\{}{\}}{0pt}{}{#1}{#2}}
\newcommand{\stirlingcycle}[2]{\genfrac{[}{]}{0pt}{}{#1}{#2}}
\newcommand{\assocstirlingsubset}[3]{{\genfrac{\{}{\}}{0pt}{}{#1}{#2}}_{\! \ge #3}}
\newcommand{\genstirlingsubset}[4]{{\genfrac{\{}{\}}{0pt}{}{#1}{#2}}_{\! #3,#4}}
\newcommand{\irredstirlingsubset}[2]{{\genfrac{\{}{\}}{0pt}{}{#1}{#2}}^{\!\rm irr}}
\newcommand{\euler}[2]{\genfrac{\langle}{\rangle}{0pt}{}{#1}{#2}}
\newcommand{\eulergen}[3]{{\genfrac{\langle}{\rangle}{0pt}{}{#1}{#2}}_{\! #3}}
\newcommand{\eulersecond}[2]{\left\langle\!\! \euler{#1}{#2} \!\!\right\rangle}
\newcommand{\eulersecondgen}[3]{{\left\langle\!\! \euler{#1}{#2} \!\!\right\rangle}_{\! #3}}
\newcommand{\binomvert}[2]{\genfrac{\vert}{\vert}{0pt}{}{#1}{#2}}
\newcommand{\binomsquare}[2]{\genfrac{[}{]}{0pt}{}{#1}{#2}}
\newcommand{\doublebinom}[2]{\left(\!\! \binom{#1}{#2} \!\!\right)}


\newenvironment{sarray}{
             \textfont0=\scriptfont0
             \scriptfont0=\scriptscriptfont0
             \textfont1=\scriptfont1
             \scriptfont1=\scriptscriptfont1
             \textfont2=\scriptfont2
             \scriptfont2=\scriptscriptfont2
             \textfont3=\scriptfont3
             \scriptfont3=\scriptscriptfont3
           \renewcommand{\arraystretch}{0.7}
           \begin{array}{l}}{\end{array}}

\newenvironment{scarray}{
             \textfont0=\scriptfont0
             \scriptfont0=\scriptscriptfont0
             \textfont1=\scriptfont1
             \scriptfont1=\scriptscriptfont1
             \textfont2=\scriptfont2
             \scriptfont2=\scriptscriptfont2
             \textfont3=\scriptfont3
             \scriptfont3=\scriptscriptfont3
           \renewcommand{\arraystretch}{0.7}
           \begin{array}{c}}{\end{array}}


\newcommand*\circled[1]{\tikz[baseline=(char.base)]{
  \node[shape=circle,draw,inner sep=1pt] (char) {#1};}}
\newcommand{\ostar}{{\circledast}}
\newcommand{\ostarN}{{\,\circledast_{\vphantom{\dot{N}}N}\,}}
\newcommand{\ostarPsi}{{\,\circledast_{\vphantom{\dot{\Psi}}\Psi}\,}}
\newcommand{\starN}{{\,\ast_{\vphantom{\dot{N}}N}\,}}
\newcommand{\starpsi}{{\,\ast_{\vphantom{\dot{\bpsi}}\!\bpsi}\,}}
\newcommand{\starone}{{\,\ast_{\vphantom{\dot{1}}1}\,}}
\newcommand{\startwo}{{\,\ast_{\vphantom{\dot{2}}2}\,}}
\newcommand{\starinfty}{{\,\ast_{\vphantom{\dot{\infty}}\infty}\,}}
\newcommand{\starT}{{\,\ast_{\vphantom{\dot{T}}T}\,}}

\newcommand*{\Scale}[2][4]{\scalebox{#1}{$#2$}}

\newcommand*{\Scaletext}[2][4]{\scalebox{#1}{#2}} 

\clearpage 

\tableofcontents

\clearpage

\section{Introduction and statement of main results}

In a seminal 1980 paper, Flajolet \cite{Flajolet_80}
showed that the coefficients in the Taylor expansion
of the generic Stieltjes-type (resp.\ Jacobi-type) continued fraction
--- which he called the {\em Stieltjes--Rogers}\/
 (resp.\ {\em Jacobi--Rogers}\/) {\em polynomials}\/ ---
can be interpreted as the generating polynomials
for Dyck (resp.\ Motzkin) paths with specified height-dependent weights.
Very recently it was independently discovered by several authors
\cite{Fusy_15,Oste_15,Josuat-Verges_18,Sokal_totalpos}
that Thron-type continued fractions also have an interpretation of this kind:
namely, their Taylor coefficients
--- which we call, by analogy, the {\em Thron--Rogers polynomials}\/ ---
can be interpreted as the generating polynomials 
for Schr\"oder paths with specified height-dependent weights.

In a recent paper \cite{latpath_SRTR}
we presented an infinite sequence
of generalizations of the Stieltjes--Rogers and Thron--Rogers polynomials,
which are parametrized by an integer $m \ge 1$
and reduce to the classical Stieltjes--Rogers and Thron--Rogers polynomials
when $m=1$;
they are the generating polynomials of $m$-Dyck and $m$-Schr\"oder paths,
respectively, with height-dependent weights,
and are also the Taylor coefficients of certain branched continued fractions.
We proved that these generalizations all possess the fundamental property of
coefficientwise Hankel-total positivity \cite{Sokal_flajolet,Sokal_totalpos},
jointly in all the (infinitely many) indeterminates.
These facts were known when $m = 1$ \cite{Sokal_flajolet,Sokal_totalpos}
but were new when $m > 1$.
By specializing the indeterminates we were able to give many examples
of Hankel-totally positive sequences
whose generating functions do not possess nice classical continued fractions.
(The concept of Hankel-total positivity \cite{Sokal_flajolet,Sokal_totalpos}
 will be explained in more detail later in this Introduction.)

In particular, in \cite[section~12]{latpath_SRTR}
we introduced the multivariate Eulerian polynomials
and Eulerian symmetric functions:
these are generating polynomials for increasing trees and forests
of various types (see below for precise definitions),
which vastly extend the classical univariate Eulerian
and $r$th-order Eulerian polynomials;
we proved their coefficientwise Hankel-total positivity.
Here we would like to refine this analysis by considering
(among other things) the row-generating polynomials:
this leads to defining multivariate Lah polynomials
and Lah symmetric functions,
which extend the classical univariate Lah polynomials.
So let us begin by reviewing briefly some well-known univariate
combinatorial polynomials;
then we define our multivariate and symmetric-function extensions.

Recall first that the {\em Bell number}\/ $B_n$
is the number of partitions of an $n$-element set into nonempty blocks;
by convention $B_0 = 1$.
Refining this, the {\em Stirling subset number}\/
(also called {\em Stirling number of the second kind}\/)
$\stirlingsubset{n}{k}$ is
the number of partitions of an $n$-element set into $k$ nonempty blocks;
by convention $\stirlingsubset{0}{k} = \delta_{k0}$.
The {\em Bell polynomials}\/ $B_n(x)$ are then defined as
$B_n(x) = \sum_{k=0}^n \stirlingsubset{n}{k} x^k$.\footnote{
   See \cite[A008277/A048993]{OEIS} for further information
   on the Stirling subset numbers and Bell polynomials.
}

Similarly, the {\em Lah number}\/ $L_n$ 
is the number of partitions of an $n$-element set
into nonempty linearly ordered blocks (also called {\em lists}\/);
we set $L_0 = 1$.
Refining this, the {\em Lah number}\/ $L(n,k)$ is
the number of partitions of an $n$-element set
into $k$ nonempty linearly ordered blocks;
we set $L(0,k) = \delta_{k0}$.
The Lah numbers also have the explicit expression
\be
   L(n,k)  \;=\;  {n! \over k!} \binom{n-1}{n-k}
           \;=\;  \begin{cases}
                      \delta_{k0}       & \textrm{if $n=0$} \\[2mm]
                      {n! \over k!}  \binom{n-1}{k-1}
                                        & \textrm{if $n \ge 1$}
                  \end{cases}
\ee
The {\em Lah polynomials}\/ $L_n(x)$ are then defined as
$L_n(x) = \sum_{k=0}^n L(n,k) \: x^k$.\footnote{
   See \cite[A008297/A105278/A271703/A111596/A066667]{OEIS}
   for further information on the Lah numbers and Lah polynomials.
}

More generally, let $x$ and $\bfw = \{w_m\}_{m \ge 1}$ be indeterminates,
and let $P_n(x,\bfw)$ be the generating polynomial
for partitions of an $n$-element set into nonempty blocks
in~which each block of cardinality $m$ gets a weight $x w_m$:
\be
   P_n(x,\bfw)
   \;\eqdef\;
   \sum_{\pi \in \Pi_n} x^{|\pi|} \prod_{B \in \pi} w_{|B|}
   \;.
 \label{def.Pnxw}
\ee
(In particular, the empty set has a unique partition into nonempty blocks
--- namely, the partition with zero blocks --- so that $P_0(x,\bfw) = 1$.)
Then the Bell polynomials correspond to $w_m = 1$,
while the Lah polynomials correspond to $w_m = m!$.
It is not difficult to show that the polynomials $P_n(x,\bfw)$
have the exponential generating function
\be
   \sum_{n=0}^\infty P_n(x,\bfw) \, {t^n \over n!}
   \;=\;
   \exp \biggl( x \sum_{m=1}^\infty w_m \, {t^m \over m!} \biggr)
   \;.
\ee
The polynomials $P_n(x,\bfw)$ are also known
\cite[pp.~133--134]{Comtet_74}
as the
{\em complete Bell polynomials}\/ ${\bf Y}_n(xw_1,\ldots,xw_n)$.

Let us now express the Bell and Lah polynomials
in terms of a different combinatorial object,
namely, unordered forests of increasing ordered trees.
Recall first \cite[pp.~294--295]{Stanley_86}
that an {\em ordered tree}\/ (also called {\em plane tree}\/)
is a rooted tree in which the children of each vertex are linearly ordered.
An {\em unordered forest of ordered trees}\/
is an unordered collection of ordered trees.
An {\em increasing ordered tree}\/ is an ordered tree
in which the vertices carry distinct labels from a linearly ordered set
(usually some set of integers) in such a way that
the label of each child is greater than the label of its parent;
otherwise put, the labels increase along every path downwards from the root.
An {\em unordered forest of increasing ordered trees}\/
is an unordered forest of ordered trees with the same type of labeling.

Now let $\bphi = (\phi_i)_{i \ge 0}$ be indeterminates,
and let $L_{n,k}(\bphi)$ be the generating polynomial for
unordered forests of increasing ordered trees on the vertex set $[n]$,
having $k$ components (i.e.\ $k$ trees),
in which each vertex with $i$ children gets a weight $\phi_i$.
Clearly $L_{n,k}(\bphi)$ is a homogeneous polynomial of degree $n$
with nonnegative integer coefficients;
it is also quasi-homogeneous of degree $n-k$
when $\phi_i$ is assigned weight~$i$.
The first few polynomials $L_{n,k}(\bphi)$
[specialized for simplicity to $\phi_0 = 1$] are
%
%
%
\vspace*{2mm}
\begin{table}[H]
\centering
\footnotesize
\begin{tabular}{c|cccccc}
$n \setminus k$ & 0 & 1 & 2 & 3 & 4 & 5 \\
\hline
0 & 1 &  &  &  &  &  \\
1 & 0 & 1 &  &  &  &  \\
2 & 0 & $\phi_1$ & 1 &  &  &  \\
3 & 0 & $\phi_1^2 + 2 \phi_2$ & $3 \phi_1$ & 1 &  &  \\
4 & 0 & $\phi_1^3 + 8 \phi_1 \phi_2 + 6 \phi_3$ & $7 \phi_1^2 + 8 \phi_2$ & $6 \phi_1$ & 1 &  \\
5 & 0 & $\phi_1^4 + 22 \phi_1^2 \phi_2 + 16 \phi_2^2 + 42 \phi_1 \phi_3 + 24 \phi_4$ & $15 \phi_1^3 + 60 \phi_1 \phi_2 + 30 \phi_3$ & $25 \phi_1^2 + 20 \phi_2$ & $10 \phi_1$ & 1 \\
\end{tabular}
\end{table}
\noindent
(see also the Appendix for $n \le 7$).
Now let $y$ be an additional indeterminate,
and define the row-generating polynomials
$L_n(\bphi,y) = \sum_{k=0}^n L_{n,k}(\bphi) \: y^k$.
Then $L_n(\bphi,y)$ is quasi-homogeneous of degree $n$
when $\phi_i$ is assigned weight~$i$ and $y$ is assigned weight~1.
We call $L_{n,k}(\bphi)$ and $L_n(\bphi,y)$
the \textbfit{generic Lah polynomials},
and we call the lower-triangular matrix $\sfL = (L_{n,k}(\bphi))_{n,k \ge 0}$
the \textbfit{generic Lah triangle}.
Here $\bphi = (\phi_i)_{i \ge 0}$ are in the first instance indeterminates,
so that $L_{n,k}(\bphi) \in \Z[\bphi]$ and $L_n(\bphi,y) \in \Z[\bphi,y]$; 
but we can then, if we wish, substitute specific values for $\bphi$
in any commutative ring $R$,
leading to values $L_{n,k}(\bphi) \in R$ and $L_n(\bphi,y) \in R[y]$.
When doing this we use the same notation $L_{n,k}(\bphi)$ and $L_n(\bphi,y)$,
as the desired interpretation for $\bphi$ should be clear from the context.

Note, finally, that an
unordered forest of increasing ordered trees on the vertex set $[n]$,
with $k$~components, can be obtained by first choosing a partition of $[n]$
into $k$~nonempty blocks, and then constructing an increasing ordered tree
on each block.  It follows that the generic Lah polynomial $L_n(\bphi,y)$
equals the set-partition polynomial $P_n(y,\bfw)$
evaluated at $w_m = L_{m,1}(\bphi)$.

Now let $X = (x_i)_{i \ge 1}$ be indeterminates,
and let $\bee = (e_n(X))_{n \ge 0}$ and $\bh = (h_n(X))_{n \ge 0}$
be the elementary symmetric functions
and complete homogeneous symmetric functions, respectively;
they are elements of the ring $\Z[[X]]_{\rm sym}$ of symmetric functions,
considered as a subring of the formal-power-series ring $\Z[[X]]$.
We then define the \textbfit{Lah symmetric functions of positive type} by
$L_{n,k}^{(\infty)+}(X) = L_{n,k}(\bee)$ and
$L_n^{(\infty)+}(X,y) = L_n(\bee,y)$,
and the \textbfit{Lah symmetric functions of negative type} by
$L_{n,k}^{(\infty)-}(X) = L_{n,k}(\bh)$ and
$L_n^{(\infty)-}(X,y) = L_n(\bh,y)$.
Also, for any integer $r \ge 1$ we can imagine specializing $X$
by setting $x_i = 0$ for $i > r$;
we then define the
\textbfit{multivariate Lah polynomials of positive type} by
$L_{n,k}^{(r)+}(x_1,\ldots,x_r) = L_{n,k}(\bee(x_1,\ldots,x_r))$ and 
$L_n^{(r)+}(x_1,\ldots,x_r;y) = L_n(\bee(x_1,\ldots,x_r),y)$,
and the \textbfit{multivariate Lah polynomials of negative type} by
$L_{n,k}^{(r)-}(x_1,\ldots,x_r) = L_{n,k}(\bh(x_1,\ldots,x_r))$ and 
$L_n^{(r)-}(x_1,\ldots,x_r;y) = L_n(\bh(x_1,\ldots,x_r),y)$.\footnote{
   In \cite[section~12]{latpath_SRTR} we considered these quantities
   only for $n,k \ge 1$, and we used the notations
   $\scrq_{n,k}^{(\infty)}$, $\scrq_{n,k}^{(r)}$,
   $\scrq_{n,k}^{(\infty)-}$, $\scrq_{n,k}^{(r)-}$
   (with $n,k \ge 0$)
   for what we are now calling
   $L_{n+1,k+1}^{(\infty)+}$, $L_{n+1,k+1}^{(r)+}$,
   $L_{n+1,k+1}^{(\infty)-}$, $L_{n+1,k+1}^{(r)-}$,
   respectively;
   we called these the {\em Eulerian symmetric functions}\/
   and {\em multivariate Eulerian polynomials}\/.
   We now think that it might be preferable to reserve the term ``Eulerian''
   for quantities associated to trees, and to use instead the term ``Lah''
   for quantities associated to forests.
}
In the Appendix we report the Lah symmetric functions
$L_{n,k}^{(\infty)+}$ and $L_{n,k}^{(\infty)-}$ for $n \le 7$
in terms of the monomial symmetric functions $m_\lambda$.

These multivariate Lah polynomials and symmetric functions
can also be interpreted as generating polynomials for increasing
$r$-ary and multi-$r$-ary trees and forests ($1 \le r \le \infty$).
Let us recall first \cite[p.~295]{Stanley_86}
the recursive definition of an {\em $r$-ary tree}\/ ($1 \le r < \infty$):
it is either empty or else consists of a root together with
an ordered list of $r$ subtrees, each of which is an $r$-ary tree
(which may be empty).  We draw an edge from each vertex to the root
of each of its nonempty subtrees;  an edge from a vertex to the root
of its $i$th subtree will be called an {\em $i$-edge}\/.
Similarly, we can define recursively an {\em $\infty$-ary tree}\/:
it is either empty or else consists of a root together with
an ordered list of subtrees indexed by the positive integers $\PP$,
each of which is an $\infty$-ary tree (which may be empty)
{\em and only finitely many of which are nonempty}\/;
we define $i$-edges as before.\footnote{
   Please note that such a graph is necessarily finite
   (as always, the recursion is carried out only finitely many times).
}
But we can now view $r$-ary trees from a slightly different point of view:
an $r$-ary (resp.\ $\infty$-ary) tree is simply an ordered tree
in which each edge carries a label $i \in [r]$ (resp.\ $i \in \PP$)
and the edges emanating outwards from each vertex consist, in order,
of zero or one edges labeled 1, then zero or one edges labeled 2,
and so forth; an edge with label $i$ will be called an $i$-edge.
Let us now consider the generating polynomial for
unordered forests of increasing $\infty$-ary trees on the vertex set $[n]$,
having $k$~components, in which each $i$-edge gets a weight $x_i$.
Since the choice of labels on the edges emanating outwards from a vertex $v$
can be made independently for each $v$,
this is equivalent to evaluating the generating polynomial $L_{n,k}(\bphi)$
at $\phi_i = e_i(X)$;
in other words, it is the Lah symmetric function of positive type
$L_{n,k}^{(\infty)+}(X) = L_{n,k}(\bee)$.
Similarly, $L_n^{(\infty)+}(X,y) = L_n(\bee,y)$
is the generating polynomial for
unordered forests of increasing $\infty$-ary trees on the vertex set $[n]$,
in which each $i$-edge gets a weight $x_i$
and each tree (or equivalently, each root) gets a weight $y$.
And if we set $x_i = 0$ for $i > r$ so as to obtain $r$-ary trees or forests,
we get the multivariate Lah polynomials
$L_{n,k}^{(r)+}(x_1,\ldots,x_r) = L_{n,k}(\bee(x_1,\ldots,x_r))$
and $L_n^{(r)+}(x_1,\ldots,x_r;y) = L_n(\bee(x_1,\ldots,x_r),y)$.

The multivariate Lah polynomials and Lah symmetric functions of negative type
can be interpreted in a similar way.
We begin by adopting the reinterpretation of $r$-ary and $\infty$-ary trees
as ordered trees with labeled edges,
and then consider \cite[section~10.3.2]{latpath_SRTR}
the variant in which the number of edges
of each label emanating from a given vertex,
instead of being ``zero or one'', is ``zero or more'':
we call this a \emph{multi-$r$-ary}
(resp.\ \emph{multi-$\infty$-ary}) \emph{tree}.
We now consider the generating polynomial for
unordered forests of increasing multi-$\infty$-ary trees
on the vertex set $[n]$,
having $k$~components, in which each $i$-edge gets a weight $x_i$.
This is equivalent to evaluating the generating polynomial $L_{n,k}(\bphi)$
at $\phi_i = h_i(X)$;
in other words, it is the Lah symmetric function of negative type
$L_{n,k}^{(\infty)-}(X) = L_{n,k}(\bh)$.
Similarly, $L_n^{(\infty)-}(X,y) = L_n(\bh,y)$
is the generating polynomial for
unordered forests of increasing multi-$\infty$-ary trees
on the vertex set $[n]$,
in which each $i$-edge gets a weight $x_i$
and each tree (or equivalently, each root) gets a weight $y$.
And if we set $x_i = 0$ for $i > r$ so as to obtain
multi-$r$-ary trees or forests,
we get the multivariate Lah polynomials
$L_{n,k}^{(r)-}(x_1,\ldots,x_r) = L_{n,k}(\bh(x_1,\ldots,x_r))$
and $L_n^{(r)-}(x_1,\ldots,x_r;y) = L_n(\bh(x_1,\ldots,x_r),y)$.

Let us now consider the further specialization
of the multivariate Lah polynomials of positive type
to $x_1 = \ldots = x_r = 1$, corresponding to $\phi_i = \binom{r}{i}$.
It is well known \cite[p.~24]{Stanley_86}
that the number of increasing binary trees on $n$ vertices is $n!$,
and more generally that the number of increasing $r$-ary trees
on $n$ vertices is the multifactorial $F_n^{(r-1)}$
\cite[p.~30, Example~1]{Bergeron_92},
where
\be
   F_n^{(r)}  \;\eqdef\;  \prod_{j=0}^{n-1} (1+jr)
   \;.
 \label{def.multifactorial}
\ee
Therefore, the univariate $r$th-order Lah polynomials of positive type,
$L_n^{(r)+}(\bone;y)$, coincide with the set-partition polynomials $P_n(y,\bfw)$
defined in \reff{def.Pnxw} when we set $w_m = F_m^{(r-1)}$.
In particular, for $r=1$ we have $w_m = 1$
and obtain the Bell polynomials $B_n(y)$;
for $r=2$ we have $w_m = m!$
and obtain the univariate Lah polynomials $L_n(y)$;
for $r=3,4,5$ we have $w_m = (2m-1)!!, (3m-2)!!!, (4m-3)!!!!$
and obtain the row-generating polynomials
of \cite[A035342, A035469, A049029]{OEIS}.

In a similar way, we can specialize
the multivariate Lah polynomials of negative type
to $x_1 = \ldots = x_r = 1$, corresponding to $\phi_i = \binom{r+i-1}{i}$.
It is known \cite[p.~30, Corollary~1(iv)]{Bergeron_92}
that the number of increasing multi-unary trees on $n$ vertices
is $(2n-3)!!$;
more generally,
it was observed in \cite[section~12.3]{latpath_SRTR}
that the number of increasing multi-$r$-ary trees
on $n$ vertices is the shifted multifactorial $\widetilde{F}_{n-1}^{(r+1)}$,
where
\be
   \widetilde{F}_n^{(r)}  \;\eqdef\; \prod\limits_{j=0}^{n-1} [(r-1)+jr]
   \;.
\ee
Therefore, the univariate $r$th-order Lah polynomials of negative type,
$L_n^{(r)-}(\bone;y)$, coincide with the set-partition polynomials $P_n(y,\bfw)$
defined in \reff{def.Pnxw} when we set $w_m = \widetilde{F}_{m-1}^{(r+1)}$.
In particular, for $r=1$ we have $w_m = (2m-3)!!$
and obtain a variant of the Bessel polynomials \cite[A001497]{OEIS}
(see also Example~\ref{exam.phi=1} below);
for $r=2,3$ we have $w_m = (3m-4)!!!, (4m-5)!!!!$
and obtain the row-generating polynomials of \cite[A004747, A000369]{OEIS}.

Let us now explain how all this relates to total positivity.
Recall first that a finite or infinite matrix of real numbers is called
{\em totally positive}\/ (TP) if all its minors are nonnegative,
and {\em totally positive of order~$r$} (TP${}_r$)
if all its minors of size $\le r$ are nonnegative.
Background information on totally positive matrices can be found
in \cite{Karlin_68,Gantmacher_02,Pinkus_10,Fallat_11};
they have application to many fields of pure and applied mathematics.
In particular, it is known
\cite[Th\'eor\`eme~9]{Gantmakher_37} \cite[section~4.6]{Pinkus_10}
that a Hankel matrix $(a_{i+j})_{i,j \ge 0}$
of real numbers is totally positive if and only if the underlying sequence
$(a_n)_{n \ge 0}$ is a Stieltjes moment sequence
(i.e.\ the moments of a positive measure on $[0,\infty)$).
And a Toeplitz matrix $(a_{i-j})_{i,j \ge 0}$ of real numbers
(where $a_n = 0$ for $n < 0$)
is totally positive if and only if its ordinary generating function
can be written as
\be
   \sum_{n=0}^\infty a_n t^n
   \;=\;
   C e^{\gamma t} t^m \prod_{i=1}^\infty {1 + \alpha_i t  \over  1 - \beta_i t}
 \label{eq.thm.aissen}
\ee
with $m \in \N$, $C,\gamma,\alpha_i,\beta_i \ge 0$,
$\sum \alpha_i < \infty$ and $\sum \beta_i < \infty$
\cite[Theorem~5.3, p.~412]{Karlin_68}.

But this is only the beginning of the story,
because we are here principally concerned,
not with sequences and matrices of real numbers,
but with sequences and matrices of polynomials
(with integer or real coefficients) in one or more indeterminates~$\bfx$:
they will typically be generating polynomials that enumerate
some combinatorial objects with respect to one or more statistics.
We equip the polynomial ring $\R[\bfx]$ with the coefficientwise
partial order:  that is, we say that $P$ is nonnegative
(and write $P \myge 0$)
in case $P$ is a polynomial with nonnegative coefficients.
We then say that a matrix with entries in $\R[\bfx]$ is
\textbfit{coefficientwise totally positive}
if all its minors are polynomials with nonnegative coefficients;
and analogously for coefficientwise total positivity of order~$r$.
We say that a sequence $\ba = (a_n)_{n \ge 0}$ with entries in $\R[\bfx]$
is \textbfit{coefficientwise Hankel-totally positive}
(resp.\ \textbfit{coefficientwise Toeplitz-totally positive})
if its associated infinite Hankel (resp.\ Toeplitz) matrix
is coefficientwise totally positive;
and likewise for the versions of order $r$.
Similar definitions apply to the formal-power-series ring $\R[[\bfx]]$.
Most generally, we can consider sequences and matrices
with values in an arbitrary partially ordered commutative ring
(a precise definition will be given in Section~\ref{subsec.totalpos.prelim});
total positivity, Hankel-total positivity and Toeplitz-total positivity
are then defined in the obvious way.

Let us also explain some partial orders on the ring of symmetric functions.
Let~$R$ be a commutative ring
and let $\bfX = (x_i)_{i \ge 1}$ be a countably infinite collection
of indeterminates.
Then let $R[[\bfX]]_{\rm sym}$ be the ring of symmetric functions
with coefficients in $R$
\cite[Chapter~1]{Macdonald_95} \cite[Chapter~7]{Stanley_99};
it is a subring of the formal-power-series ring $R[[\bfX]]$.
(It goes without saying that ``function'' is a misnomer;
 these are formal power series.)
Now let $R$ carry a partial order $\scrp$.
When the coefficientwise order on $R[[\bfX]]$
is restricted to  $R[[\bfX]]_{\rm sym}$,
it becomes the \textbfit{monomial order}:
an element $f \in R[[\bfX]]_{\rm sym}$ is monomial-nonnegative
if and only~if it can be written as a (finite) nonnegative linear combination
of monomial symmetric functions $m_\lambda(\bfX)$.
However, the ring of symmetric functions can also be equipped
with a stronger order, namely the \textbfit{Schur order}:
an element $f \in R[[\bfX]]_{\rm sym}$ is Schur-nonnegative
if and only~if it can be written as a (finite) nonnegative linear combination
of Schur functions $s_\lambda(\bfX)$.
This indeed defines a partial order compatible with the ring structure,
since any product of Schur functions
is a nonnegative linear combination of Schur functions
(Littlewood--Richardson coefficients \cite[Section~7.A1.3]{Stanley_99}).
And it is strictly stronger than the monomial order,
because every Schur function
is a nonnegative linear combination of monomial symmetric functions
(Kostka numbers \cite[eq.~(7.35), p.~311]{Stanley_99})
but not conversely (e.g.\ $m_2 = s_2 - s_{11}$).

We can now state our main result:

\begin{theorem}[Total positivity of Lah matrices and Lah polynomials]
   \label{thm1.1}
\hfill\break\noindent
Fix $1 \le r \le \infty$.
Let $R$ be a partially ordered commutative ring,
and let $\bphi = (\phi_i)_{i \ge 0}$ be a sequence in $R$
that is Toeplitz-totally positive of order $r$.
Then:
\begin{itemize}
   \item[(a)]  The lower-triangular matrix
       $\sfL(\bphi) = (L_{n,k}(\bphi))_{n,k \ge 0}$
       is totally positive of order $r$ in the ring $R$.
   \item[(b)]  The sequence $\bL(\bphi) = (L_n(\bphi,y))_{n \ge 0}$
       is Hankel-totally positive of order $r$ in the ring $R[y]$
       equipped with the coefficientwise order.
   \item[(c)]  The sequence $\bL^\triangle(\bphi) = (L_{n+1,1}(\bphi))_{n \ge 0}$
       is Hankel-totally positive of order $r$ in the ring $R$.
\end{itemize}
\end{theorem}

Specializing this to $\bphi = \bee(X)$ or $\bh(X)$
and using the Jacobi--Trudi identity
\cite[Theorem~7.16.1 and Corollary~7.16.2]{Stanley_99},
we obtain:

\begin{corollary}[Total positivity of Lah symmetric functions]
   \label{cor1.2}
\hfill\break\noindent
\vspace*{-6mm}
\begin{itemize}
   \item[(a)]  The unit-lower-triangular matrices
      $\sfL^{(\infty)+} = (L_{n,k}(\bee(X)))_{n,k \ge 0}$ and
      $\sfL^{(\infty)-} =$ \linebreak $(L_{n,k}(\bh(X)))_{n,k \ge 0}$
      are totally positive with respect to the Schur order
      on the ring of symmetric functions (with coefficients in $\Z$).
   \item[(b)]  The sequences
      $\bL^{(\infty)+} = (L_n(\bee(X),y))_{n \ge 0}$ and
      $\bL^{(\infty)-} = (L_n(\bh(X),y))_{n \ge 0}$
      are Hankel-totally positive with respect to the Schur order
      on the ring of symmetric functions (with coefficients in $\Z$)
      and the coefficientwise order on polynomials in $y$.
   \item[(c)]  The sequences
      $\bL^{(\infty)+\triangle} = (L_{n+1,1}(\bee(X)))_{n \ge 0}$ and
      $\bL^{(\infty)-\triangle} = (L_{n+1,1}(\bh(X)))_{n \ge 0}$
      are Hankel-totally positive with respect to the Schur order
      on the ring of symmetric functions (with coefficients in $\Z$).
\end{itemize}
\end{corollary}
 
Weakening this result from the Schur order to the monomial order,
and then further specializing by setting $x_i = 0$ for $i > r$,
we obtain:

\begin{corollary}[Total positivity of multivariate Lah polynomials]
\nopagebreak
   \label{cor1.3}
\hfill\break\noindent
\nopagebreak
\vspace*{-6mm}
\begin{itemize}
   \item[(a)]  The unit-lower-triangular matrices
      $\sfL^{(r)+} = (L_{n,k}(\bee(x_1,\ldots,x_r)))_{n,k \ge 0}$ and
      $\sfL^{(r)-} = (L_{n,k}(\bh(x_1,\ldots,x_r)))_{n,k \ge 0}$
      are totally positive with respect to the coefficientwise order
      on the polynomial ring $\Z[x_1,\ldots,x_r]$.
   \item[(b)]  The sequences
      $\bL^{(r)+} = (L_n(\bee(x_1,\ldots,x_r),y))_{n \ge 0}$ and
      $\bL^{(r)-} =$ \qquad\hfill\linebreak
         $(L_n(\bh(x_1,\ldots,x_r),y))_{n \ge 0}$
      are Hankel-totally positive with respect to the coefficientwise order
      on the polynomial ring $\Z[x_1,\ldots,x_r,y]$.
   \item[(c)]  The sequences
      $\bL^{(r)+\triangle} = (L_{n+1,1}(\bee(x_1,\ldots,x_r)))_{n \ge 0}$ and
      $\bL^{(r)-\triangle} =$ \qquad\hfill\linebreak 
         $(L_{n+1,1}(\bh(x_1,\ldots,x_r)))_{n \ge 0}$
      are Hankel-totally positive with respect to the coefficientwise order
      on the polynomial ring $\Z[x_1,\ldots,x_r]$.
\end{itemize}
\end{corollary}
 
{\bf Remarks.}
1. In Theorem~\ref{thm1.1} and its two corollaries,
part~(c) follows trivially from part~(b)
by dividing $L_{n+1}(\bphi,y)$ by $y$
and then specializing to $y=0$.
But in Section~\ref{sec.proofs} we will introduce a generalization
where the analogue of~(c) still holds (by a different proof),
but it is unknown whether there is any analogue of~(b).

2. Corollaries~\ref{cor1.2}(a) and \ref{cor1.3}(a)
are essentially already contained in
\cite[Corollaries~12.5 and 12.25 and Remark after the proof of Lemma~12.13]{latpath_SRTR}.
But the extension in Theorem~1.1(a) to the generic Lah polynomials is new.

3.  Similarly, Corollaries~\ref{cor1.2}(c) and \ref{cor1.3}(c)
are already contained in \cite[Corollaries~12.3 and 12.24]{latpath_SRTR}.
Indeed, a slightly stronger version of Corollaries~\ref{cor1.2}(c)
and \ref{cor1.3}(c) for $\sfL^{(\infty)+}$ and $\sfL^{(r)+}$
--- in which the Hankel-totally positive sequence is extended backwards
by prepending one element ---
is given in
\cite[Theorem~12.1(a) and Corollaries~12.3 and 12.7]{latpath_SRTR}.
There is also an analogous result for the case of negative type
\cite[Theorem~12.20(a) and Corollary~12.22]{latpath_SRTR},
but it does not seem to imply (at least in any obvious way)
the negative-type case of Corollaries~\ref{cor1.2}(c) and \ref{cor1.3}(c).
\myendremark

\medskip

Our proof of Theorem~\ref{thm1.1} will be based on the method of
{\em production matrices}\/ \cite{Deutsch_05,Deutsch_09}.
We shall review this theory
in Sections~\ref{subsec.production} and \ref{subsec.totalpos.prodmat},
so now we state only the bare-bones definitions.
Let $P = (p_{ij})_{i,j \ge 0}$ be an infinite matrix
with entries in a commutative ring~$R$;
we assume that $P$ is either row-finite
(i.e.\ has only finitely many nonzero entries in each row)
or column-finite.
Now define an infinite matrix $A = (a_{nk})_{n,k \ge 0}$ by
\be
   a_{nk}  \;=\;  (P^n)_{0k}
   \;.
\ee
We call $P$ the \textbfit{production matrix}
and $A$ the \textbfit{output matrix}, and we write $A = \scro(P)$.
The two key facts here are the following \cite{Sokal_totalpos}:
if $R$ is a partially ordered commutative ring
and $P$ is totally positive of order $r$,
then $\scro(P)$ is totally positive of order $r$
and the zeroth column of $\scro(P)$ is Hankel-totally positive of order $r$.
See Section~\ref{subsec.totalpos.prodmat} for precise statements and proofs.

We shall prove Theorem~\ref{thm1.1} by explicitly constructing
the production matrix that generates the
the generic Lah triangle $\sfL = (L_{n,k}(\bphi))_{n,k \ge 0}$,
and then verifying its total positivity.
Let $\Delta = (\delta_{i+1,j})_{i,j \ge 0}$
be the matrix with 1 on the superdiagonal and 0 elsewhere.
We then have:

\begin{proposition}[Production matrix for the generic Lah triangle]
   \label{prop.prodmat}
Let $\bphi = (\phi_i)_{i \ge 0}$ and $y$ be indeterminates,
and work in the ring $\Z[\bphi,y]$.
Define the lower-Hessenberg matrix $P = (p_{ij})_{i,j \ge 0}$ by
\be
   p_{ij}
   \;=\; 
   \begin{cases}
      0   &  \textrm{if $j=0$ or $j > i+1$}  \\[1mm]
      {\displaystyle {i! \over (j-1)!}} \, \phi_{i-j+1}
          &  \textrm{if $1 \le j \le i+1$}
   \end{cases}
 \label{eq.prop.prodmat}
\ee
and the unit-lower-triangular $y$-binomial matrix $B_y$ by
\be
   (B_y)_{nk}  \;=\;  \binom{n}{k} \, y^{n-k}
   \;.
\ee
Then:
\begin{itemize}
   \item[(a)]  $P$ is the production matrix for the generic Lah triangle
       $\sfL = (L_{n,k}(\bphi))_{n,k \ge 0}$.
   \item[(b)]  $B_y^{-1} P B_y = P (I + y \Delta^{\rm T})$
       is the production matrix for $\sfL B_y$.
\end{itemize}
\end{proposition}

We will prove Proposition~\ref{prop.prodmat}(a)
by constructing a bijection from
{\em ordered}\/ forests of increasing ordered trees
to labeled partial \L{}ukasiewicz paths,
along the lines of \cite[proofs of Theorems~12.11 and 12.28]{latpath_SRTR}.
Then Proposition~\ref{prop.prodmat}(b) will follow
by a straightforward but slightly nontrivial computation
(Lemma~\ref{lemma.ByinvPBy}).

In fact, we will prove a generalization of Proposition~\ref{prop.prodmat}(a)
[and hence also of Theorem~\ref{thm1.1}(a,c)
 and Corollaries~\ref{cor1.2}(a,c) and \ref{cor1.3}(a,c)]
for some polynomials $\Lhat_{n,k}(\bphihat)$,
to be defined in Section~\ref{subsec.proofs.prodmat},
that depend on a refined set of indeterminates
$\bphihat = (\phi_i^{[L]})_{i \ge 0,\, L \ge 1}$
and that reduce to $L_{n,k}(\bphi)$ when $\phi_i^{[L]} = \phi_i$ for all $L$.
However, no analogue of Proposition~\ref{prop.prodmat}(b)
[and hence also of Theorem~\ref{thm1.1}(b)
 and Corollaries~\ref{cor1.2}(b) and \ref{cor1.3}(b)]
appears to exist for these more general polynomials.

\medskip

{\bf Remark.}
If we were to work in the ring $\Q[\bphi,y]$ instead of $\Z[\bphi,y]$,
we would have $P = D T_\infty(\bphi) D^{-1} \Delta$
where $T_\infty(\bphi)$ is the infinite lower-triangular Toeplitz matrix
associated to the sequence $\bphi$,
and $D = \diag\bigl( (n!)_{n \ge 0} \bigr)$.
\myendremark

\medskip

Now return to the situation of Theorem~\ref{thm1.1}.
If the ring $R$ contains the rationals (with their usual order),
it follows from $P = D T_\infty(\bphi) D^{-1} \Delta$
that $P$ is totally positive of order $r$
whenever $\bphi$ is Toeplitz-totally positive of order $r$;
and the same holds for $B_y^{-1} P B_y = P (I + y \Delta^{\rm T})$.
And even if $R$ does not contain the rationals,
it~turns out that the same conclusions are true,
as we can show with a bit more work (Lemma~\ref{lemma.diagmult.TP}).
Theorem~\ref{thm1.1} is then an immediate consequence
of Proposition~\ref{prop.prodmat} and Lemma~\ref{lemma.diagmult.TP}
together with the general theory of production matrices and total positivity
(Section~\ref{subsec.totalpos.prodmat}).

Now fix an integer $r \ge 1$, and let us consider the
multivariate Lah polynomials of positive type
by specializing the production matrix \reff{eq.prop.prodmat}
to $\phi_n = e_n(x_1,\ldots,x_r)$.
Recall that the product of two lower-triangular Toeplitz matrices
corresponds to the convolution of their generating sequences,
or equivalently the product of their ordinary generating functions;
and since
$\sum\limits_{n=0}^\infty e_n(x_1,\ldots,x_r) \, t^n
 = \prod\limits_{i=1}^r (1 + x_i t)$,
it follows that the lower-triangular Toeplitz matrix
$T_\infty(\bee(x_1,\ldots,x_r))$
has the factorization $T_\infty(\bee) = L_1 \cdots L_r$
where $L_i = L(1,1,\ldots;x_i,x_i,\ldots)$
is the lower-bidiagonal Toeplitz matrix with 1 on the diagonal
and $x_i$ on the subdiagonal.
Therefore, $D T_\infty(\bee) D^{-1}$
has the factorization $L'_1 \cdots L'_r$
where $L'_i = D L_i D^{-1} = L(1,1,\ldots;x_i,2x_i,3x_i,\ldots)$;
the production matrix
$P(x_1,\ldots,x_r) = D T_\infty(\bee) D^{-1} \Delta$
has the factorization $L'_1 \cdots L'_r \Delta$;
and the modified production matrix
$B_y^{-1} P B_y = P (I + y \Delta^{\rm T})$ has the factorization
\be
   B_y^{-1} \, P(x_1,\ldots,x_r) \, B_y
   \;=\;
   L'_1 \,\cdots\, L'_r \, (\Delta + yI)
 \label{eq.lah.prodmat.y}
\ee
(since $\Delta \Delta^{\rm T} = I$).
On the other hand, \reff{eq.lah.prodmat.y}
is precisely the production matrix for an $r$-branched S-fraction
with coefficients
\be
   \balpha \;=\; (\alpha_i)_{i \ge r}
      \;=\; y,x_1,\ldots,x_r,y,2x_1,\ldots,2x_r,y,3x_1,\ldots,3x_r,\ldots
 \label{eq.lah.alphas}
\ee
(see~\cite[Propositions~7.2 and 8.2(b)]{latpath_SRTR}
 and eq.~\reff{eq.prop.contraction} below).
Since the zeroth column of the matrix $\sfL B_y$
is given by the Lah polynomials $L_n(\bphi,y)$,
it follows that the multivariate Lah polynomials of positive type
$L_n^{(r)+}(x_1,\ldots,x_r;y)$
are given by an $r$-branched S-fraction with coefficients \reff{eq.lah.alphas}:

\begin{theorem}[Branched S-fraction for multivariate Lah polynomials of positive type]
   \label{thm.lah.S-fraction}
We have $L_n^{(r)+}(x_1,\ldots,x_r;y) = S_n^{(r)}(\balpha)$
where the coefficients $\balpha$ are given by \reff{eq.lah.alphas}.
[Here $S_n^{(r)}(\balpha)$ is the $r$-Stieltjes--Rogers polynomial
 of order $n$;  see Section~\ref{subsec.mSR} for the precise definition.]
\end{theorem}

\medskip

{\bf Remarks.}
1.  For $r=1$ the multivariate Lah polynomials of positive type
are simply the homogenized Bell polynomials
$L_n^{(1)+}(x_1;y) = x_1^n B_n(y/x_1)$,
and this is the well-known classical S-fraction
with coefficients $\balpha = y,x_1,y,2x_1,y,3x_1,\ldots$.

2. Since $L_n^{(r)+}(x_1,\ldots,x_r;y)$ is invariant under
permutations of $x_1,\ldots,x_r$, it is actually represented by
$r!$ different $r$-branched S-fractions in which the coefficients $\balpha$
are obtained from \reff{eq.lah.alphas} by permuting $x_1,\ldots,x_r$.
This illustrates the nonuniqueness of $r$-branched S-fractions for $r \ge 2$
\cite{latpath_SRTR}.

3. In Section~\ref{sec.euler-gauss} we will also give a completely independent
proof of Theorem~\ref{thm.lah.S-fraction}, based on the Euler--Gauss
recurrence method.
\myendremark

\medskip

For the multivariate Lah polynomials of negative type,
the lower-triangular Toeplitz matrix $T_\infty(\bh(x_1,\ldots,x_r))$
has the factorization $T_\infty(\bh) = \Ltilde_1 \cdots \Ltilde_r$
where $\Ltilde_i =$ \linebreak $L(1,1,\ldots;-x_i,-x_i,\ldots)^{-1}$;
so $D T_\infty(\bh) D^{-1}$
has the factorization $\Ltilde'_1 \cdots \Ltilde'_r$
where
$\Ltilde'_i = D \Ltilde_i D^{-1} = L(1,1,\ldots;-x_i,-2x_i,-3x_i,\ldots)^{-1}$;
and the production matrix $D T_\infty(\bh) D^{-1} \Delta$
has the factorization $\Ltilde'_1 \cdots \Ltilde'_r \Delta$.
But the matrices $\Ltilde_i$ and $\Ltilde'_i$
are dense lower-triangular, not lower-bidiagonal,
so we do not see any way of interpreting this as the production matrix
of a branched S-fraction.
Indeed, we have verified
that the multivariate Lah {\em numbers}\/ of negative type
$L_n^{(r)-}(1,\ldots,1;1)$
{\em cannot}\/ be expressed as an $m$-branched S-fraction
of the following types:
\begin{itemize}
   \item For $r=1,2,3$, the numbers $L_n^{(r)-}(1,\ldots,1;1)$
      cannot be expressed as a 2-branched S-fraction
      with {\em nonnegative integer}\/ coefficients:
      this was verified by exhaustive computer search using
      $n \le 7,8,7$, respectively.
   \item For $m > r+1$, the numbers $a_n = L_n^{(r)-}(1,\ldots,1;1)$
      cannot be expressed as an $m$-branched S-fraction
      with {\em positive integer}\/ coefficients:
      this is simply because $a_0 = a_1 = 1$ and $a_2 = r+1$,
      while the $m$-Stieltjes--Rogers polynomial
      $S_2^{(m)}(\balpha)$ equals
      $\alpha_m (\alpha_m + \alpha_{m+1} + \cdots + \alpha_{2m-1})$.
\end{itemize}
For $(r,m) = (1,3)$ and $(2,3)$,
our computations (as far as we were able to go)
were unable to give either a proof of nonexistence
(with nonnegative integer coefficients)
or a comprehensible candidate for $\balpha$.

\bigskip

Although the present paper is a follow-up to our paper \cite{latpath_SRTR},
we have endeavored, for the convenience of the reader,
to make it as self-contained as possible.
We have therefore begun (Section~\ref{sec.prelim})
with a brief review of the key definitions and results from \cite{latpath_SRTR}
that will be needed in the sequel.
We then proceed as follows:
In Section~\ref{sec.proofs} --- which is the technical heart of the paper ---
we prove Proposition~\ref{prop.prodmat},
from which we deduce Theorem~\ref{thm1.1};
indeed, we state and prove a generalization involving
a refined set of indeterminates.
In Section~\ref{sec.differential}
we give expressions for the multivariate Lah polynomials
of positive and negative type
in terms of the action of certain first-order linear differential operators.
In Section~\ref{sec.euler-gauss} we give a second proof
of Theorem~\ref{thm.lah.S-fraction},
based on the differential operators and the Euler--Gauss recurrence method.
In Section~\ref{sec.decorated} we interpret
the multivariate Lah polynomials of positive type $L_n^{(r)+}(\bfx,y)$
as generating polynomials for partitions of the set $[n]$
in~which each block is ``decorated'' with an additional structure.
In Section~\ref{sec.exponential} we compute explicit expressions
for the generic Lah polynomials $L_{n,k}(\bphi)$
by using exponential generating functions.
In the Note Added (Section~\ref{sec.exp_riordan})
we give an alternate proof of Proposition~\ref{prop.prodmat}(a),
using the theory of exponential Riordan arrays.
In the Appendix we report the generic Lah polynomials
and Lah symmetric functions for $n \le 7$.

\section{Preliminaries}   \label{sec.prelim}

Here we review some definitions and results from \cite{latpath_SRTR}
that will be needed in the sequel.

\subsection{Partially ordered commutative rings and total positivity}
   \label{subsec.totalpos.prelim}

In this paper all rings will be assumed to have an identity element 1
and to be nontrivial ($1 \ne 0$).

A \textbfit{partially ordered commutative ring} is a pair $(R,\scrp)$ where
$R$ is a commutative ring and $\scrp$ is a subset of $R$ satisfying
\begin{itemize}
   \item[(a)]  $0,1 \in \scrp$.
   \item[(b)]  If $a,b \in \scrp$, then $a+b \in \scrp$ and $ab \in \scrp$.
   \item[(c)]  $\scrp \cap (-\scrp) = \{0\}$.
\end{itemize}
We call $\scrp$ the {\em nonnegative elements}\/ of $R$,
and we define a partial order on $R$ (compatible with the ring structure)
by writing $a \le b$ as a synonym for $b-a \in \scrp$.
Please note that, unlike the practice in real algebraic geometry
\cite{Brumfiel_79,Lam_84,Prestel_01,Marshall_08},
we do {\em not}\/ assume here that squares are nonnegative;
indeed, this property fails completely for our prototypical example,
the ring of polynomials with the coefficientwise order,
since $(1-x)^2 = 1-2x+x^2 \not\myge 0$.

Now let $(R,\scrp)$ be a partially ordered commutative ring
and let $\bfx = \{x_i\}_{i \in I}$ be a collection of indeterminates.
In the polynomial ring $R[\bfx]$ and the formal-power-series ring $R[[\bfx]]$,
let $\scrp[\bfx]$ and $\scrp[[\bfx]]$ be the subsets
consisting of polynomials (resp.\ series) with nonnegative coefficients.
Then $(R[\bfx],\scrp[\bfx])$ and $(R[[\bfx]],\scrp[[\bfx]])$
are partially ordered commutative rings;
we refer to this as the \textbfit{coefficientwise order}
on $R[\bfx]$ and $R[[\bfx]]$.

A (finite or infinite) matrix with entries in a
partially ordered commutative ring
is called \textbfit{totally positive} (TP) if all its minors are nonnegative;
it is called \textbfit{totally positive of order~$\bm{r}$} (TP${}_r$)
if all its minors of size $\le r$ are nonnegative.
It follows immediately from the Cauchy--Binet formula that
the product of two TP (resp.\ TP${}_r$) matrices is TP
(resp.\ TP${}_r$).\footnote{
   For infinite matrices, we need some condition to ensure that
   the product is well-defined.
   For instance, the product $AB$ is well-defined whenever
   $A$ is row-finite (i.e.\ has only finitely many nonzero entries in each row)
   or $B$ is column-finite.
}
This fact is so fundamental to the theory of total positivity
that we shall henceforth use it without comment.

We say that a sequence $\ba = (a_n)_{n \ge 0}$
with entries in a partially ordered commutative ring
is \textbfit{Hankel-totally positive} 
(resp.\ \textbfit{Hankel-totally positive of order~$\bm{r}$})
if its associated infinite Hankel matrix
$H_\infty(\ba) = (a_{i+j})_{i,j \ge 0}$
is TP (resp.\ TP${}_r$).
We say that $\ba$
is \textbfit{Toeplitz-totally positive} 
(resp.\ \textbfit{Toeplitz-totally positive of order~$\bm{r}$})
if its associated infinite Toeplitz matrix
$T_\infty(\ba) = (a_{i-j})_{i,j \ge 0}$
(where $a_n \eqdef 0$ for $n < 0$)
is TP (resp.\ TP${}_r$).\footnote{
   When $R = \R$, Toeplitz-totally positive sequences are traditionally called
   {\em P\'olya frequency sequences}\/ (PF),
   and Toeplitz-totally positive sequences of order $r$
   are called {\em P\'olya frequency sequences of order $r$}\/ (PF${}_r$).
   See \cite[chapter~8]{Karlin_68} for a detailed treatment.
}

We will need an easy fact about the total positivity of special matrices:

\begin{lemma}[Bidiagonal matrices]
  \label{lemma.bidiagonal}
Let $A$ be a matrix with entries in a partially ordered commutative ring,
with the property that all its nonzero entries belong to two consecutive
diagonals.
Then $A$ is totally positive if and only if all its entries are nonnegative.
\end{lemma}

\proof
The nonnegativity of the entries (i.e.\ TP${}_1$)
is obviously a necessary condition for TP.
Conversely, for a matrix of this type it is easy to see that
every nonzero minor is simply a product of some entries.
\qed

\subsection{Production matrices}   \label{subsec.production}

The method of production matrices \cite{Deutsch_05,Deutsch_09}
has become in recent years an important tool in enumerative combinatorics.
In the special case of a tridiagonal production matrix,
this construction goes back to Stieltjes' \cite{Stieltjes_1889,Stieltjes_1894}
work on continued fractions:
the production matrix of a classical S-fraction or J-fraction is tridiagonal.
In~the present paper, by contrast,
we shall need production matrices that are lower-Hessenberg
(i.e.\ vanish above the first superdiagonal)
but are not in general tridiagonal.
We therefore begin by reviewing briefly
the basic theory of production matrices.
The important connection of production matrices with total positivity
will be treated in the next subsection.

Let $P = (p_{ij})_{i,j \ge 0}$ be an infinite matrix
with entries in a commutative ring $R$.
In~order that powers of $P$ be well-defined,
we shall assume that $P$ is either row-finite
(i.e.\ has only finitely many nonzero entries in each row)
or column-finite.

Let us now define an infinite matrix $A = (a_{nk})_{n,k \ge 0}$ by
\be
   a_{nk}  \;=\;  (P^n)_{0k}
 \label{def.iteration}
\ee
(in particular, $a_{0k} = \delta_{0k}$).
Writing out the matrix multiplications explicitly, we have
\be
   a_{nk}
   \;=\;
   \sum_{i_1,\ldots,i_{n-1}}
      p_{0 i_1} \, p_{i_1 i_2} \, p_{i_2 i_3} \,\cdots\,
        p_{i_{n-2} i_{n-1}} \, p_{i_{n-1} k}
   \;,
 \label{def.iteration.walk}
\ee
so that $a_{nk}$ is the total weight for all $n$-step walks in $\N$
from $i_0 = 0$ to $i_n = k$, in~which the weight of a walk is the
product of the weights of its steps, and a step from $i$ to $j$
gets a weight $p_{ij}$.
Yet another equivalent formulation is to define the entries $a_{nk}$
by the recurrence
\be
   a_{nk}  \;=\;  \sum_{i=0}^\infty a_{n-1,i} \, p_{ik}
   \qquad\hbox{for $n \ge 1$}
 \label{def.iteration.bis}
\ee
with the initial condition $a_{0k} = \delta_{0k}$.

We call $P$ the \textbfit{production matrix}
and $A$ the \textbfit{output matrix},
and we write $A = \scro(P)$.
Note that if $P$ is row-finite, then so is $\scro(P)$;
if $P$ is lower-Hessenberg, then $\scro(P)$ is lower-triangular;
if $P$ is lower-Hessenberg with invertible superdiagonal entries,
then $\scro(P)$ is lower-triangular with invertible diagonal entries;
and if $P$ is unit-lower-Hessenberg
(i.e.\ lower-Hessenberg with entries 1 on the superdiagonal),
then $\scro(P)$ is unit-lower-triangular.
In all the applications in this paper, $P$ will be lower-Hessenberg.

The matrix $P$ can also be interpreted as the adjacency matrix
for a weighted directed graph on the vertex set $\N$
(where the edge $ij$ is omitted whenever $p_{ij}  = 0$).
Then $P$ is row-finite (resp.\ column-finite)
if and only if every vertex has finite out-degree (resp.\ finite in-degree).

This iteration process can be given a compact matrix formulation.
Let $\Delta = (\delta_{i+1,j})_{i,j \ge 0}$
be the matrix with 1 on the superdiagonal and 0 elsewhere.
Then for any matrix $M$ with rows indexed by $\N$,
the product $\Delta M$ is simply $M$ with its zeroth row removed
and all other rows shifted upwards by 1.
(Some authors use the notation $\overline{M} \eqdef \Delta M$.)
The recurrence \reff{def.iteration.bis} can then be written as
\be
   \Delta \, \scro(P)  \;=\;  \scro(P) \, P
   \;.
 \label{def.iteration.bis.matrixform}
\ee
It follows that if $A$ is a row-finite matrix
that has a row-finite inverse $A^{-1}$
and has first row $a_{0k} = \delta_{0k}$,
then $P = A^{-1} \Delta A$ is the unique matrix such that $A = \scro(P)$.
This holds, in particular, if $A$ is lower-triangular with
invertible diagonal entries and $a_{00} = 1$;
then $A^{-1}$ is lower-triangular
and $P = A^{-1} \Delta A$ is lower-Hessenberg.
And if $A$ is unit-lower-triangular,
then $P = A^{-1} \Delta A$ is unit-lower-Hessenberg.

We shall repeatedly use the following easy fact:

\begin{lemma}[Production matrix of a product]
   \label{lemma.production.AB}
Let $P = (p_{ij})_{i,j \ge 0}$ be a row-finite matrix
(with entries in a commutative ring $R$),
with output matrix $A = \scro(P)$;
and let $B = (b_{ij})_{i,j \ge 0}$
be a lower-triangular matrix with invertible (in $R$) diagonal entries.
Then
\be
   AB \;=\;  b_{00} \, \scro(B^{-1} P B)
   \;.
\ee
That is, up to a factor $b_{00}$,
the matrix $AB$ has production matrix $B^{-1} P B$.
\end{lemma}

\proof
Since $P$ is row-finite, so is $A = \scro(P)$;
then the matrix products $AB$ and $B^{-1} P B$
arising in the lemma are well-defined.  Now
\be
   a_{nk}
   \;=\;
   \sum_{i_1,\ldots,i_{n-1}}
      p_{0 i_1} \, p_{i_1 i_2} \, p_{i_2 i_3} \,\cdots\,
        p_{i_{n-2} i_{n-1}} \, p_{i_{n-1} k}
   \;,
\ee
while
\be
   \scro(B^{-1} P B)_{nk}
   \;=\;
   \sum_{j,i_1,\ldots,i_{n-1},i_n}
      (B^{-1})_{0j} \,
      p_{j i_1} \, p_{i_1 i_2} \, p_{i_2 i_3} \,\cdots\,
        p_{i_{n-2} i_{n-1}} \, p_{i_{n-1} i_n} \, b_{i_n k}
   \;.
\ee
But $B$ is lower-triangular with invertible diagonal entries,
so $(B^{-1})_{0j} = b_{00}^{-1} \delta_{j0}$.
It follows that $AB = b_{00} \, \scro(B^{-1} P B)$.
\qed

We will also need the following easy lemma:

\begin{lemma}[Production matrix of a down-shifted matrix]
   \label{lemma.down-shifted}
Let $P = (p_{ij})_{i,j \ge 0}$ be a row-finite or column-finite matrix
(with entries in a commutative ring $R$),
with output matrix $A = \scro(P)$;
and let $c$ be an element of $R$.
Now define
\be
   Q
   \;=\;
   \left[
   \begin{array}{c|c@{\hspace*{2mm}}c@{\hspace*{2mm}}c}
        0      &   c & 0 & \cdots \\
        \hline
        0      &     &   &        \\
        0      &     & P &        \\[-1mm]
        \vdots &     &   &        \\
   \end{array}
   \right]
   \;=\;
   c {\bf e}_{01} \,+\, \Delta^{\rm T} P \Delta
\ee
and
\be
   B
   \;=\;
   \left[
   \begin{array}{c|c@{\hspace*{2mm}}c@{\hspace*{2mm}}c}
        1      &   0 & 0 & \cdots \\
        \hline
        0      &     &   &        \\
        0      &     & cA &        \\[-1mm]
        \vdots &     &   &        \\
   \end{array}
   \right]
   \;=\;
   {\bf e}_{00} \,+\, c \Delta^{\rm T} A \Delta
   \;.
\ee
Then $B = \scro(Q)$.
\end{lemma}

\proof
We use \reff{def.iteration.walk} and its analogue for $Q$:
\be
   \scro(Q)_{nk}
   \;=\;
   \sum_{i_1,\ldots,i_{n-1}}
      q_{0 i_1} \, q_{i_1 i_2} \, q_{i_2 i_3} \,\cdots\,
        q_{i_{n-2} i_{n-1}} \, q_{i_{n-1} k}
   \;.
 \label{def.iteration.walk.BQ}
\ee
In \reff{def.iteration.walk.BQ}, the only nonzero contributions come from
$i_1 = 1$, with $q_{01} = c$;
and then we must also have $i_2,i_3,\ldots \ge 1$ and $k \ge 1$,
with $q_{ij} = p_{i-1,j-1}$.
Hence $\scro(Q)_{nk} = c a_{n-1,k-1}$ for $n \ge 1$.
\qed


\subsection{Production matrices and total positivity}
   \label{subsec.totalpos.prodmat}

Let $P = (p_{ij})_{i,j \ge 0}$ be a matrix with entries in a
partially ordered commutative ring $R$.
We will use $P$ as a production matrix;
let $A = \scro(P)$ be the corresponding output matrix.
As before, we assume that $P$ is either row-finite or column-finite.

When $P$ is totally positive, it turns out \cite{Sokal_totalpos}
that the output matrix $\scro(P)$ has {\em two}\/ total-positivity properties:
firstly, it is totally positive;
and secondly, its zeroth column is Hankel-totally positive.
Since \cite{Sokal_totalpos} is not yet publicly available,
we shall present briefly here (with proof) the main results
that will be needed in the sequel.

The fundamental fact that drives the whole theory is the following:

\begin{proposition}[Minors of the output matrix]
   \label{prop.iteration.homo}
Every $k \times k$ minor of the output matrix $A = \scro(P)$
can be written as a sum of products of minors of size $\le k$
of the production matrix $P$.
\end{proposition}

In this proposition the matrix elements $\bfp = \{p_{ij}\}_{i,j \ge 0}$
should be interpreted in the first instance as indeterminates:
for instance, we can fix a row-finite or column-finite set
$S \subseteq \N \times \N$
and define the matrix $P^S = (p^S_{ij})_{i,j \in \N}$ with entries
\be
   p^S_{ij}
   \;=\;
   \begin{cases}
       p_{ij}  & \textrm{if $(i,j) \in S$} \\[1mm]
       0       & \textrm{if $(i,j) \notin S$}
   \end{cases}
\ee
Then the entries (and hence also the minors) of both $P$ and $A$
belong to the polynomial ring $\Z[\bfp]$,
and the assertion of Proposition~\ref{prop.iteration.homo} makes sense.
Of course, we can subsequently specialize the indeterminates $\bfp$
to values in any commutative ring $R$.

\proofof{Proposition~\ref{prop.iteration.homo}}
Consider any minor of $A$ involving only the rows 0 through $N$.
We will prove the assertion of the Proposition by induction on $N$.
The statement is obvious for $N=0$.
For $N \ge 1$, let $A_N$ be the matrix consisting of rows
0 through $N-1$ of $A$, and let $A'_N$ be the matrix consisting of rows
1 through $N$ of $A$.  Then we have
\be
   A'_N  \;=\;  A_N P
   \;.
 \label{eq.proof.prop.iteration}
\ee
If the minor in question does not involve row 0,
then obviously it involves only rows 1 through $N$.
If the minor in question does involve row 0,
then it is either zero (in case it does not involve column~0)
or else equal to a minor of $A$ (of one size smaller)
that involves only rows 1 through $N$
(since $a_{0k} = \delta_{0k}$).
Either way it is a minor of $A'_N$;
but by \reff{eq.proof.prop.iteration} and the Cauchy--Binet formula,
every minor of $A'_N$ is a sum of products of minors (of the same size)
of $A_N$ and $P$.
This completes the inductive step.
\qed

If we now specialize the indeterminates $\bfp$
to values in some partially ordered commutative ring $R$,
we can immediately conclude:

\begin{theorem}[Total positivity of the output matrix]
   \label{thm.iteration.homo}
Let $P$ be an infinite matrix that is either row-finite or column-finite,
with entries in a partially ordered commutative ring $R$.
If $P$ is totally positive of order~$r$, then so is $A = \scro(P)$.
\end{theorem}

\medskip

{\bf Remarks.}
1.  In the case $R = \R$, Theorem~\ref{thm.iteration.homo}
is due to Karlin \cite[pp.~132--134]{Karlin_68};
see also \cite[Theorem~1.11]{Pinkus_10}.
Karlin's proof is different from ours.

2.  Our quick inductive proof of Proposition~\ref{prop.iteration.homo}
follows an idea of Zhu \cite[proof of Theorem~2.1]{Zhu_13},
which was in turn inspired in part by Aigner \cite[pp.~45--46]{Aigner_99}.
The same idea recurs in recent work of several authors
\cite[Theorem~2.1]{Zhu_14}
\cite[Theorem~2.1(i)]{Chen_15a}
\cite[Theorem~2.3(i)]{Chen_15b}
\cite[Theorem~2.1]{Liang_16}
\cite[Theorems~2.1 and 2.3]{Chen_19}
\cite{Gao_non-triangular_transforms}.
However, all of these results concerned only special cases:
\cite{Aigner_99,Zhu_13,Chen_15b,Liang_16}
treated the case in which the production matrix $P$ is tridiagonal;
\cite{Zhu_14} treated a (special) case in which $P$ is upper bidiagonal;
\cite{Chen_15a} treated the case in which
$P$ is the production matrix of a Riordan array;
\cite{Chen_19,Gao_non-triangular_transforms}
treated (implicitly) the case in which $P$ is upper-triangular and Toeplitz.
But the argument is in fact completely general, as we have just seen;
there is no need to assume any special form for the matrix $P$.
\myendremark

\bigskip

Now define 
$\scroo_0(P)$ to be the zeroth-column sequence of $\scro(P)$, i.e.
\be
   \scroo_0(P)_n  \;\eqdef\;  \scro(P)_{n0}  \;\eqdef\;  (P^n)_{00}
   \;.
 \label{def.scroo0}
\ee
Then the Hankel matrix of $\scroo_0(P)$ has matrix elements
\begin{eqnarray}
   & &
   \!\!\!\!\!\!\!
   H_\infty(\scroo_0(P))_{nn'}
   \;=\;
   \scroo_0(P)_{n+n'}
   \;=\;
   (P^{n+n'})_{00}
   \;=\;
   \sum_{k=0}^\infty (P^n)_{0k} \, (P^{n'})_{k0}
   \;=\;
          \nonumber \\
   & &
   \sum_{k=0}^\infty (P^n)_{0k} \, ((P^{\rm T})^{n'})_{0k}
   \;=\;
   \sum_{k=0}^\infty \scro(P)_{nk} \, \scro(P^{\rm T})_{n'k}
   \;=\;
   \big[ \scro(P) \, {\scro(P^{\rm T})}^{\rm T} \big]_{nn'}
   \;.
   \qquad
\end{eqnarray}
(Note that the sum over $k$ has only finitely many nonzero terms:
 if $P$ is row-finite, then there are finitely many nonzero $(P^n)_{0k}$,
 while if $P$ is column-finite,
 there are finitely many nonzero $(P^{n'})_{k0}$.)
We have therefore proven:

\begin{lemma}[Identity for Hankel matrix of the zeroth column]
   \label{lemma.hankel.karlin}
Let $P$ be a row-finite or column-finite matrix
with entries in a commutative ring $R$.
Then
\be
   H_\infty(\scroo_0(P))
   \;=\;
   \scro(P) \, {\scro(P^{\rm T})}^{\rm T}
   \;.
\ee
\end{lemma}

{\bf Remark.}
If $P$ is row-finite, then $\scro(P)$ is row-finite;
$\scro(P^{\rm T})$ need not be row- or column-finite,
but the product $\scro(P) \, {\scro(P^{\rm T})}^{\rm T}$
is anyway well-defined.
\myendremark

\medskip

Combining Proposition~\ref{prop.iteration.homo}
with Lemma~\ref{lemma.hankel.karlin} and the Cauchy--Binet formula,
we obtain:

\begin{corollary}[Hankel minors of the zeroth column]
   \label{cor.iteration2}
Every $k \times k$ minor of the infinite Hankel matrix
$H_\infty(\scroo_0(P)) = ((P^{n+n'})_{00})_{n,n' \ge 0}$
can be written as a sum of products
of the minors of size $\le k$ of the production matrix $P$.
\end{corollary}

And specializing the indeterminates $\bfp$
to nonnegative elements in a partially ordered commutative ring,
in such a way that $P$ is row-finite or column-finite,
we deduce:

\begin{theorem}[Hankel-total positivity of the zeroth column]
   \label{thm.iteration2bis}
Let $P = (p_{ij})_{i,j \ge 0}$ be an infinite row-finite or column-finite
matrix with entries in a partially ordered commutative ring $R$,
and define the infinite Hankel matrix
$H_\infty(\scroo_0(P)) = ((P^{n+n'})_{00})_{n,n' \ge 0}$.
If $P$ is totally positive of order~$r$, then so is $H_\infty(\scroo_0(P))$.
\end{theorem}

\subsection[$m$-Stieltjes--Rogers polynomials]{$\bm{m}$-Stieltjes--Rogers
   polynomials}  \label{subsec.mSR}

Throughout this subsection we fix an integer $m \ge 1$.
We recall \cite{Aval_08,Cameron_16,Prodinger_16,latpath_SRTR}
that an \textbfit{$\bm{m}$-Dyck path}
is a path in the upper half-plane $\Z \times \N$,
starting and ending on the horizontal axis,
using steps $(1,1)$ [``rise'' or ``up step'']
and $(1,-m)$ [``$m$-fall'' or ``down step''].
More generally, an \textbfit{$\bm{m}$-Dyck path at level $\bm{k}$}
is a path in $\Z \times \N_{\ge k}$,
starting and ending at height $k$,
using steps $(1,1)$ and $(1,-m)$.
Since the number of up steps must equal $m$ times the number of down steps,
the length of an $m$-Dyck path must be a multiple of $m+1$.

Now let $\balpha = (\alpha_i)_{i \ge m}$ be an infinite set of indeterminates.
Then \cite{latpath_SRTR}
the \textbfit{$\bm{m}$-Stieltjes--Rogers polynomial} of order~$n$,
denoted $S^{(m)}_n(\balpha)$, is the generating polynomial
for $m$-Dyck paths of length~$(m+1)n$ in which each rise gets weight~1
and each $m$-fall from height~$i$ gets weight $\alpha_i$.
Clearly $S_n^{(m)}(\balpha)$ is a homogeneous polynomial
of degree~$n$ with nonnegative integer coefficients.

Let $f_0(t) = \sum_{n=0}^\infty S^{(m)}_n(\balpha) \, t^n$
be the ordinary generating function for $m$-Dyck paths with these weights;
and more generally, let $f_k(t)$ be the ordinary generating function
for $m$-Dyck paths at level $k$ with these same weights.
(Obviously $f_k$ is just $f_0$ with each $\alpha_i$ replaced by $\alpha_{i+k}$;
 but we shall not explicitly use this fact.)
Then straightforward combinatorial arguments \cite{latpath_SRTR}
lead to the functional equation
\be
   f_k(t)  \;=\;  1 \:+\: \alpha_{k+m} t \, f_k(t) \, f_{k+1}(t) \,\cdots\, f_{k+m}(t)
 \label{eq.mSRfk.1}
\ee
or equivalently
\be
   f_k(t)  \;=\;  {1 \over 1 \:-\: \alpha_{k+m} t \, f_{k+1}(t) \,\cdots\, f_{k+m}(t)}
   \;.
 \label{eq.mSRfk.2}
\ee
Iterating \reff{eq.mSRfk.2}, we see immediately that $f_k$
is given by the branched continued fraction
\begin{subeqnarray}
   f_k(t)
   & = &
   \cfrac{1}
         {1 \,-\, \alpha_{k+m} t
            \prod\limits_{i_1=1}^{m}
                 \cfrac{1}
            {1 \,-\, \alpha_{k+m+i_1} t
               \prod\limits_{i_2=1}^{m}
               \cfrac{1}
            {1 \,-\, \alpha_{k+m+i_1+i_2} t
               \prod\limits_{i_3=1}^{m}
               \cfrac{1}{1 - \cdots}
            }
           }
         }
%
      \slabel{eq.fk.mSfrac.a} \\[2mm]
   & = &
\Scale[0.6]{
   \cfrac{1}{1 - \cfrac{\alpha_{k+m} t}{
     \Biggl( 1 - \cfrac{\alpha_{k+m+1} t}{
        \Bigl( 1  - \cfrac{\alpha_{k+m+2} t}{(\cdots) \,\cdots\, (\cdots)} \Bigr)
        \,\cdots\,
        \Bigl( 1  - \cfrac{\alpha_{k+2m+1} t}{(\cdots) \,\cdots\, (\cdots)} \Bigr)
       }
     \Biggr)
     \,\cdots\,
     \Biggl( 1 - \cfrac{\alpha_{k+2m} t}{
        \Bigl( 1  - \cfrac{\alpha_{k+2m+1} t}{(\cdots) \,\cdots\, (\cdots)} \Bigr)
        \,\cdots\,
        \Bigl( 1  - \cfrac{\alpha_{k+3m} t}{(\cdots) \,\cdots\, (\cdots)} \Bigr)
       }
     \Biggr)
    }
   }
}
     \nonumber \\
 \slabel{eq.fk.mSfrac.b}
 \label{eq.fk.mSfrac}
\end{subeqnarray}
and in particular that $f_0$ is given by
the specialization of \reff{eq.fk.mSfrac} to $k=0$.
We shall call the right-hand side of \reff{eq.fk.mSfrac}
an \textbfit{$\bm{m}$-branched Stieltjes-type continued fraction},
or \textbfit{$\bm{m}$-S-fraction} for short.

\medskip

{\bf Remark.}
In truth, we hardly ever use the branched continued fraction
\reff{eq.fk.mSfrac};
instead, we work directly with the $m$-Dyck paths
and/or with the recurrence \reff{eq.mSRfk.1}/\reff{eq.mSRfk.2}
that their generating functions satisfy.
\myendremark

\medskip

We now generalize these definitions as follows.
A \textbfit{partial $\bm{m}$-Dyck path}
is a path in the upper half-plane $\Z \times \N$,
starting on the horizontal axis but ending anywhere,
using steps $(1,1)$ [``rise'']
and $(1,-m)$ [``$m$-fall''].
A partial $m$-Dyck path starting at $(0,0)$
must stay always within the set
$V_m = \{ (x,y) \in \Z \times \N \colon\: x=y \bmod m+1 \}$.

Now let $\balpha = (\alpha_i)_{i \ge m}$ be an infinite set of indeterminates,
and let $S^{(m)}_{n,k}(\balpha)$ be the generating polynomial
for partial $m$-Dyck paths from $(0,0)$ to ${((m+1)n,(m+1)k)}$
in~which each rise gets weight~1
and each $m$-fall from height~$i$ gets weight $\alpha_i$.
We call the $S^{(m)}_{n,k}$ the
\textbfit{generalized $\bm{m}$-Stieltjes--Rogers polynomials}.
Obviously $S^{(m)}_{n,k}$ is nonvanishing only for $0 \le k \le n$,
and $S^{(m)}_{n,n} = 1$.
We therefore have an infinite unit-lower-triangular array
$\sfS^{(m)} = \big( S^{(m)}_{n,k}(\balpha) \big)_{n,k \ge 0}$
in which the first ($k=0$) column displays
the ordinary $m$-Stieltjes--Rogers polynomials $S^{(m)}_{n,0} = S^{(m)}_n$.

The production matrix for the triangle $\sfS^{(m)}$
was found in \cite[sections~7.1 and 8.2]{latpath_SRTR}.
We begin by defining some special matrices $M = (m_{ij})_{i,j \ge 0}$:
\begin{itemize}
   \item $L(s_1,s_2,\ldots)$ is the lower-bidiagonal matrix
       with 1 on the diagonal and $s_1,s_2,\ldots$ on the subdiagonal:
\be
   L(s_1,s_2,\ldots)
   \;=\;
   \begin{bmatrix}
      1  &     &     &     &    \\
      s_1 & 1  &     &     &    \\
          & s_2 & 1  &     &    \\
          &     & s_3 & 1  &    \\
          &     &     & \ddots & \ddots
   \end{bmatrix}
   \;.
 \label{def.L}
\ee
   \item $U^\star(s_1,s_2,\ldots)$ is the upper-bidiagonal matrix
       with 1 on the superdiagonal and $s_1,s_2,\ldots$ on the diagonal:
\be
   U^\star(s_1,s_2,\ldots)
   \;=\;
   \begin{bmatrix}
      s_1 & 1   &     &     &     &    \\
          & s_2 & 1   &     &     &    \\
          &     & s_3 & 1   &     &    \\
          &     &     & s_4 & 1   &    \\
          &     &     &     & \ddots & \ddots
   \end{bmatrix}
   \;.
 \label{def.Ustar}
\ee
\end{itemize}
Then the production matrix for the triangle $\sfS^{(m)}$ is
\begin{eqnarray}
   P^{(m)\mathrm{S}}(\balpha)
   & \eqdef &
   L(\alpha_{m+1}, \alpha_{2m+2}, \alpha_{3m+3}, \ldots)
   \:
   L(\alpha_{m+2}, \alpha_{2m+3}, \alpha_{3m+4}, \ldots)
   \:\cdots\:
   \hspace*{1cm}
       \nonumber \\
   & & \qquad
   L(\alpha_{2m}, \alpha_{3m+1}, \alpha_{4m+2}, \ldots)
   \:
   U^\star(\alpha_m, \alpha_{2m+1}, \alpha_{3m+2}, \ldots)
   \;,
   \hspace*{1cm}
 \label{eq.prop.contraction}
\end{eqnarray}
that is, the product of $m$ factors $L$ and one factor $U^\star$
\cite[Proposition~8.2]{latpath_SRTR}.

Finally, we proved the following fundamental results on total positivity
\cite[Theorems~9.8 and 9.10]{latpath_SRTR}:\footnote{
   In fact, we gave two independent proofs of these results:
   a graphical proof, based on the Lindstr\"om--Gessel--Viennot lemma;
   and an algebraic proof, based on the production matrix
   \reff{eq.prop.contraction}.
}
\begin{itemize}
   \item[(a)]  For each integer $m \ge 1$, the lower-triangular matrix
$\sfS^{(m)} = \big( S^{(m)}_{n,k}(\balpha) \big)_{n,k \ge 0}$
of generalized $m$-Stieltjes--Rogers polynomials
is totally positive in the polynomial ring $\Z[\balpha]$
equipped with the coefficientwise partial order.
   \item[(b)]  For each integer $m \ge 1$,
the sequence $\bS^{(m)} = ( S^{(m)}_n(\balpha) )_{n \ge 0}$
of $m$-Stieltjes--Rogers polynomials
is a Hankel-totally positive sequence in the polynomial ring $\Z[\balpha]$
equipped with the coefficientwise partial order.
\end{itemize}
Of course, we can then substitute for $\balpha$
any sequence of nonnegative elements
of any partially ordered commutative ring $R$,
and the resulting matrix $\sfS^{(m)}$ (resp.\ sequence $\bS^{(m)}$)
will be totally positive (resp.\ Hankel-totally positive) in $R$.


\section{Proofs of main results}
   \label{sec.proofs}

In this section we prove Theorem~\ref{thm1.1} and its corollaries,
by the following steps:
First we prove Proposition~\ref{prop.prodmat}(a),
which asserts that the matrix $P$ defined in \reff{eq.prop.prodmat}
is the production matrix for the generic Lah triangle
$\sfL = (L_{n,k}(\bphi))_{n,k \ge 0}$.
Next we prove the matrix identity $B_y^{-1} P B_y = P (I + y \Delta^{\rm T})$.
Once this is done, Lemma~\ref{lemma.production.AB} implies
that $P (I + y \Delta^{\rm T})$ is the production matrix for $\sfL B_y$,
which completes the proof of Proposition~\ref{prop.prodmat}(b).
Finally, we show that if the Toeplitz matrix $T_\infty(\bphi)$ is TP${}_r$,
then so is $P$
(Lemma~\ref{lemma.diagmult.TP} and Corollary~\ref{cor.diagmult.TP}).
Then Theorem~\ref{thm1.1} follows from
the general theory of production matrices and total positivity
(Theorems~\ref{thm.iteration.homo} and \ref{thm.iteration2bis}).

In fact, we will prove a generalization of Proposition~\ref{prop.prodmat}(a)
[and hence also of Theorem~\ref{thm1.1}(a,c)
 and Corollaries~\ref{cor1.2}(a,c) and \ref{cor1.3}(a,c)]
for some polynomials $\Lhat_{n,k}(\bphihat)$
that depend on a refined set of indeterminates
$\bphihat = (\phi_i^{[L]})_{i \ge 0,\, L \ge 1}$
and that reduce to $L_{n,k}(\bphi)$ when $\phi_i^{[L]} = \phi_i$ for all $L$.
However, no analogue of Proposition~\ref{prop.prodmat}(b)
[and hence also of Theorem~\ref{thm1.1}(b)
 and Corollaries~\ref{cor1.2}(b) and \ref{cor1.3}(b)]
appears to exist for these more general polynomials.

\subsection{A generalization of Proposition~\ref{prop.prodmat}(a)}
   \label{subsec.proofs.prodmat}

In this subsection, we shall state and prove a generalization
of Proposition~\ref{prop.prodmat}(a).
We begin by introducing the notion of {\em level}\/
of a vertex in a forest, as was done in \cite{latpath_SRTR}:

\begin{definition}[Level of a vertex]
   \label{def.level}
Let $F$ be a forest of increasing trees on a totally ordered vertex set,
with $k$ trees.\footnote{
   Here the forest can be either ordered or unordered;
   and the trees in the forest can be either ordered or unordered
   (in the sense that the children at each vertex
    can be either ordered or unordered).
   Neither of these orderings, if present,
   will play any role in the definition of ``level''.
}
For each vertex $j$ in $F$, let $r_j$ be the number of trees in $F$
that contain at least one vertex $\le j$.
Then the \textbfit{level} of the vertex~$j$ in the forest~$F$,
denoted $\lev_F(j)$,
is the number of children of vertices $< j$ whose labels are $>j$,
plus $k+1-r_j$. 
\end{definition}

\noindent
Note that $1 \le r_j \le k$, and hence $\lev_F(j) \ge 1$.

\bigskip

{\bf Remark.}
This definition of ``level'' is slightly different
from the one given in \cite{latpath_SRTR},
since our forests here have $k$ trees rather than $k+1$
as in \cite{latpath_SRTR},
and our levels here are $\ge 1$ rather than $\ge 0$.
\myendremark

\medskip

We can now define a generalization of our generic Lah triangle, as follows:
Let $\bphihat = (\phi_i^{[L]})_{i \ge 0,\, L \ge 1}$ be indeterminates,
and let $\Lhat_{n,k}(\bphihat)$ be the generating polynomial
for unordered forests of increasing ordered trees on the vertex set $[n]$,
having $k$ trees, in~which each vertex with $i$ children and level $L$
gets a weight $\phi_i^{[L]}$.
We shall refer to the lower-triangular matrix
$\sfLhat = (\Lhat_{n,k}(\bphihat))_{n,k \ge 0}$
as the \textbfit{refined generic Lah triangle}.
Of course, when $\phi_i^{[L]} = \phi_i$ for all $L$,
it reduces to the original generic Lah triangle
$\sfL = (L_{n,k}(\bphi))_{n,k \ge 0}$.

We shall see later that
the polynomial $\Lhat_{n,k}(\bphihat)$ has a factor
$\phi_0^{[1]} \phi_0^{[2]} \cdots \phi_0^{[k]}$.
So we can, if we wish, pull this factor out,
and consider also the lower-triangular array
$\sfLtilde = \big( \Lhat_{n,k}(\bphihat) /
                      (\phi_0^{[1]} \phi_0^{[2]} \cdots \phi_0^{[k]})
             \big)_{\! n,k \ge 0}$.

Finally, it turns out that, in proving the formula for the production matrix,
it~is most convenient to work with {\em ordered}\/ forests
of increasing ordered trees, not unordered ones.
Since the trees of our forests are labeled and hence distinguishable,
the generating polynomial for ordered forests
on the vertex set $[n]$ with $k$ components
is simply $k!$ times the generating polynomial for unordered forests.
So we will begin by finding the production matrix $P^{\rm ord}$
for the triangle
$\sfLhat^{\rm ord} = (k! \, \Lhat_{n,k}(\bphihat))_{n,k \geq 0}$;
then we will deduce from it the production matrix $P$
for the triangle $\sfLhat = (\Lhat_{n,k}(\bphihat))_{n,k \geq 0}$,
and the production matrix $\Ptilde$ for the triangle
$\sfLtilde = \big( \Lhat_{n,k}(\bphihat) /
                      (\phi_0^{[1]} \phi_0^{[2]} \cdots \phi_0^{[k]})
             \big)_{\! n,k \ge 0}$.

We now claim that the following generalization of
Proposition~\ref{prop.prodmat}(a) holds:

\begin{proposition}[Production matrix for the refined generic Lah triangle]
   \label{prop.prodmat.phiL}
\hfill\break\noindent
Let $\bphihat = (\phi_i^{[L]})_{i \ge 0,\, L \ge 1}$ be indeterminates,
and work in the ring $\mathbb{Z}[\bphihat]$.
Define the lower-Hessenberg matrices
$P^{\rm ord} = (p^{\rm ord}_{ij})_{i,j \ge 0}$,
$P = (p_{ij})_{i,j \ge 0}$ and $\Ptilde = (\ptilde_{ij})_{i,j \ge 0}$ by
\begin{eqnarray}
   p^{\rm ord}_{ij}
   & = &
   \begin{cases}
      0   &  \textrm{if $j=0$ or $j > i+1$}  \\[1mm]
      j \: \phi_{i-j+1}^{[j]}
          &  \textrm{if $1 \le j \le i+1$}
   \end{cases}
 \label{eq.prop.prodmat.phiL.0}
        \\[4mm]
   p_{ij}
   & = &
   \begin{cases}
      0   &  \textrm{if $j=0$ or $j > i+1$}  \\[1mm]
      {\displaystyle {i! \over (j-1)!}} \: \phi_{i-j+1}^{[j]}
          &  \textrm{if $1 \le j \le i+1$}
   \end{cases}
 \label{eq.prop.prodmat.phiL.1}
        \\[4mm]
   \ptilde_{ij}
   & = &
   \begin{cases}
      0   &  \textrm{if $j=0$ or $j > i+1$}  \\[1mm]
      {\displaystyle {i! \over (j-1)!}}
                \; \phi_0^{[j+1]} \cdots \phi_0^{[i]}
                \; \phi_{i-j+1}^{[j]}
          &  \textrm{if $1 \le j \le i$}   \\[4mm]
      1   &  \textrm{if $j=i+1$}
   \end{cases}
 \label{eq.prop.prodmat.phiL.2}
\end{eqnarray}
Then:
\begin{itemize}
   \item[(a)]  $P^{\rm ord}$ is the production matrix for the triangle
       $\sfLhat^{\rm ord} = (k! \, \Lhat_{n,k}(\bphihat))_{n,k \geq 0}$.
   \item[(b)]  $P$ is the production matrix for the
       refined generic Lah triangle
       $\sfLhat = (\Lhat_{n,k}(\bphihat))_{n,k \geq 0}$.
   \item[(c)]  $\Ptilde$ is the production matrix for the triangle
       $\sfLtilde = \big( \Lhat_{n,k}(\bphihat) /
                             (\phi_0^{[1]} \phi_0^{[2]} \cdots \phi_0^{[k]})
                    \big)_{\! n,k \ge 0}$.
\end{itemize}
\end{proposition}

As preparation for the proof of Proposition~\ref{prop.prodmat.phiL},
we recall the definition of the \textbfit{depth-first-search labeling}
of an ordered forest of ordered trees.
(The more precise name is {\em preorder traversal}\/,
 i.e.\ parent first, then children in order from left to right,
 carried out recursively starting at the root.)
The recursive definition can be found in \cite[pp.~33--34]{Stanley_99},
but there is a simple informal description:
for each tree, we walk clockwise around the tree,
starting at the root, and label the vertices in the order
in which they are first seen;
this is done successively for the trees of the forest, in the given order.
Note that, in the depth-first-search labeling,
all the children of a vertex~$j$ will have labels $> j$;
that is, the depth-first-search labeling is a (very special)
increasing labeling.
Note also that, in the depth-first-search labeling, if $r < r'$,
then all the vertices of the $r$th tree will have labels smaller than
all the vertices of the $r'$th tree;
of course, this property need {\em not}\/ hold
in a general increasing labeling.

Finally, we recall that a \textbfit{partial \L{}ukasiewicz path}
(in our definition)
is a path in the upper half-plane $\Z \times \N$
using steps $(1,s)$ with $-\infty < s \le 1$,
while a \textbfit{reversed partial \L{}ukasiewicz path}
is a path in the upper half-plane $\Z \times \N$
using steps $(1,s)$ with $-1 \le s < \infty$.

\proofof{Proposition~\ref{prop.prodmat.phiL}}
We will construct a bijection from the set of {\em ordered}\/ forests
of increasing ordered trees on the vertex set $[n]$, with $k$ components,
to a set of $\bfL$-labeled reversed partial \L{}ukasiewicz paths
from $(0,k)$ to $(n,0)$,
where the label sets $\bfL$ will be defined below.
The case $n=k=0$ is trivial, so we assume $n,k \ge 1$.

\begin{figure}[t]
\begin{center}
    \includegraphics[scale=1]{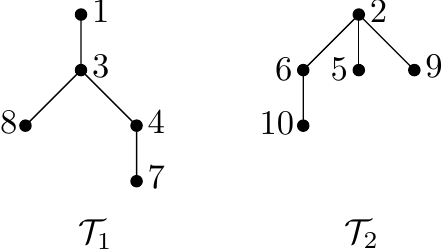}
    \hspace*{45pt}
    \includegraphics[scale=1]{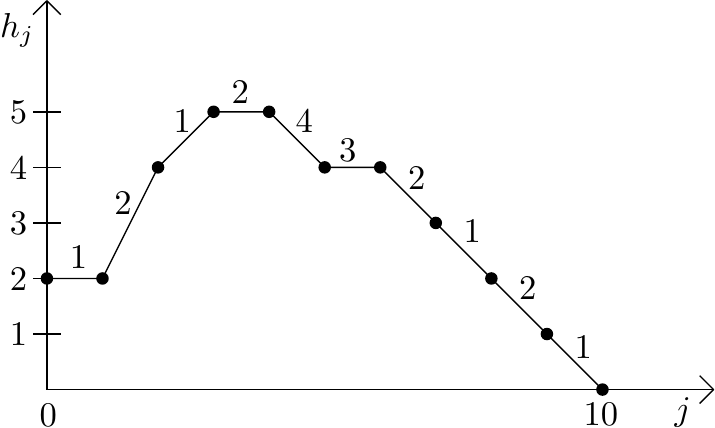}
\end{center}
\caption{An ordered forest of two increasing ordered trees
         on the vertex set $[10]$,
         and its image under the bijection.
         We put the label $\xi_j$ above the step $s_j$.
        }
\label{fig.bijection}
\end{figure}

Given an ordered forest $F$ of increasing ordered trees
on the vertex set $[n]$, with $k$ components,
we define a labeled reversed partial \L{}ukasiewicz path $(\omega,\xi)$
of length~$n$ as follows (see Figure~\ref{fig.bijection} for an example):

\bigskip

{\bf Definition of the path $\bomega$.}
The path $\omega$ starts at height $h_0 = k$
and takes steps $s_1,\ldots,s_n$ with $s_i = \deg(i) - 1$,
where $\deg(i)$ is the number of children of vertex $i$.
Therefore, the heights $h_0,\ldots,h_n$ are
\be
   h_j  \;=\;  k \,+\, \sum_{i=1}^j [\deg(i) - 1]
   \;.
\ee
Since $\sum_{i=1}^n \deg(i) = n-k$, we have $h_n = 0$.
We will show later that $h_1,\ldots,h_{n-1} \ge 1$.

\bigskip

{\bf Definition of the labels $\bxi$.}
The label $\xi_j$ is, by definition,
1 plus the number of vertices $> j$
that are either children of $\{1,\ldots,j-1\}$ or roots
and that precede $j$ in the depth-first-search order.\footnote{
   Here the depth-first-search order could be replaced by any chosen order
   on the vertices of $F$ that commutes with truncation.
   The key property we need is that
   the order on the truncated forest $F_{j-1}$ to be defined below
   is the restriction of the order on the full forest $F$.
}
Obviously $\xi_j$ is an integer $\ge 1$;
we will show later that $\xi_j \le h_{j-1}$ (Corollary~\ref{cor.xij}).

\bigskip

{\bf Interpretation of the heights $\bm{h_j}$.}
Recall that $r_j$ is the number of trees in $F$ that contain
at least one of the vertices $\{1,\ldots,j\}$.
We then claim:

\begin{lemma}
   \label{lemma.hj}
For $1 \le j \le n$,
the height $h_j = k+ \sum_{i=1}^j s_i$ has the following interpretations:
\begin{itemize}
   \item[(a)]  $h_j$ is the number of children of the vertices $\{1,\ldots,j\}$
      whose labels are $> j$, plus $k-r_j$.
   \item[(b)]  $h_{j-1}$ is the number of children of the vertices
      $\{1,\ldots,j-1\}$ whose labels are $> j$, plus $k+1-r_j$.
      That is, $h_{j-1}$ is the level of the vertex $j$
      as given in Definition~\ref{def.level}.
\end{itemize}
In particular, $h_j > k-r_j$ whenever $j$
is not the highest-numbered vertex of its tree,
and $h_j \ge k-r_j$ always.
\end{lemma}

\proof
By induction on $j$.
For the base case $j=1$, the claims are clear
since $r_1 = 1$, $h_0 = k$ and $h_1 = k + \deg(1) - 1$.

For $j>1$, the vertex $j$ is either the child of another node,
or the root of a tree.  We consider these two cases separately:

(i) Suppose that $j$ is the child of another node
(obviously numbered ${\le j-1}$).
By the inductive hypothesis~(a),
$h_{j-1}$ is the number of children of the vertices ${\{1,\ldots,j-1\}}$
whose labels are $\ge j$, plus $k-r_{j-1}$;
and since one of these children is $j$,
it follows that
$h_{j-1} - 1$ is the number of children of the vertices $\{1,\ldots,j-1\}$
whose labels are $> j$, plus $k-r_{j-1}$.
Now vertex $j$ has $\deg(j)$ children, all of which have labels $> j$;
so $h_j = h_{j-1} + s_j = h_{j-1} - 1 + \deg(j)$
is the number of children of the vertices $\{1,\ldots,j\}$
whose labels are $> j$, plus $k-r_{j-1}$.
Since $r_j = r_{j-1}$, the preceding two sentences prove
claims (b) and (a), respectively.

(ii) Suppose that $j$ is a root.
By the inductive hypothesis~(a),
$h_{j-1}$ is the number of children of the vertices $\{1,\ldots,j-1\}$
whose labels are $\ge j$, plus $k-r_{j-1}$;
and since $j$ is not one of these children,
it follows that
$h_{j-1}$ is also the number of children of the vertices $\{1,\ldots,j-1\}$
whose labels are $> j$, plus $k-r_{j-1}$.
Now vertex $j$ has $\deg(j)$ children, all of which have labels $> j$;
so $h_j = h_{j-1} +  \deg(j) - 1$
is the number of children of the vertices $\{1,\ldots,j\}$
whose labels are $> j$, plus $k-r_{j-1}-1$.
Since $r_j = r_{j-1} +1$, the preceding two sentences prove 
claims (b) and (a), respectively.
\qed

It follows from Lemma~\ref{lemma.hj}(b)
that $h_0,\ldots,h_{n-1} \ge 1$ and $h_n = 0$.
So the path $\omega$ is indeed a reversed partial \L{}ukasiewicz path
from $(0,k)$ to $(n,0)$,
which reaches level 0 only at the last step.

\begin{corollary}
   \label{cor.xij}
$\xi_j \le h_{j-1}$.
\end{corollary}

\proof
By Lemma~\ref{lemma.hj}(b), the number of vertices $> j$
that are children of ${ \{1,\ldots,j-1\} }$
is $h_{j-1} - (k+1 - r_j)$.
The number of vertices $> j$ that are roots is at most $k-r_j$
(since any tree containing a vertex $\le j$ necessarily has its root $\le j$).
So $\xi_j \le 1 + h_{j-1} - (k+1 - r_j) + (k-r_j) = h_{j-1}$.
\qed

\bigskip

{\bf The inverse bijection.}
We claim that this mapping $F \mapsto (\omega,\xi)$ is a bijection
from the set of ordered forests of increasing ordered trees
on the vertex set $[n]$ with $k$ components
to the set of labeled reversed partial \L{}ukasiewicz paths
from $(0,k)$ to $(n,0)$
that reach level 0 only at the last step,
with integer labels satisfying $1 \le \xi_j \le h_{j-1}$.
To prove this, we explain the inverse mapping.

Given a labeled reversed partial \L{}ukasiewicz path $(\omega,\xi)$,
where $\omega$ reaches level 0 only at the last step,
and $1 \le \xi_j \le h_{j-1}$ for all $j$,
we build up the ordered forest $F$ vertex-by-vertex:
after stage $j$ we will have an ordered forest $F_j$
in~which some of the vertices are labeled $1,\ldots,j$
and some others are unnumbered ``vacant slots''.
The starting forest $F_0$ has $k$ singleton components,
each of which is a vacant slot (these components are of course ordered).
We now ``read'' the path step-by-step, from $j=1$ through $j=n$.
When we read a step $s_j$ with label $\xi_j$,
we insert a new vertex $j$ into one of the vacant slots of $F_{j-1}$:
namely, the $\xi_j$th vacant slot
in the depth-first-search order of $F_{j-1}$.
We also create $s_j + 1$ new vacant slots that are children of~$j$.
This defines $F_j$.
Since $F_0$ has $k = h_0$ vacant slots,
and at stage $j$ we remove one vacant slot and add $s_j + 1$ new ones,
it follows by induction that $F_j$ has $h_j$ vacant slots.
(In particular, the placement of the vertex $j$
 into the  $\xi_j$th vacant slot of $F_{j-1}$ is well-defined,
 since $1 \le \xi_j \le h_{j-1}$ by hypothesis.)
Since by hypothesis the path $\omega$ satisfies
$h_0,\ldots,h_{n-1} \ge 1$ and $h_n = 0$,
it follows that each forest $F_0,\ldots,F_{n-1}$ has at least one vacant slot,
while the forest $F_n$ has no vacant slot.
We define $F = F_n$.

It is fairly clear that this insertion algorithm
defines a map $(\omega,\xi) \mapsto F$
that is indeed the inverse of the mapping $F \mapsto (\omega,\xi)$
defined previously:
this follows from the proof of Lemma~\ref{lemma.hj}
and the definition of the insertion algorithm.

\bigskip

{\bf Computation of the weights.}
We want to enumerate ordered forests of increasing ordered trees
on the vertex set $[n]$ with $k$ components,
in~which each vertex at level $L$ with $i$ children
gets a weight $\phi_i^{[L]}$.
We use the bijection to push these weights from the forests
to the labeled reversed partial \L{}ukasiewicz paths.
Given a forest $F$, each vertex $j \in [n]$
contributes a weight $\phi_{\deg(j)}^{[\lev(j)]}$.
Under the bijection, this vertex 
is mapped to a step $s_j = \deg(j) - 1$
from height $h_{j-1} = \lev(j)$ to height $h_j = h_{j-1} + s_j$.
Therefore, the weight in the labeled path $(\omega,\xi)$
corresponding to this vertex is $\phi_{s_j+1}^{[h_{j-1}]}$,
and the weight of the labeled path $(\omega,\xi)$
is the product of these weights over $1 \le j \le n$.

Now we sum over the labels $\xi$ to get the total weight
for each path $\omega$:
summing over $\xi_j$ gives a factor $h_{j-1}$.
Therefore, the weight in the reversed partial \L{}ukasiewicz path
for a step~$s$ ($-1 \le s < \infty$) starting from height~$h$ will be 
\be
   W(s,h) \;=\;  h \, \phi_{s+1}^{[h]}
   \;.
 \label{eq.Wsh}
\ee
Note that $W(s,0) = 0$;
this implements automatically the constraint that
the reversed partial \L{}ukasiewicz path is not allowed to reach level 0
before the last step.

We now want to read the path $\omega$ backwards,
so that it becomes an ordinary partial \L{}ukasiewicz path $\omegahat$
from $(0,0)$ to $(n,k)$.
A step~$s$ starting at height~$h$ in $\omega$
becomes a step~$s' = -s$ starting at height~$h' = h+s$ in $\omegahat$.
Therefore, in the ordinary partial \L{}ukasiewicz path $\omegahat$,
the weight will be
\be
   W'(s',h')  \;=\;  (h'+s') \, \phi_{1-s'}^{[h'+s']}
   \;.
 \label{eq.Wsh.prime}
\ee
That is, a step from height $i$ to height $j$ gets a weight
\be
   p^{\rm ord}_{ij}  \;=\;  W(j-i,i)  \;=\;  j \, \phi_{i-j+1}^{[j]}
   \;.
\ee
(Note that $p^{\rm ord}_{i0} = 0$, i.e.\ steps to level 0 are forbidden.)
Then $P^{\rm ord} = (p^{\rm ord}_{ij})_{i,j \ge 0}$
is the production matrix for the triangle
$(k! \, \Lhat_{n,k}(\bphihat))_{n,k \geq 0}$
that enumerates ordered forests.
This proves Proposition~\ref{prop.prodmat.phiL}(a).

We then apply Lemma~\ref{lemma.production.AB}
with $B = \diag\big( (1/k!)_{k \ge 0} \big)$,
working temporarily in the ring $\Q[\bphihat]$.
It follows that the production matrix $P$
for the triangle $(\Lhat_{n,k}(\bphihat))_{n,k \geq 0}$
is given by $P = B^{-1} P^{\rm ord} B$,
which is precisely \reff{eq.prop.prodmat.phiL.1}.
This proves Proposition~\ref{prop.prodmat.phiL}(b).

It also follows that the polynomial $\Lhat_{n,k}(\bphihat)$ has a factor
$\phi_0^{[1]} \phi_0^{[2]} \cdots \phi_0^{[k]}$,
since every partial \L{}ukasiewicz path from $(0,0)$ to $(n,k)$
must have rises $0 \to 1$, $1 \to 2$, \ldots, $k-1 \to k$.

To prove Proposition~\ref{prop.prodmat.phiL}(c),
we apply Lemma~\ref{lemma.production.AB} once again,
this time with
$\Btilde = \diag\big( (1/\phi_0^{[1]} \cdots \phi_0^{[k]})_{k \ge 0} \big)$,
working temporarily in the ring $\Z[\bphihat,\bphihat^{-1}]$.
Of course the matrix elements of $\Ptilde = \Btilde^{-1} P \Btilde$
lie in the subring $\Z[\bphihat] \subseteq \Z[\bphihat,\bphihat^{-1}]$.

This completes the proof of Proposition~\ref{prop.prodmat.phiL}.
\qed

\medskip

{\bf Remark.}  The reasoning here using Lemma~\ref{lemma.production.AB}
corresponds, at the level of \L{}ukasiewicz paths,
to pairing each $\ell$-fall $i \to i-\ell$ ($\ell \ge 1$)
with the corresponding rises $i-\ell \to i-\ell+1 \to \ldots \to i$
and then transferring the weights (or part of the weights)
from those rises to the $\ell$-fall
(as was done in \cite{latpath_SRTR}).
\myendremark

\medskip

\proofof{Proposition~\ref{prop.prodmat}{\rm (a)}}
Specialize Proposition~\ref{prop.prodmat.phiL}(b) to the case
$\phi_i^{[L]} = \phi_i$ for all $L$.
\qed

\bigskip

Note now that the triangular arrays
$\sfLhat^{\rm ord}$, $\sfLhat$, $\sfLtilde$ are each of the form
$   \left[
   \begin{array}{c|c@{\hspace*{2mm}}c@{\hspace*{2mm}}c}
        1      &   0 & 0 & \cdots \\
        \hline
        0      &     &   &        \\
        0      &     & L^\wedge &        \\[-1mm]
        \vdots &     &   &        \\
   \end{array}
   \right]
$.
So it is of some interest to find the production matrix
for the corresponding submatrices $L^\wedge$.
Let us use the following notation:
For any matrix $M = (m_{ij})_{i,j \ge 0}$,
write $M^\wedge \eqdef \Delta M \Delta^{\rm T}$
for $M$ with its zeroth row and column removed.
We then have:

\begin{corollary}
   \label{cor.prodmat.phiL}
Let $\bphihat = (\phi_i^{[L]})_{i \ge 0,\, L \ge 1}$ be indeterminates,
and work in the ring $\mathbb{Z}[\bphihat]$.
Define the lower-Hessenberg matrices $P^{\rm ord}$, $P$, $\Ptilde$
and the lower-triangular matrices
$\sfLhat^{\rm ord}$, $\sfLhat$, $\sfLtilde$
as in Proposition~\ref{prop.prodmat.phiL}.
Then:
\begin{itemize}
   \item[(a)]  $(\sfLhat^{\rm ord})^\wedge =
                    \phi_0^{[1]} \, \scro((P^{\rm ord})^\wedge)$.
   \item[(b)]  $\sfLhat^\wedge = \phi_0^{[1]} \, \scro(P^\wedge)$.
   \item[(c)]  $\sfLtilde^\wedge = \scro(\Ptilde^\wedge)$.
\end{itemize}
\end{corollary}

\proof
(a)  By Proposition~\ref{prop.prodmat.phiL}(a) we have
$\sfLhat^{\rm ord} = \scro(P^{\rm ord})$.
Now use Lemma~\ref{lemma.down-shifted}
with $P,Q,c,A,B$ replaced by
$(P^{\rm ord})^\wedge, P^{\rm ord}, \phi_0^{[1]}, (\sfLhat^{\rm ord})^\wedge,
 \sfLhat^{\rm ord}$, respectively.

(b) and (c) are analogous.
\qed

\smallskip

{\bf Remark.}
In \cite[sections~12.2 and 12.3]{latpath_SRTR}
we considered the output matrix $\sfLhat^\wedge$ rather than $\sfLhat$,
and therefore obtained the production matrix $P^\wedge$ rather than $P$.
Also, in that paper we considered only the cases $\bphi = \bee$ and $\bh$;
but the proof for the generic case $\bphi$, given here,
is completely analogous and indeed slightly simpler.
\qed

\subsection[Identity for $B_y^{-1} P B_y$ and proof of
   Proposition~\ref{prop.prodmat}(b)]{Identity for $\bm{B_y^{-1} P B_y}$
            and proof of Proposition~\ref{prop.prodmat}(b)}

We now wish to prove the following identity:

\begin{lemma}[Identity for $B_y^{-1} P B_y$]
   \label{lemma.ByinvPBy}
Let $\bphi = (\phi_i)_{i \ge 0}$ and $y$ be indeterminates,
and work in the ring $\Z[\bphi,y]$.
Define the lower-Hessenberg matrix $P = (p_{ij})_{i,j \ge 0}$ by
\be
   p_{ij}
   \;=\; 
   \begin{cases}
      0   &  \textrm{if $j=0$ or $j > i+1$}  \\[1mm]
      {\displaystyle {i! \over (j-1)!}} \, \phi_{i-j+1}
          &  \textrm{if $1 \le j \le i+1$}
   \end{cases}
 \label{eq.lemma.ByinvPBy.1}
\ee
and the unit-lower-triangular $y$-binomial matrix $B_y$ by
\be
   (B_y)_{nk}  \;=\;  \binom{n}{k} \, y^{n-k}
   \;.
\ee
Let $\Delta = (\delta_{i+1,j})_{i,j \ge 0}$
be the matrix with 1 on the superdiagonal and 0 elsewhere.
Then
\be
   B_y^{-1} P B_y   \;=\;   P (I + y \Delta^{\rm T})
   \;.
 \label{eq.lemma.ByinvPBy}
\ee
\end{lemma}

\proof
It is easy to see, using the Chu--Vandermonde identity,
that $B_y B_z = B_{y+z}$ and hence that $B_y^{-1} = B_{-y}$.
Therefore
\be
   (B_y^{-1} P B_y)_{ij}
   \;=\;
   \sum_{k,\ell} (-1)^{i+k} \, {i! \over k! \, (i-k)!} \, y^{i-k}
        \: {k! \over (\ell-1)!} \, \phi_{k-\ell+1}
        \: {\ell! \over j! \, (\ell-j)!} \, y^{\ell-j}
   \;.
 \label{eq.proof.lemma.ByinvPBy.1}
\ee
Then the coefficient of $\phi_m$ in this is (setting $k=\ell+m-1$)
\begin{subeqnarray}
   & & \!\!\!\!
   [\phi_m] \, (B_y^{-1} P B_y)_{ij}
   \;=\;
   (-1)^{i+m-1} \, {i! \over j!} \, y^{i-j+1-m}
      \sum_\ell (-1)^\ell \, {\ell \over (i-\ell-m+1)! \, (\ell-j)!}
      \qquad \nonumber \\[-1mm] \\[1mm]
   & & \qquad =\;
   (-1)^{i+m-1} \, {i! \over j! \, (i-j+1-m)!} \, y^{i-j+1-m}
      \sum_\ell (-1)^\ell \, \ell \, \binom{i-j+1-m}{\ell-j}
   \;.
      \nonumber \\
 \slabel{eq.proof.lemma.ByinvPBy.2.b}
 \label{eq.proof.lemma.ByinvPBy.2}
\end{subeqnarray}
Now
\be
   \sum_\ell (-1)^\ell \, \binom{i-j+1-m}{\ell-j} \, x^\ell
   \;=\;
   (-x)^j \, (1-x)^{i-j+1-m}
   \;,
\ee
so that
\begin{subeqnarray}
   \sum_\ell (-1)^\ell \, \ell \, \binom{i-j+1-m}{\ell-j}
   & = &
   {d \over dx}
   \bigl[ (-x)^j \, (1-x)^{i-j+1-m} \bigr]
   \biggr|_{x=1}
        \\[2mm]
   & = &
   (-1)^j \, \bigl[ j \delta_{m,i-j+1} \,-\, \delta_{m,i-j} \bigr]
   \;.
\end{subeqnarray}
Substituting this into \reff{eq.proof.lemma.ByinvPBy.2.b} gives
\begin{subeqnarray}
   (B_y^{-1} P B_y)_{ij}
   & = &
   {i! \over j!} \, \bigl[ j \phi_{i-j+1} \,+\, y \phi_{i-j} \bigr]
         \\[2mm]
   & = &
   p_{ij} \,+\, y p_{i,j+1}
   \;,
\end{subeqnarray}
which is precisely \reff{eq.lemma.ByinvPBy}.
\qed

\bigskip

We can now prove Proposition~\ref{prop.prodmat}(b):

\proofof{Proposition~\ref{prop.prodmat}{\rm (b)}}
By Proposition~\ref{prop.prodmat}(a),
the matrix $P$ defined in \reff{eq.prop.prodmat}/\reff{eq.lemma.ByinvPBy.1}
is the production matrix for the generic Lah triangle $\sfL$.
By Lemma~\ref{lemma.production.AB},
the production matrix for $\sfL B_y$ is then $B_y^{-1} P B_y$;
and by Lemma~\ref{lemma.ByinvPBy} this equals $P (I + y \Delta^{\rm T})$.
\qed

\subsection{Total positivity of the production matrix}

We shall use the following general lemma:

\begin{lemma}
   \label{lemma.diagmult.TP}
Let $A = (a_{ij})_{i,j \ge 0}$ be a lower-triangular matrix
with entries in a partially ordered commutative ring $R$,
and let $\bd = (d_i)_{i \ge 1}$.
Define the lower-triangular matrix $B = (b_{ij})_{i,j \ge 0}$ by
\be
   b_{ij}  \;=\;  d_{j+1} d_{j+2} \cdots d_i \, a_{ij}
   \;.
\ee
Then:
\begin{itemize}
   \item[(a)] If $A$ is TP${}_r$ and $\bd$ are indeterminates,
      then $B$ is TP${}_r$ in the ring $R[\bd]$ equipped with
      the coefficientwise order.
   \item[(b)] If $A$ is TP${}_r$ and $\bd$ are nonnegative elements of $R$,
      then $B$ is TP${}_r$ in the ring $R$.
\end{itemize}
\end{lemma}

\proof
(a) Let $\bd = (d_i)_{i \ge 1}$ be commuting indeterminates,
and let us work in the ring $R[\bd,\bd^{-1}]$
equipped with the coefficientwise order.
Let $D = \diag(1,\, d_1,\, d_1 d_2,\, \ldots)$.
Then $D$ is invertible, and both $D$ and
$D^{-1} = \diag(1,\, d_1^{-1},\, d_1^{-1} d_2^{-1},\, \ldots)$
have nonnegative elements.
It follows that $B = D A D^{-1}$ is TP${}_r$ in the ring $R[\bd,\bd^{-1}]$
equipped with the coefficientwise order.
But the matrix elements $b_{ij}$
actually belong to the subring $R[\bd] \subseteq R[\bd,\bd^{-1}]$.
So $B$ is TP${}_r$ in the ring $R[\bd]$
equipped with the coefficientwise order.

(b) follows from (a) by specializing indeterminates.
\qed

\begin{corollary}
   \label{cor.diagmult.TP}
Let $(\phi_i^{[L]})_{i \ge 0,\, L \ge 1}$ be
elements of a partially ordered commutative ring~$R$,
with $\phi_i^{[L]} \eqdef 0$ for $i < 0$.
Suppose that the lower-triangular matrix
$\Phi = (\phi_{i-j}^{[j+1]})_{i,j \ge 0}$
is TP${}_r$.
Then the matrices $P^{\rm ord}$ and $P$
defined by \reff{eq.prop.prodmat.phiL.0}/\reff{eq.prop.prodmat.phiL.1}
are also TP${}_r$.
\end{corollary}

\proof
(a) We have $P^{\rm ord} = \Phi \Delta D$
where $D = \diag\big( (j)_{j \ge 0} \big)$.

(b) Applying Lemma~\ref{lemma.diagmult.TP}(b) with $A = \Phi$
and $d_i = i!$,
we see that the lower-triangular matrix $P' = (p'_{ij})_{i,j \ge 0}$
with entries $p'_{ij} = (i!/j!) \, \phi_{i-j}^{[j+1]}$ is TP${}_r$.
But then $P = P' \Delta$ is also TP${}_r$.
\qed

%

The doubly-indexed sequence
$(\phi_i^{[L]})_{i \ge 0,\, L \ge 1}$ is very general,
but precisely because of its generality it is somewhat difficult to work with:
indeed, the corresponding matrix $\Phi$ is a
{\em completely arbitrary}\/ lower-triangular matrix,
for which it may or may not be feasible to determine its total positivity.
It is therefore of interest to consider specializations
for which the total positivity may be proven more easily.
One such specialization is the following:
Let $\bphi= (\phi_i)_{i \geq 0}$ and $\bfc= (c_L)_{L \ge 1}$
be two sequences of indeterminates,
and set $\phi_i^{[L]} = \phi_i c_L$;
we denote this specialization by the shorthand $\bphihat = \bphi \bfc$.
We then have the following easy fact:

\begin{lemma}
   \label{lemma.phic}
Let $\bphi= (\phi_i)_{i \geq 0}$ be a sequence
in a partially ordered commutative ring~$R$,
with $\phi_i \eqdef 0$ for $i < 0$,
and let $\bc = (c_L)_{L \ge 1}$;
and define the lower-triangular matrix
$\Phi = (\phi_{i-j} c_{j+1})_{i,j \ge 0}$
Then:
\begin{itemize}
   \item[(a)] If $\bphi$ is Toeplitz-TP${}_r$
      and $\bc$ are indeterminates,
      then $\Phi$ is TP${}_r$ in the ring $R[\bfc]$ equipped with
      the coefficientwise order.
   \item[(b)] If $\bphi$ is Toeplitz-TP${}_r$
      and $\bc$ are nonnegative elements of $R$,
      then $\Phi$ is TP${}_r$ in the ring $R$.
\end{itemize}
\end{lemma}

\proof
We have $\Phi = T_\infty(\bphi) \, \diag\big( (c_{j+1})_{j \ge 0} \bigr)$.
\qed

By combining Lemma~\ref{lemma.phic} and Corollary~\ref{cor.diagmult.TP},
we deduce, under the same hypotheses, the TP${}_r$ property
for the production matrices $P^{\rm ord}$ and $P$
evaluated at $\bphihat = \bphi \bfc$.

\subsection{Generalization of Theorem~1.1(a,c)}

We can now state and prove a generalization of Theorem~1.1(a,c):

\begin{theorem}[Total positivity of the refined generic Lah polynomials]
   \label{thm1.1.generalized}
\hfill\break\noindent
Let $\bphihat = (\phi_i^{[L]})_{i \ge 0,\, L \ge 1}$ be
elements of a partially ordered commutative ring~$R$,
with $\phi_i^{[L]} \eqdef 0$ for $i < 0$,
such that the lower-triangular matrix $\Phi = (\phi_{i-j}^{[j+1]})_{i,j \ge 0}$
is TP${}_r$.
Then:
\begin{itemize}
   \item[(a)]  The lower-triangular matrix
       $\sfLhat(\bphihat) = (\Lhat_{n,k}(\bphihat))_{n,k \ge 0}$
       is TP${}_r$.
   \item[(c)]  The sequence $\bLhat^\triangle(\bphihat) =
                     (\Lhat_{n+1,1}(\bphihat))_{n \ge 0}$
       is Hankel-TP${}_r$.
\end{itemize}
\end{theorem}

\proof
(a) By Corollary~\ref{cor.diagmult.TP}, the production matrix $P$
defined in \reff{eq.prop.prodmat.phiL.1} is TP${}_r$.
By Proposition~\ref{prop.prodmat.phiL}(b),
the corresponding output matrix is $\sfLhat = \scro(P)$.
So Theorem~\ref{thm.iteration.homo} implies that $\sfLhat$ is TP${}_r$.

(c) By Corollary~\ref{cor.diagmult.TP}, the matrix $P$ is TP${}_r$;
hence so is $P^\wedge = \Delta P \Delta^{\rm T}$.
By Corollary~\ref{cor.prodmat.phiL}(b),
we have $\sfLhat^\wedge = \phi_0^{[1]} \, \scro(P^\wedge)$.
So Theorem~\ref{thm.iteration2bis} implies that
the zeroth column of $\sfLhat^\wedge$ is Hankel-TP${}_r$.
But that is precisely $(\Lhat_{n+1,1}(\bphihat))_{n \ge 0}$.
\qed

\medskip

{\bf Remark.}
Since an {\em arbitrary}\/ lower-triangular matrix $\Phi$
can be written in the form $\Phi = (\phi_{i-j}^{[j+1]})_{i,j \ge 0}$,
it follows that $\Phi \mapsto \sfLhat(\bphihat)$
is a well-defined polynomial mapping
of the lower-triangular matrices into themselves,
which preserves total positivity of each order $r$.
However, this mapping seems rather complicated,
even when restricted to Toeplitz matrices $\Phi$
(see the comments in Section~\ref{sec.exponential} below).
It would be interesting to better understand this mapping
from an algebraic point of view.
\myendremark

\medskip

As an immediate consequence of Lemma~\ref{lemma.phic}
and Theorem~\ref{thm1.1.generalized}, we have:

\begin{corollary}
   \label{cor.thm1.1.generalized}
Let $\bphi= (\phi_i)_{i \geq 0}$ be a sequence
in a partially ordered commutative ring~$R$,
and let $\bfc = (c_L)_{L \ge 1}$ be indeterminates.
If $\bphi$ is Toeplitz-TP${}_r$, then:
\begin{itemize}
   \item[(a)]  The lower-triangular matrix
       $\sfLhat(\bphi,\bfc) = (\Lhat_{n,k}(\bphi\bfc))_{n,k \ge 0}$
       is TP${}_r$ in the ring $R[\bfc]$ equipped with
      the coefficientwise order.
   \item[(c)]  The sequence $\bLhat^\triangle(\bphi,\bfc) =
                     (\Lhat_{n+1,1}(\bphi\bfc))_{n \ge 0}$
       is Hankel-TP${}_r$ in the ring $R[\bfc]$ equipped with
      the coefficientwise order.
\end{itemize}
\end{corollary}

\proofof{Theorem~\ref{thm1.1}{\rm (a,c)}}
Specialize Corollary~\ref{cor.thm1.1.generalized}(a,c) to $\bfc = \bone$.
\qed

\subsection{Completion of the proofs}

\proofof{Theorem~\ref{thm1.1}{\rm (b)}}
Since the zeroth column of the matrix $\sfL B_y$
is given by the Lah polynomials $L_n(\bphi,y)$,
Theorem~\ref{thm1.1}(b) is an immediate consequence
of Proposition~\ref{prop.prodmat}(b), Corollary~\ref{cor.diagmult.TP}
and Theorem~\ref{thm.iteration2bis}.
\qed

\proofof{Corollary~\ref{cor1.2}}
The Jacobi--Trudi identity
\cite[Theorem~7.16.1 and Corollary~7.16.2]{Stanley_99}
expresses all the Toeplitz minors of $\bee$ or $\bh$
as skew Schur functions.
Furthermore, every skew Schur function
is a nonnegative linear combination of Schur functions
(Littlewood--Richardson coefficients \cite[Section~7.A1.3]{Stanley_99}).
So the sequences $\bee$ and $\bh$ are Toeplitz-totally positive
with respect to the Schur order.
Corollary~\ref{cor1.2} is then an immediate consequence
of Theorem~\ref{thm1.1}.
\qed

Furthermore, by using Corollary~\ref{cor.thm1.1.generalized} here
in place of Theorem~\ref{thm1.1}(a,c),
we obtain a generalization of Corollary~\ref{cor1.2}(a,c),
whose precise statement we leave to the reader;
and likewise for Corollary~\ref{cor1.3}.

\section{Differential operators for the multivariate Lah polynomials}
   \label{sec.differential}

\subsection{Differential operator for positive type}

In \cite[Proposition~12.6]{latpath_SRTR} we gave expressions
for the multivariate Eulerian polynomials of positive type
in terms of the action of certain first-order linear differential operators.
Translated to our current notation, we proved the following:\footnote{
   Strictly speaking, what we proved in \cite[Proposition~12.6]{latpath_SRTR},
   when translated to our current notation,
   puts $\delta_{n1} \delta_{k1}$ instead of $\delta_{n0} \delta_{k0}$
   in \reff{eq.scrqnk.differential}, and holds only for $n,k \ge 1$.
   But it is then easy to see that also \reff{eq.scrqnk.differential}
   holds as written for $n,k \ge 0$.

   Also, the statement in \cite[Proposition~12.6]{latpath_SRTR}
   applied only to $r \ge 2$.
   But for $r=1$ we have
   $L_{n,k}^{(1)+}(x_1) = \stirlingsubset{n}{k} x_1^{n-k}$,
   so that \reff{eq.scrqnk.differential} is the well-known recurrence
   for the Stirling subset numbers (note that $\scrd_1 = 0$).
}

\begin{proposition}
{$\!\!\!$ \bf \protect\cite[Proposition~12.6]{latpath_SRTR}\ }
   \label{prop.MVeulerian.differential}
For every integer ${r \ge 1}$, we have
\begin{eqnarray}
   L_{n,1}^{(r)+}(\bfx)  & = &
      \Bigl( \scrd_r \,+\, \sum_{i=1}^r x_i \Bigr)^{\! n-1} \: 1
      \qquad\hbox{for $n \ge 1$}
          \label{eq.scrq.differential}  \\
   L_{n,k}^{(r)+}(\bfx)  & = &
      \Bigl( \scrd_r \,+\, k \sum_{i=1}^r x_i \Bigr)
      L_{n-1,k}^{(r)+}(\bfx)  \:+\:  L_{n-1,k-1}^{(r)+}(\bfx)
      \:+\: \delta_{n0} \delta_{k0}
      \quad\hbox{for $n,k \ge 0$}
      \qquad
          \label{eq.scrqnk.differential}
\end{eqnarray}
where
\be
   \scrd_r
   \;=\;
   \sum_{i=1}^r \biggl( \! x_i
                        \!\!\!\!
                            \sum_{\begin{scarray}
                                    1 \le j \le r \\
                                    j \ne i
                                  \end{scarray}}
                        \!\!\!\!\!\!
                            x_j \!
                \biggr) \, {\partial \over \partial x_i}
   \;.
 \label{def.dm}
\ee
\end{proposition}

Now we would like to extend this to give a
differential expression for the row-generating polynomials
$L_n^{(r)+}({\bf x},y)$:

\begin{proposition}
   \label{prop.diff.op}
For every integer $r\geq 1$, we have:
\begin{equation}
   L_n^{(r)+}({\bf x},y)  \;=\;  (\scrdtilde_r +y)^n \: 1
   \;,
 \label{eq.prop.diff.op}
\end{equation}
where
\be
   \scrdtilde_r
   \;=\;
   \scrd_r \:+\: \sum_{i=1}^r  x_i y \, {\partial \over \partial y}
   \;.
  \label{def.diff.opp}
\ee
and $\scrd_r$ is defined in \reff{def.dm}.
\end{proposition}

\proof
Multiply \reff{eq.scrqnk.differential} by $y^k$ and sum over $k$:
the factor $k$ becomes $y \, \partial/\partial y$, and we have
\be
   L_n^{(r)+}(\bfx,y)
   \;=\;
      \Bigl( \scrd_r \,+\, \sum_{i=1}^r x_i y \, {\partial \over \partial y}
                     \,+\, y  \Bigr)
      L_{n-1}^{(r)+}(\bfx,y)
      \:+\: \delta_{n0} \delta_{k0}
   \;.
\ee
Iterating this yields \reff{eq.prop.diff.op}.
\qed

\subsection{Differential operator for negative type}

Similarly, in \cite[Proposition~12.26]{latpath_SRTR} we gave expressions
for the multivariate Eulerian polynomials of negative type
in terms of the action of certain first-order linear differential operators.
Translated to our current notation, we proved the following:\footnote{
   Strictly speaking, what we proved in \cite[Proposition~12.26]{latpath_SRTR},
   when translated to our current notation,
   puts $\delta_{n1} \delta_{k1}$ instead of $\delta_{n0} \delta_{k0}$
   in \reff{eq.scrqnk.differential.multi}, and holds only for $n,k \ge 1$.
   But it is then easy to see that also \reff{eq.scrqnk.differential.multi}
   holds as written for $n,k \ge 0$.
}

\begin{proposition}
{$\!\!\!$ \bf \protect\cite[Proposition~12.26]{latpath_SRTR}\ }
   \label{prop.MVeulerian.differential.multi}
For every integer ${r \ge 1}$, we have
\begin{eqnarray}
   L_{n,1}^{(r)-}(\bfx)  & = &
      \Bigl( \scrd^-_r \,+\, \sum_{i=1}^r x_i \Bigr)^{\! n-1} \: 1
      \qquad\hbox{for $n \ge 1$}
          \label{eq.scrq.differential.multi}  \\
   L_{n,k}^{(r)-}(\bfx)  & = &
      \Bigl( \scrd^-_r \,+\, k \sum_{i=1}^r x_i \Bigr)
      L_{n-1,k}^{(r)-}(\bfx)  \:+\:  L_{n-1,k-1}^{(r)-}(\bfx)
      \:+\: \delta_{n0} \delta_{k0}
      \quad\hbox{for $n,k \ge 0$}
      \qquad
          \label{eq.scrqnk.differential.multi}
\end{eqnarray}
where
\be
   \scrd^-_r
   \;=\;
   \sum_{i=1}^r \biggl( \! x_i^2  \:+\:
                           x_i \sum_{j=1}^r x_j \!
                \biggr) \, {\partial \over \partial x_i}
   \;.
 \label{def.dm-}
\ee
\end{proposition}

Now we would like to extend this to give a
differential expression for the row-generating polynomials
$L_n^{(r)-}({\bf x},y)$,
analogously to what we did for the positive type.
The result is:

\begin{proposition}
   \label{prop.diff.op.minus}
For every integer $r\geq 1$, we have:
\begin{equation}
   L_n^{(r)-}({\bf x},y)  \;=\;  (\scrdtilde_r^- +y)^n \: 1
   \;,
 \label{eq.prop.diff.op.minus}
\end{equation}
where
\be
   \scrdtilde_r^-
   \;=\;
   \scrd^-_r \:+\: \sum_{i=1}^r  x_i y \, {\partial \over \partial y}
   \;.
  \label{def.diff.opp.minus}
\ee
and $\scrd^-_r$ is defined in \reff{def.dm-}.
\end{proposition}

\proof
Identical to the proof of Proposition~\ref{prop.diff.op},
but with $\scrd_r$ replaced by $\scrd^-_r$.
\qed

\section{Proof of Theorem~\ref{thm.lah.S-fraction} by the Euler--Gauss recurrence method}
   \label{sec.euler-gauss}

In the Introduction we explained how,
for the multivariate Lah polynomials of {\em positive}\/ type,
the matrix $B_y^{-1} P B_y = P (I + y \Delta^{\rm T})$,
which is the production matrix for $\sfL B_y$,
has the bidiagonal factorization \reff{eq.lah.prodmat.y};
and we explained how this in turn implies,
by virtue of \reff{eq.prop.contraction},
that the multivariate Lah polynomials of positive type $L_n^{(r)+}(\bfx,y)$
are given by an $r$-branched S-fraction with coefficients
\be
   \balpha \;=\; (\alpha_i)_{i \ge r}
   \;=\;
   y,x_1,\ldots,x_r,y,2x_1,\ldots,2x_r,y,3x_1,\ldots,3x_r, \,\ldots
   \;,
 \label{eq.lah.alphas.bis}
\ee
as stated in Theorem~\ref{thm.lah.S-fraction}.

In this section we would like to give a second (and completely independent)
proof of Theorem~\ref{thm.lah.S-fraction},
based on the Euler--Gauss recurrence method for proving continued fractions,
generalized to $m$-S-fractions as in \cite[Proposition~2.3]{latpath_SRTR}.
Let us recall briefly the method:
if $(g_k(t))_{k \ge -1}$ are formal power series with constant term 1
(with coefficients in some commutative ring $R$)
satisfying a recurrence
\be
   g_k(t) - g_{k-1}(t)  \;=\; \alpha_{k+m} t \, g_{k+m}(t)
   \qquad\hbox{for } k \ge 0
 \label{eq.recurrence.gkm.0.bis}
\ee
for some coefficients $\balpha = (\alpha_i)_{i \ge m}$ in $R$,
then $g_0(t)/g_{-1}(t) = \sum_{n=0}^\infty S^{(m)}_n(\balpha) \, t^n$,
where $S^{(m)}_n(\balpha)$ is the $m$-Stieltjes--Rogers polynomial
evaluated at the specified values $\balpha$.

As in \cite[sections~12.2.4 and 12.3.4]{latpath_SRTR},
we will apply this method with the choice $g_{-1}(t) = 1$.
We need to find series $(g_k(t))_{k \ge 0}$ with constant term 1
satisfying \reff{eq.recurrence.gkm.0.bis},
where here $m=r$.
Let us write $g_k(t) = \sum_{n=0}^\infty g_{k,n} \, t^n$
and define $g_{k,n} = 0$ for $n < 0$.
Then \reff{eq.recurrence.gkm.0.bis} can be written as
\be
   g_{k,n} \,-\, g_{k-1,n}  \;=\;  \alpha_{k+r} \, g_{k+r,n-1}
   \qquad\hbox{for } k,n \ge 0
   \;.
 \label{eqrec.gkn}
\ee
Here are the required $g_{k,n}$:

\begin{proposition}[Euler--Gauss recurrence for multivariate Lah polynomials
of positive type]
   \label{prop.euler-gauss.quasi-affine}
Let $\bfx = (x_1,\ldots,x_r)$ be indeterminates;
we work in the ring $R = \Z[\bfx]$.
Set $g_{k,n} = \delta_{n0}$ for $k < 0$,
and then define $g_{k,n}$ for $k,n \ge 0$ by the recurrence
\be
   g_{k,n}
   \;=\;
   \Bigl( \scrdtilde_r \,+\, \sum_{i=1}^r \alpha_{k+i} \Bigr) g_{k,n-1}
       \:+\:  g_{k-r,n}
 \label{def.gkn}
\ee
where $\scrdtilde_r$ is given by \reff{def.diff.opp}
and $\balpha$ are given by
\be
   \balpha \;=\; (\alpha_i)_{i \ge 0}
      \;=\;
   \underbrace{0,\ldots,0}_\text{$r$ times},
       y,x_1,\ldots,x_r,y,2x_1,\ldots,2x_r,y,3x_1,\ldots,3x_r, \,\ldots
 \label{eq.lah.alphas.bis2}
\ee
or in detail
\be
   \alpha_i
   \;=\;
   \begin{cases}
        \Big\lfloor\displaystyle\frac{i}{r+1} \Big\rfloor \, x_{j+1}
                   & \textrm{if $i \equiv j \bmod r\!+\!1$
                                 with $0 \le j \le r\!-\!1$}
                             \\[3mm]
        y          & \textrm{if $i \equiv r \bmod r\!+\!1$}
   \end{cases}
 \label{eq.lah.alphas.bis3}
\ee
Then:
\begin{itemize}
   \item[(a)]  $g_{k,0} = 1$ for all $k \in \Z$.
\\[-5mm]
   \item[(b)]  $(g_{k,n})$ satisfies the recurrence
      \reff{eqrec.gkn} for all $k \ge -r$ and $n \ge 0$.
\\[-5mm]
   \item[(c)]  $S^{(r)}_n(\balpha) = g_{0,n} = (\scrdtilde_r + y)^n \, 1$.
\end{itemize}
\end{proposition}

%
%

\proof
(a)  We see trivially using \reff{def.gkn} that $g_{k,0} = 1$
for all $k \in \Z$, i.e.\ $g_k(t)$ has constant term 1.

(b)  We will now prove that the recurrence \reff{eqrec.gkn} holds.
The proof will be by an outer induction on $k$ (in steps of $r$)
in which we encapsulate an inner induction on~$n$.
The base cases $k = -r,\ldots,-1$ for the outer induction
hold trivially because $g_{k,n} = g_{k-1,n} = \delta_{n0}$ for $k < 0$
and $\alpha_0,\ldots,\alpha_{r-1} = 0$.
We now assume that \reff{eqrec.gkn} holds for a given $k$ and all $n \ge 0$;
we want to prove that it still holds when we replace $k$  by $k+r$, i.e.\ that
\begin{equation}
   g_{k+r,n} \:-\: g_{k+r-1,n} \:-\: \alpha_{k+2r} \, g_{k+2r,n-1}
   \;=\;
   0
   \quad \hbox{for all $n \ge 0$}
   \;.
 \label{eqrec.gkn.inn}
\end{equation}
We will prove \reff{eqrec.gkn.inn} by induction on $n$.

Clearly \reff{eqrec.gkn.inn} holds for $n=0$
because $g_{k+r,0} = g_{k+r-1,0} = 1$ [from part~(a)]
and $g_{k+2r,-1} = 0$ (by definition of the $g$'s).

When  $n > 0$, we use \reff{def.gkn}
on each of the three $g$'s on the left-hand side of
\reff{eqrec.gkn.inn}, giving
\begin{subeqnarray}
   g_{k+r,n}
   & = &
   \Bigl( \scrdtilde_r \,+\, \sum_{i=1}^r \alpha_{k+r+i} \Bigr)  g_{k+r,n-1}
       \:+\:  g_{k,n}
       \\
   g_{k+r-1,n}
   & = &
   \Bigl( \scrdtilde_r \,+\, \sum_{i=1}^r \alpha_{k+r-1+i} \Bigr)  g_{k+r-1,n-1}
       \:+\:  g_{k-1,n}
       \\
   \alpha_{k+2r} \, g_{k+2r,n-1}
   & = &
   \alpha_{k+2r} \,
   \Bigl( \scrdtilde_r \,+\, \sum_{i=1}^r \alpha_{k+2r+i} \Bigr)  g_{k+2r,n-2}
       \:+\:  \alpha_{k+2r} \, g_{k+r,n-1}
   \qquad
\end{subeqnarray}
We can then rewrite the left-hand-side of \reff{eqrec.gkn.inn} as
\begin{eqnarray}
   \text{LHS of \reff{eqrec.gkn.inn}}
   & = &
   g_{k,n}-g_{k-1,n}-\alpha_{k+2r}g_{k+r,n-1}
               \nonumber \\[1mm]
   &  &
   \quad+\,
   \sum_{i=1}^r \alpha_{k+r+i} \,
     (g_{k+r,n-1}-g_{k+r-1,n-1}-\alpha_{k+2r}g_{k+2r,n-2})
               \nonumber \\[1mm]
   &  &
   \quad+\,
   \scrdtilde_r (g_{k+r,n-1}-g_{k+r-1,n-1})
      \:-\:
      \alpha_{k+2r} \, \scrdtilde_r g_{k+2r,n-2}
               \nonumber \\[1mm]
   &  &
   \quad-\,
   \sum_{i=1}^r (\alpha_{k+r-1+i}-\alpha_{k+r+i}) \, g_{k+r-1,n-1}
               \nonumber \\[1mm]
   &  &
   \quad-\,
   \sum_{i=1}^r (\alpha_{k+2r+i}-\alpha_{k+r+i})
      \, \alpha_{k+2r} \, g_{k+2r,n-2}
   \;.
 \label{eq.LHS.eqrec.gkn.inn}
\end{eqnarray}
On the right-hand side of \reff{eq.LHS.eqrec.gkn.inn}, the first line is  
\begin{eqnarray}
   & &
   g_{k,n}-g_{k-1,n}-\alpha_{k+2r}g_{k+r,n-1}
         \nonumber \\[1mm]
   & &
   \qquad =\;
   g_{k,n}-g_{k-1,n}-\alpha_{k+r}g_{k+r,n-1} +
      (\alpha_{k+r}-\alpha_{k+2r}) g_{k+r,n-1}
         \nonumber \\[1mm]
   & &
   \qquad =\;
   (\alpha_{k+r}-\alpha_{k+2r})g_{k+r,n-1}
   \;,
\end{eqnarray}
where the first equality is trivial
and the second one comes from the induction hypothesis on~$k$. 
The second line of \reff{eq.LHS.eqrec.gkn.inn}
is zero by the induction hypothesis on~$n$. 
The fourth line of \reff{eq.LHS.eqrec.gkn.inn} is a telescoping sum over $i$,
yielding simply $-(\alpha_{k+r}-\alpha_{k+2r}) \, g_{k+r-1,n-1}$.
For the third line of \reff{eq.LHS.eqrec.gkn.inn}, we use the fact that
$\scrdtilde_r$ is a pure first-order differential operator;
then the Leibniz rule implies that the third line equals
\begin{eqnarray}
   & &
   \scrdtilde_r (g_{k+r,n-1}-g_{k+r-1,n-1})
      \:-\: \alpha_{k+2r}\scrdtilde_r g_{k+2r,n-2}
                \nonumber \\[1mm]
   & &
   \qquad =\;
   \scrdtilde_r (g_{k+r,n-1}-g_{k+r-1,n-1}-\alpha_{k+2r}g_{k+2r,n-2})
      \:+\: g_{k+2r,n-2} \scrdtilde_r \alpha_{k+2r}
                \nonumber \\[1mm]
   & &
   \qquad =\;
   g_{k+2r,n-2} \scrdtilde_r \alpha_{k+2r}
   \;,
\end{eqnarray}
where the last equality comes from the induction hypothesis on $n$.
We now need to do a distinction of cases to compute
$\scrdtilde_r \alpha_{k+2r}$.
If $k +2r\equiv r \bmod r+1$, we have $\alpha_{k+2r}=y$, and so 
\be
   \scrdtilde_r \alpha_{k+2r}
   \;=\;
   \scrdtilde_r y
   \;=\;
   y \sum_{i=1}^r x_i 
   \;=\;
   \alpha_{k+2r} \sum_{i=1}^r x_i
   \;.
\ee
On the other hand, if $k +2r\equiv j \bmod r+1$ with $j \neq r$,
we then have $\alpha_{k+2r}=\lfloor \frac{ k+2r}{r+1}\rfloor x_{j+1}$,
and so 
\be
   \scrdtilde_r \alpha_{k+2r}
   \;=\; \scrdtilde_r \Bigl( \Bigl\lfloor \frac{ k+2r}{r+1} \Bigr\rfloor
                             x_{j+1} \Bigr)
   \;=\; \Bigl\lfloor \frac{ k+2r}{r+1} \Bigr\rfloor x_{j+1}
         \!\!\!\!
          \sum_{\begin{scarray} 
                  1 \le i \le r \\ 
                  i \neq j+1 
                \end{scarray}}
         \!\!\!\!
         x_i
   \;=\; \alpha_{k+2r}
         \!\!\!\!
          \sum_{\begin{scarray} 
                  1 \le i \le r \\ 
                  i \neq j+1 
                \end{scarray}}
         \!\!\!\!
         x_i
   \;.
\ee
Finally, the fifth line of of \reff{eq.LHS.eqrec.gkn.inn} is 
\begin{eqnarray}
   & &
-\sum_{i=1}^r (\alpha_{k+2r+i}-\alpha_{k+r+i}) \, \alpha_{k+2r} \, g_{k+2r,n-2}
        \nonumber \\[1mm]
   & &
   \qquad =\;
   - \Bigl( \alpha_{k+r}-\alpha_{k+2r} +
              \sum_{i=1}^r (\alpha_{k+2r+i}-\alpha_{k+r+i-1})
     \Bigr)\alpha_{k+2r} \, g_{k+2r,n-2}
   \qquad
\end{eqnarray}
by a change of index $i \to i-1$ in the second sum.
Again, we need to do a distinction of cases to compute the sum.
If $k +2r\equiv r \bmod r+1$, we then have
\be
   \sum_{i=1}^r \left(\alpha_{k+2r+i}-\alpha_{k+r+i-1}\right)
   \;=\;
   \sum_{i=1}^r x_i
\ee
by definition of the $\alpha$'s;
whereas when $k +2r\equiv j \bmod r+1$ with $j \neq r$,
we have 
\be
   \sum_{i=1}^r \left(\alpha_{k+2r+i}-\alpha_{k+r+i-1}\right)
   \;=\;
   \!\!\!\!
   \sum_{\begin{scarray}
            1 \le i \le r \\
            i \neq j+1
          \end{scarray}}
   \!\!\!\!
   x_i
   \;.
\ee
And we still have the term
$-(\alpha_{k+r}-\alpha_{k+2r}) \alpha_{k+2r}g_{k+2r,n-2}$.
In both cases, the sum involving the $x_i$'s
cancels between the third and fifth lines;
therefore, all that remains of the third and fifth lines is
$-(\alpha_{k+r}-\alpha_{k+2r}) \alpha_{k+2r}g_{k+2r,n-2}$.

Now, adding all the lines gives 
\be
   (\alpha_{k+r}-\alpha_{k+2r}) \,
   (g_{k+r,n-1}-g_{k+r-1,n-1}-\alpha_{k+2r}g_{k+2r,n-2})
   \;,
\ee
which vanishes by the induction hypothesis on $n$.
This concludes the inductive step in $n$ to prove \reff{eqrec.gkn.inn},
which in turn concludes the induction on $k$ and finishes the proof
of part~(b).

(c) Putting $k=0$ in \reff{def.gkn} gives
$g_{0,n} = (\scrdtilde_r + y) g_{0,n-1} + \delta_{n0}$,
which proves $g_{0,n} = (\scrdtilde_r + y)^n \, 1$;
and this equals $L_n^{(r)+}(\bfx,y)$ by Proposition~\ref{prop.diff.op}.
On the other hand, starting from the recurrence \reff{eqrec.gkn}
and applying the Euler--Gauss recurrence method
\cite[Proposition~2.3]{latpath_SRTR},
we conclude that $g_{0,n} = S_n^{(r)}(\balpha)$.
\qed

\section{Multivariate Lah polynomials in terms of decorated set partitions}
   \label{sec.decorated}

In this section we would like to interpret the
multivariate Lah polynomials of positive type $L_n^{(r)+}(\bfx,y)$
as generating polynomials for partitions of the set $[n]$
in~which each block is ``decorated'' with an additional structure,
where the nature of this structure depends on the value of $r$.

For $r=1,2$ we observed already in the Introduction how this goes:
for $r=1$ the additional structure is empty,
while for $r=2$ it is a linear ordering on the block.
More precisely:

\begin{proposition}
   \label{prop.bell}
The polynomial $L_n^{(1)+}(1,y)$ is the Bell polynomial $B_n(y)$,
that is, the generating polynomial of set partitions of $n$ elements
with a weight $y$ for each block.

More generally, $L_n^{(1)+}(x_1,y)$
is the homogenized Bell polynomial $x_1^n B_n(y/x_1)$.
\end{proposition}

\proof
By definition, the polynomial $L_n^{(1)+}(1,y)$ is the generating
polynomial of unordered forests of increasing unary trees.
Since, given a set of integer labels,
there is only one way to increasingly label a unary tree,
there is a natural bijection between increasing unary trees and sets of labels.
An unordered forest of increasing unary trees is then
an unordered collection of disjoint sets of integers, whose union is $[n]$.
But this is nothing other than a set partition of $[n]$.

The final statement follows from the fact that $L_n^{(1)+}(x_1,y)$
is homogeneous of degree $n$.
\qed

\begin{proposition}
   \label{prop.lah}
The polynomial $L_n^{(2)+}(\bone;y)$ is the Lah polynomial $L_n(y)$,
that is, the generating polynomial of set partitions of $n$ elements
into any number of nonempty lists (= linearly ordered subsets),
with a weight $y$ for each list.

More generally, the polynomial $L_n^{(2)+}(x_1,x_2;y)$
is the generating polynomial of partitions of the set $[n]$
into any number of nonempty lists, with a weight $y$ for each list,
$x_1$ for each descent in a list, and $x_2$ for each ascent in a list.
\end{proposition}

\proof
By definition, the polynomial $L_n^{(2)+}(x_1,x_2;y)$
is the generating polynomial of unordered forests
of increasing binary trees,
with a weight $y$ for each root (or equivalently, each tree)
and a weight $x_1$ (resp.\ $x_2$) for each left (resp.\ right) child.
Now, a classical bijection \cite[pp.~23--25]{Stanley_86}
sends increasing binary trees to permutations.
Since a permutation is the same thing as a list,
applying this classical bijection to each tree of the forest
maps bijectively an unordered forest of increasing binary trees
to a set of lists whose union is $[n]$,
where the number of trees in the forest equals the number of lists.

The second claim comes out naturally, as this classical bijection maps
a left (resp.\ right) child of the tree to a descent (resp.\ ascent)
in the resulting permutation.
\qed

This can be generalized to any positive integer $r$ by using the concept
of {\em Stirling permutation}\/ \cite{Gessel_78,Gessel_78a,Park_94a,Park_94b}
as discussed in \cite[section~12.5]{latpath_SRTR}.
Recall that a word $\bfw = w_1 \cdots w_L$
on a totally ordered alphabet $\mathbb{A}$
is called a \textbfit{Stirling word}
if $i < j < k$ and $w_i = w_k$ imply $w_j \ge w_i$:
that is, between any two occurrences of any letter $a$,
only letters that are larger than or equal to $a$ are allowed.
(Equivalently, between any two successive occurrences of the letter $a$,
only letters that are larger than $a$ are allowed.)
Now let $\mathbb{A}$ be a totally ordered alphabet
of finite cardinality $\ell$, and let $r$ be a nonnegative integer;
we denote by $r \mathbb{A}$ the multiset
consisting of $r$ copies of each letter $a \in \mathbb{A}$.
A~\textbfit{permutation} of $r \mathbb{A}$ is a word $w_1 \cdots w_{r\ell}$
containing exactly $r$ copies of each letter $a \in \mathbb{A}$;
it is called a \textbfit{Stirling permutation} of $r \mathbb{A}$
if it is also a Stirling word.

Now let $n$ be a nonnegative integer.
We define an \textbfit{$\bm{r}$-Stirling set partition} of $[n]$
to be a set partition of $[n]$ in which each block $B$
is decorated by a Stirling permutation of $rB$,
where the total order on $B$ is of course the one inherited from
the usual total order on the integers.
In particular, when $r=0$, the decoration is empty,
and we get back to classical set partitions;
and when $r=1$, we get a partition of the set $[n]$
in which each block is decorated by a permutation of the letters of that block,
or in other words, a partition of the set $[n]$ into nonempty lists.

\begin{proposition}
   \label{prop.stirling}
The polynomial $L_n^{(r)+}(\bone, y)$  is the generating polynomial
of $(r-1)$-Stirling set partitions of $[n]$, with a weight $y$ for each block. 

More generally, the polynomial $L_n^{(r)+}(\bfx, y)$
is the generating polynomial of $(r-1)$-Stirling set partitions of $[n]$,
with a weight $y$ for each block,
a weight $x_i$ ($1 \le i \le r-1$)
for each time the $i$th occurrence of a letter is the end of a descent,
and a weight $x_r$ for each time the last occurrence of a letter
is the beginning of an ascent.
\end{proposition}

\proof
By definition, the polynomial $L_n^{(r)+}(\bfx, y)$ is the
generating polynomial for unordered forests of increasing $r$-ary trees,
with a weight $y$ for each root and a weight $x_i$ for each $i$-child.
Now a classical bijection \cite{Gessel_78a,Janson_11,Kuba_09}
(see also \cite[section~12.5]{latpath_SRTR})
sends increasing $r$-ary trees on the vertex set $[n]$
to Stirling permutations of the multiset $(r-1) [n]$,
such that\footnote{
   See \cite[Lemma~12.34]{latpath_SRTR}
   after some slight translation of notation.
}:
\begin{itemize}
   \item[1)]  A vertex $j$ has a 1-child if and only if
      the first occurrence of the letter $j$ in the word $\bfw$
      is the end of a descent.
   \item[2)]  A vertex $j$ has an $r$-child if and only if
      the last occurrence of the letter $j$ in the word $\bfw$
      is the beginning of an ascent.
   \item[3)]  A vertex $j$ has an $i$-child ($2 \le i \le r-1$) if and only if
      in the word $\bfw$, between the $(i-1)$st and $i$th occurrences
      of the letter $j$ there is a nonempty subword;
      or equivalently, the $(i-1)$st occurrence of the letter $j$
      is the beginning of an ascent;
      or equivalently, the $i$th occurrence of the letter $j$
      is the end of a descent.
\end{itemize}
Applying this bijection to each tree of the forest maps bijectively
unordered forests of increasing $r$-ary trees on the vertex set $[n]$
to a collection of Stirling permutations of multisets $(r-1)B_i$,
where the $B_i$ taken together form a partition of the set $[n]$.
This collection is nothing other than
an $(r-1)$-Stirling set partition of $[n]$.
\qed

\section{Exponential generating functions}  \label{sec.exponential}

By using exponential generating functions together with
the Lagrange inversion formula, we can obtain explicit expressions
for the generic Lah polynomials $L_{n,k}(\bphi)$.
The method is due to Bergeron, Flajolet and Salvy \cite{Bergeron_92};
see also \cite[Chapter~5, especially pp.~364--365]{Bergeron_98}
and \cite[section~12.2.1]{latpath_SRTR}.
We will use Lagrange inversion in the following form \cite{Gessel_16}:
If $A(u)$ is a formal power series
with coefficients in a commutative ring $R$ containing the rationals,
then there exists a unique formal power series $f(t)$
with zero constant term satisfying
\be
   f(t)  \;=\;  t \, A(f(t))
   \;,
\ee
and it is given by
\be
   [t^n] \, f(t)  \;=\;  {1 \over n} \, [u^{n-1}] \, A(u)^n
     \quad\hbox{for $n \ge 1$}
   \;;
\ee
and more generally, if $H(u)$ is any formal power series, then
\be
   [t^n] \, H(f(t))  \;=\;  {1 \over n} \, [u^{n-1}] \, H'(u) \, A(u)^n
     \quad\hbox{for $n \ge 1$}
   \;.
 \label{eq.lagrange.H}
\ee

Let $\bphi = (\phi_i)_{i \ge 0}$ and $y$ be indeterminates;
we will employ formal power series with coefficients in
$\Q[\bphi]$ or $\Q[\bphi,y]$.
Recall that $L_{n,k}(\bphi)$ is the generating polynomial
for unordered forests of increasing ordered trees on $n$ total vertices
with $k$ components,
in which each vertex with $i$ children gets a weight $\phi_i$;
in particular, $L_{n,1}(\bphi)$ is the generating polynomial
for increasing ordered trees.
And $L_n(\bphi,y) = \sum_{k=0}^n L_{n,k}(\bphi) \: y^k$
are the row-generating polynomials.
Define now the exponential generating function for trees:
\be
   U(t)  \;=\;  \sum_{n=1}^\infty L_{n,1}(\bphi) \, {t^n \over n!}
   \;.
 \label{def.U}
\ee
It is easy to see that the exponential generating function
for $k$-component unordered forests is
\be
   {U(t)^k \over k!}  \;=\;  \sum_{n=0}^\infty L_{n,k}(\bphi) \, {t^n \over n!}
   \;.
 \label{eq.Uk}
\ee
Multiplying this by $y^k$ and summing over $k$ then gives
the exponential generating function for the row-generating polynomials:
\be
   e^{y U(t)}  \;=\;  \sum_{n=0}^\infty L_n(\bphi,y) \, {t^n \over n!}
   \;.
\ee

Here is the key step:
standard enumerative arguments \cite[Theorem~1]{Bergeron_92}
show that $U(t)$ satisfies the ordinary differential equation
\be
   U'(t) \;=\; \Phi(U(t))
   \;,
 \label{eq.bergeron.ODE}
\ee
where $\Phi(w) \eqdef \sum_{k=0}^\infty \phi_k w^k$
is the ordinary generating function for $\bphi$.
At this point it is convenient to specialize to $\phi_0 = 1$;
at the end we can restore the missing factors of $\phi_0$
by recalling that $L_{n,k}(\bphi)$ is homogeneous of degree $n$ in $\bphi$.
(Alternatively, we could keep $\phi_0$ and work over the ring
 $\Q[\bphi,\phi_0^{-1}]$ instead of $\Q[\bphi]$.)
We can now rewrite the differential equation \reff{eq.bergeron.ODE}
as the implicit equation
\be
   t  \;=\; \int\limits_0^{U(t)} {dw \over \Phi(w)}
   \;.
 \label{eq.bergeron}
\ee
Introducing $\Psi(w) \eqdef 1/\Phi(w) \eqdef 1 + \sum_{i=1}^\infty \psi_i w^i$,
we then have
\be
   t   \;=\; U(t) \, \widehat{\Psi}(U(t))
   \qquad\hbox{where}\qquad
   \widehat{\Psi}(z) \;=\; 1 + \sum_{i=1}^\infty {\psi_i \over i+1} \, z^i
   \;.
 \label{eq.bergeron.fZ}
\ee 
Solving $U(t) = t/\widehat{\Psi}(U(t))$
by Lagrange inversion \reff{eq.lagrange.H}
with $A(u) = 1/\widehat{\Psi}(u)$ and $H(u) = u^k/k!$ gives
\begin{subeqnarray}
   & &
   \hspace*{-1.45cm}
   L_{n,k}(\bphi) \Bigr|_{\phi_0 = 1}
   \;=\;  n! \: [t^n] \, {U(t)^k \over k!}
   \;=\; {(n-1)! \over (k-1)!} \: [z^{n-k}] \, \widehat{\Psi}(z)^{-n}
        \\[2mm]
   & & \;
   =\;
   {(n-1)! \over (k-1)!}
    \!\! \sum_{\begin{scarray}
                   l_1, l_2, \ldots \ge 0 \\
                   \sum i l_i = n-k
                \end{scarray}}
          \!\!\!\!
          \binom{-n}{-n-\sum l_i,\, l_1,\, l_2,\, \ldots}
          \prod_{i=1}^\infty \Bigl( {\psi_i \over i+1} \Bigr) ^{\! l_i}
        \\[2mm]
   & & \;
   =\;
   {(n-1)! \over (k-1)!}
    \!\! \sum_{\begin{scarray}
                   l_1, l_2, \ldots \ge 0 \\
                   \sum i l_i = n-k
                \end{scarray}}
          \!\!
          (-1)^{\sum l_i} \,
          \binom{n+\sum l_i - 1}{n-1,\, l_1,\, l_2,\, \ldots}
          \prod_{i=1}^\infty \Bigl( {\psi_i \over i+1} \Bigr) ^{\! l_i}
   \,.
   \;
 \slabel{eq.bergeron.Zn.c}
 \label{eq.bergeron.Zn}
\end{subeqnarray}


\begin{example}[Forests of increasing multi-unary trees]
   \label{exam.phi=1}
\rm
The multivariate Lah polynomials of negative type
$L_{n,k}^{(r)-}(x_1,\ldots,x_r) = L_{n,k}(\bh(x_1,\ldots,x_r))$
specialized to $r=1$ and $x_1 = 1$
--- which count unordered forests of increasing multi-unary trees ---
correspond to $\phi = \bone$,
hence $\psi_1 = -1$ and $\psi_i = 0$ for $i \ge 2$.
It follows that in \reff{eq.bergeron.Zn.c} we have
$l_1 = n-k$ and $l_i = 0$ for $i \ge 2$, hence
\be
   L_{n,k}(\bone)
   \;=\;
   {(2n-k-1)! \over 2^{n-k} \, (n-k)! \, (k-1)!}
   \;.
\ee
These are a shifted version of the coefficients
of the (reversed) Bessel polynomials \cite[A001497]{OEIS}.
\myendremark
\end{example}

\begin{example}[$r=\infty$]
   \label{exam.r=infty}
\rm
If we consider $r$-ary or multi-$r$-ary trees and take $r \to\infty$,
then after an appropriate rescaling we get $\phi_i = 1/i!$ and $\Phi(w) = e^w$.
Solving \reff{eq.bergeron.ODE} gives $U(t) = -\log(1-t)$ and hence
\be
   L_{n,k}(\bphi)
   \;=\;
   {n! \over k!} \; [t^n] \, (-\log(1-t))^k
   \;=\;
   \stirlingcycle{n}{k}
   \;,
\ee
the Stirling cycle numbers \cite[A132393]{OEIS}.
\myendremark
\end{example}

\medskip

{\bf Final remark.}  Here \reff{eq.bergeron.Zn.c} gives a nice explicit
expression for $L_{n,k}(\bphi)$,
but it~is in~terms of the coefficients $\psi_i$ in $\Psi(w) = 1/\Phi(w)$,
not directly in terms of the $\bphi$.
Indeed, if one computes from \reff{eq.bergeron.Zn.c}
the polynomials $L_{n,k}(\bphi)$,
one finds some coefficients that have modestly (but not hugely) large
prime factors:
for instance, one of the terms in $L_{11,1}(\bphi)$ is
$24950808\, \phi_1^2 \phi_2 \phi_3^2$,
where $24950808 = 2^3 \cdot 3^3 \cdot 115513$;
and one of the terms in $L_{13,1}(\bphi)$ is
$2318149824\, \phi_1^3 \phi_2^3 \phi_3$,
where $2318149824 = 2^6 \cdot 3 \cdot 12073697$.
This suggests that the polynomials $L_{n,k}(\bphi)$
{\em might}\/ not have any simple explicit expression.
Or alternatively, they might have a simple explicit expression,
but with coefficients that are given by sums and not just by products.
We leave it as an open problem to find such an expression.
\myendremark

\section{Note Added: Exponential Riordan arrays}
   \label{sec.exp_riordan}

After completing this paper we realized that the key
Proposition~\ref{prop.prodmat}(a),
which expresses the production matrix of the generic Lah triangle
and which we proved combinatorially in Section~\ref{subsec.proofs.prodmat}
by bijection onto labeled partial \L{}ukasiewicz paths,
can also be proven algebraically by using the theory of
exponential Riordan arrays \cite{Deutsch_04,Deutsch_09,Barry_16}.
Here we would like to present briefly this alternate proof.

Let $R$ be a commutative ring containing the rationals,
and let $F(t) = \sum_{n=0}^\infty f_n t^n/n!$
and $G(t) = \sum_{n=1}^\infty g_n t^n/n!$ be formal power series
with coefficients in $R$; we set $g_0 = 0$.
Then the \textbfit{exponential Riordan array} associated to the pair $(F,G)$
--- or equivalently to the pair of sequences
$\bff = (f_n)_{n \ge 0}$ and $\bg = (g_n)_{n \ge 1}$ ---
is the infinite lower-triangular matrix
$\scrr[F,G] = (\scrr[F,G]_{nk})_{n,k \ge 0}$ defined by
\be
   \scrr[F,G]_{nk}
   \;=\;
   {n! \over k!} \:
   [t^n] \, F(t) G(t)^k
   \;.
\ee
That is, the $k$th column of $\scrr[F,G]$
has exponential generating function $F(t) G(t)^k/k!$.
Please note that the diagonal elements of $\scrr[F,G]$
are $\scrr[F,G]_{nn} = f_0 g_1^n$,
so the matrix $\scrr[F,G]$ is invertible
in the ring $R^{\N \times \N}_{\rm lt}$ of lower-triangular matrices
if and only if $f_0$ and $g_1$ are invertible in $R$.

We shall use an easy but important result that is sometimes called
the \emph{fundamental theorem of exponential Riordan arrays} (FTERA):

\begin{lemma}[Fundamental theorem of exponential Riordan arrays]
   \label{lemma.FETRA}
Let $\bb = (b_n)_{n \ge 0}$ be a sequence with
exponential generating function $B(t) = \sum_{n=0}^\infty b_n t^n/n!$.
Considering $\bb$ as a column vector and letting $\scrr[F,G]$
act on it by matrix multiplication, we obtain a sequence $\scrr[F,G] \bb$
whose exponential generating function is $F(t) \, B(G(t))$.
\end{lemma}

\proof
We compute
\begin{subeqnarray}
   \sum_{k=0}^n \scrr[F,G]_{nk} \, b_k
   & = &
   \sum_{k=0}^\infty {n! \over k!} \, [t^n] \, F(t) G(t)^k \, b_k
            \\[2mm]
   & = &
   n! \: [t^n] \: F(t) \sum_{k=0}^\infty b_k \, {G(t)^k \over k!}
            \\[2mm]
   & = &
   n! \: [t^n] \: F(t) \, B(G(t))
   \;.
\end{subeqnarray}
\qed

We can now determine the production matrix of an exponential
Riordan array $\scrr[F,G]$:

\begin{theorem}[Production matrices of exponential Riordan arrays]
   \label{thm.riordan.exponential.production}
Let $L$ be a lower-triangular matrix
(with entries in a commutative ring $R$ containing the rationals)
with invertible diagonal entries and $L_{00} = 1$,
and let $P = L^{-1} \Delta L$ be its production matrix.
Then $L$ is an exponential Riordan array
if and only~if $P = (p_{nk})_{n,k \ge 0}$ has the form
\be
   p_{nk}
   \;=\;
   {n! \over k!} \: (z_{n-k} \,+\, k \, a_{n-k+1})
 \label{eq.thm.riordan.exponential.production}
\ee
for some sequences $\ba = (a_n)_{n \ge 0}$ and $\bz = (z_n)_{n \ge 0}$
in $R$.

More precisely, $L = \scrr[F,G]$ if and only~if $P$
is of the form \reff{eq.thm.riordan.exponential.production}
where the ordinary generating functions
$A(t) = \sum_{n=0}^\infty a_n t^n$ and $Z(t) = \sum_{n=0}^\infty z_n t^n$
are connected to $F(t)$ and $G(t)$ by
\be
   G'(t) \;=\; A(G(t))  \;,\qquad
   {F'(t) \over F(t)} \;=\; Z(G(t))
 \label{eq.prop.riordan.exponential.production.1}
\ee
or equivalently
\be
   A(t)  \;=\;  G'(\bar{G}(t))  \;,\qquad
   Z(t)  \;=\;  {F'(\bar{G}(t)) \over F(\bar{G}(t))}
 \label{eq.prop.riordan.exponential.production.2}
\ee
where $\bar{G}(t)$ is the compositional inverse of $G(t)$.
\end{theorem}

\par\bigskip\noindent{\sc Proof}
(mostly contained in \cite[pp.~217--218]{Barry_16}).
Suppose that $L = \scrr[F,G]$.
The hypotheses on $L$ imply that $f_0 = 1$
and that $g_1$ is invertible in $R$;
so $G(t)$ has a compositional inverse.
Now let $P = (p_{nk})_{n,k \ge 0}$ be a matrix;
its column exponential generating functions are, by definition,
$P_k(t) = \sum_{n=0}^\infty p_{nk} \, t^n/n!$.
Applying the FTERA to each column of $P$,
we see that $\scrr[F,G] P$ is a matrix
whose column exponential generating functions
are $\big( F(t) \, P_k(G(t)) \big)_{k \ge 0}$.
On~the other hand, $\Delta \, \scrr[F,G]$
is the matrix $\scrr[F,G]$ with its zeroth row removed
and all other rows shifted upwards,
so it has column exponential generating functions
\be
   {d \over dt} \, \big( F(t) \, G(t)^k/k! \big)
   \;=\;
   {1 \over k!} \: \Big[ F'(t) \, G(t)^k
                         \:+\: k \, F(t) \, G(t)^{k-1} \, G'(t) \Big]
   \;.
\ee
Comparing these two results, we see that
$\Delta \, \scrr[F,G] = \scrr[F,G] \, P$
if and only~if
\be
   P_k(G(t))
   \;=\;
   {1 \over k!} \:
   {F'(t) \, G(t)^k \:+\: k \, F(t) \, G(t)^{k-1} \, G'(t)
    \over
    F(t)}
   \;,
\ee
or in other words
\be
   P_k(t)
   \;=\;
   {1 \over k!}  \:
      \biggl[ {F'(\bar{G}(t)) \over F(\bar{G}(t))} \, t^k
              \:+\: k \, t^{k-1} \, G'(\bar{G}(t))
      \biggr]
   \;.
\ee
Therefore
\begin{subeqnarray}
   p_{nk}
   & = &
   {n! \over k!} \: [t^n] \,
      \biggl[ {F'(\bar{G}(t)) \over F(\bar{G}(t))} \, t^k
              \:+\: k \, t^{k-1} \, G'(\bar{G}(t))
      \biggr]
    \\[2mm]
   & = &
   {n! \over k!} \:
      \biggl[ [t^{n-k}] \: {F'(\bar{G}(t)) \over F(\bar{G}(t))}
              \:+\: k \, [t^{n-k+1}] \: G'(\bar{G}(t))
      \biggr]
    \\[2mm]
   & = &
   {n! \over k!} \: (z_{n-k} \,+\, k \, a_{n-k+1})
\end{subeqnarray}
where $\ba = (a_n)_{n \ge 0}$ and $\bz = (z_n)_{n \ge 0}$
are given by \reff{eq.prop.riordan.exponential.production.2}.

Conversely, suppose that $P = (p_{nk})_{n,k \ge 0}$ has the form
\reff{eq.thm.riordan.exponential.production}.
Define $F(t)$ and $G(t)$
as the unique solutions (in the formal-power-series ring $R[[t]]$)
of the differential equations \reff{eq.prop.riordan.exponential.production.1}
with initial conditions $F(0) = 1$ and $G(0) = 0$.
Then running the foregoing computation backwards
shows that $\Delta \, \scrr[F,G] = \scrr[F,G] \, P$.
\qed

\medskip

\par\medskip\noindent{\sc Alternate Proof of Proposition~\ref{prop.prodmat}}(a).
We use the expressions for the exponential generating functions
of the generic Lah polynomials, which were determined in
Section~\ref{sec.exponential}.
{}From \reff{def.U}/\reff{eq.Uk} we see that the
generic Lah triangle $\sfL = (L_{n,k}(\bphi))_{n,k \ge 0}$
is the exponential Riordan array $\scrr[F,G]$
with $F(t) = 1$ and $G(t) = U(t)$.
Comparing \reff{eq.prop.riordan.exponential.production.1}
with \reff{eq.bergeron.ODE},
we see that $A(t) = \Phi(t)$ and $Z(t) = 0$.
The production matrix \reff{eq.thm.riordan.exponential.production}
then becomes \reff{eq.prop.prodmat},
which proves Proposition~\ref{prop.prodmat}(a).
\qed

\medskip

{\bf Remark.}
This proof shows that the generic Lah triangle
$\sfL = (L_{n,k}(\bphi))_{n,k \ge 0}$
is in fact the {\em general}\/ exponential Riordan array $\scrr[F,G]$
of the ``associated subgroup'' $F=1$,
expressed in terms of its $A$-sequence $\ba = \bphi$.
In this way, the theory of the generic Lah triangle
is {\em equivalent}\/ to the theory of
exponential Riordan arrays of the ``associated subgroup'' $\scrr[1,G]$,
expressed in the combinatorial language of
unordered forests of increasing ordered trees.
It would be interesting to work out the combinatorial interpretation
of exponential Riordan arrays $\scrr[F,G]$ with $F \neq 1$.
\myendremark

\medskip

This algebraic proof of Proposition~\ref{prop.prodmat}(a)
is arguably much simpler than the combinatorial proof
presented in Section~\ref{subsec.proofs.prodmat}.
On the other hand, the combinatorial method seems to be more powerful:
we do not see (at least at present) how to extend the algebraic proof
to obtain the more general Proposition~\ref{prop.prodmat.phiL},
which expresses the production matrix of the refined generic Lah triangle.

\section*{Acknowledgments}


One of us (A.D.S.)\ wishes to thank Xi Chen for helpful conversations
concerning exponential Riordan arrays.

This research was supported in part by
the U.K.~Engineering and Physical Sciences Research Council grant EP/N025636/1.

\clearpage

\appendix

\section{Lah polynomials for $\bm{n \le 7}$}

In this appendix we report the generic Lah polynomials $L_{n,k}(\bphi)$
for $n \le 7$ (specialized for simplicity to $\phi_0 = 1$).
We also report the Lah symmetric functions
$L_{n,k}^{(\infty)+}$ and $L_{n,k}^{(\infty)-}$ for $n \le 7$
in terms of the monomial symmetric functions $m_\lambda$;
from these the reader can easily reconstruct explicit expressions
for the multivariate Lah polynomials
$L_{n,k}^{(r)+}(x_1,\ldots,x_r)$ and $L_{n,k}^{(r)-}(x_1,\ldots,x_r)$
for any chosen value of $r$.
The conversions from $e_\lambda$ or $h_\lambda$ to $m_\lambda$
were performed using the {\tt SymFun} {\sc Mathematica} package
(version~3.1), developed by Curtis Greene and collaborators
\cite{Greene_symfun}.

\subsection{Generic Lah polynomials}

\vspace*{-5mm}

%

\begin{eqnarray*}
L_{1,1}(\bphi)  & = &  1    \\[3mm]
L_{2,1}(\bphi)  & = &  \phi_{1}    \\[0.5mm]
L_{2,2}(\bphi)  & = &  1    \\[3mm]
L_{3,1}(\bphi)  & = &  \phi_{1}^2+2 \phi_{2}    \\[0.5mm]
L_{3,2}(\bphi)  & = &  3 \phi_{1}    \\[0.5mm]
L_{3,3}(\bphi)  & = &  1    \\[3mm]
L_{4,1}(\bphi)  & = &  \phi_{1}^3+8 \phi_{1} \phi_{2}+6 \phi_{3}    \\[0.5mm]
L_{4,2}(\bphi)  & = &  7 \phi_{1}^2+8 \phi_{2}    \\[0.5mm]
L_{4,3}(\bphi)  & = &  6 \phi_{1}    \\[0.5mm]
L_{4,4}(\bphi)  & = &  1    \\[3mm]
L_{5,1}(\bphi)  & = &  \phi_{1}^4+22 \phi_{1}^2 \phi_{2}+16 \phi_{2}^2+42 \phi_{1} \phi_{3}+24 \phi_{4}    \\[0.5mm]
L_{5,2}(\bphi)  & = &  15 \phi_{1}^3+60 \phi_{1} \phi_{2}+30 \phi_{3}    \\[0.5mm]
L_{5,3}(\bphi)  & = &  25 \phi_{1}^2+20 \phi_{2}    \\[0.5mm]
L_{5,4}(\bphi)  & = &  10 \phi_{1}    \\[0.5mm]
L_{5,5}(\bphi)  & = &  1    \\[3mm]
L_{6,1}(\bphi)  & = &  \phi_{1}^5+52 \phi_{1}^3 \phi_{2}+136 \phi_{1} \phi_{2}^2+192 \phi_{1}^2 \phi_{3}+180 \phi_{2} \phi_{3}+264 \phi_{1} \phi_{4}+120 \phi_{5}    \\[0.5mm]
L_{6,2}(\bphi)  & = &  31 \phi_{1}^4+292 \phi_{1}^2 \phi_{2}+136 \phi_{2}^2+342 \phi_{1} \phi_{3}+144 \phi_{4}    \\[0.5mm]
L_{6,3}(\bphi)  & = &  90 \phi_{1}^3+240 \phi_{1} \phi_{2}+90 \phi_{3}    \\[0.5mm]
L_{6,4}(\bphi)  & = &  65 \phi_{1}^2+40 \phi_{2}    \\[0.5mm]
L_{6,5}(\bphi)  & = &  15 \phi_{1}    \\[0.5mm]
L_{6,6}(\bphi)  & = &  1    \\[3mm]
L_{7,1}(\bphi)  & = &  \phi_{1}^6+114 \phi_{1}^4 \phi_{2}+720 \phi_{1}^2 \phi_{2}^2+272 \phi_{2}^3+732 \phi_{1}^3 \phi_{3}+2304 \phi_{1} \phi_{2} \phi_{3}+540 \phi_{3}^2  \\
   & & \quad +\, 1824 \phi_{1}^2 \phi_{4}+1248 \phi_{2} \phi_{4}+1920 \phi_{1} \phi_{5}+720 \phi_{6}    \\[0.5mm]
L_{7,2}(\bphi)  & = &  63 \phi_{1}^5+1176 \phi_{1}^3 \phi_{2}+1848 \phi_{1} \phi_{2}^2+2436 \phi_{1}^2 \phi_{3}+1680 \phi_{2} \phi_{3}+2352 \phi_{1} \phi_{4}+840 \phi_{5}    \\[0.5mm]
L_{7,3}(\bphi)  & = &  301 \phi_{1}^4+1792 \phi_{1}^2 \phi_{2}+616 \phi_{2}^2+1512 \phi_{1} \phi_{3}+504 \phi_{4}    \\[0.5mm]
L_{7,4}(\bphi)  & = &  350 \phi_{1}^3+700 \phi_{1} \phi_{2}+210 \phi_{3}    \\[0.5mm]
L_{7,5}(\bphi)  & = &  140 \phi_{1}^2+70 \phi_{2}    \\[0.5mm]
L_{7,6}(\bphi)  & = &  21 \phi_{1}    \\[0.5mm]
L_{7,7}(\bphi)  & = &  1    
\end{eqnarray*}

\subsection{Lah symmetric functions of positive type}

\vspace*{-5mm}

%
\begin{eqnarray*}
   L_{1,1}^{(\infty)+}  & = &  1    \\[3mm]
   L_{2,1}^{(\infty)+}  & = &  m_{1}   \\[0.5mm]
   L_{2,2}^{(\infty)+}  & = &  1    \\[3mm]
   L_{3,1}^{(\infty)+}  & = &  m_{2}+4 m_{11}   \\[0.5mm]
   L_{3,2}^{(\infty)+}  & = &  3 m_{1}   \\[0.5mm]
   L_{3,3}^{(\infty)+}  & = &  1    \\[3mm]
   L_{4,1}^{(\infty)+}  & = &  m_{3}+11 m_{21}+36 m_{111}   \\[0.5mm]
   L_{4,2}^{(\infty)+}  & = &  7 m_{2}+22 m_{11}   \\[0.5mm]
   L_{4,3}^{(\infty)+}  & = &  6 m_{1}   \\[0.5mm]
   L_{4,4}^{(\infty)+}  & = &  1    \\[3mm]
   L_{5,1}^{(\infty)+}  & = &  m_{4}+66 m_{22}+26 m_{31}+196 m_{211}+576 m_{1111}   \\[0.5mm]
   L_{5,2}^{(\infty)+}  & = &  15 m_{3}+105 m_{21}+300 m_{111}   \\[0.5mm]
   L_{5,3}^{(\infty)+}  & = &  25 m_{2}+70 m_{11}   \\[0.5mm]
   L_{5,4}^{(\infty)+}  & = &  10 m_{1}   \\[0.5mm]
   L_{5,5}^{(\infty)+}  & = &  1    \\[3mm]
   L_{6,1}^{(\infty)+}  & = &  m_{5}+302 m_{32}+57 m_{41}+1898 m_{221}+848 m_{311}+5244 m_{2111}+14400 m_{11111}   \\[0.5mm]
   L_{6,2}^{(\infty)+}  & = &  31 m_{4}+906 m_{22}+416 m_{31}+2446 m_{211}+6576 m_{1111}   \\[0.5mm]
   L_{6,3}^{(\infty)+}  & = &  90 m_{3}+510 m_{21}+1350 m_{111}   \\[0.5mm]
   L_{6,4}^{(\infty)+}  & = &  65 m_{2}+170 m_{11}   \\[0.5mm]
   L_{6,5}^{(\infty)+}  & = &  15 m_{1}   \\[0.5mm]
   L_{6,6}^{(\infty)+}  & = &  1    \\[3mm]
   L_{7,1}^{(\infty)+}  & = &  m_{6}+2416 m_{33}+1191 m_{42}+120 m_{51}+28470 m_{222}+13644 m_{321}+3228 m_{411} \\
       & & \quad +\, 75216 m_{2211}+36240 m_{3111}+197856 m_{21111}+518400 m_{111111}   \\[0.5mm]
   L_{7,2}^{(\infty)+}  & = &  63 m_{5}+6006 m_{32}+1491 m_{41}+31794 m_{221}+15624 m_{311}+82152 m_{2111}+211680 m_{11111}   \\[0.5mm]
   L_{7,3}^{(\infty)+}  & = &  301 m_{4}+6006 m_{22}+2996 m_{31}+15316 m_{211}+38976 m_{1111}   \\[0.5mm]
   L_{7,4}^{(\infty)+}  & = &  350 m_{3}+1750 m_{21}+4410 m_{111}   \\[0.5mm]
   L_{7,5}^{(\infty)+}  & = &  140 m_{2}+350 m_{11}   \\[0.5mm]
   L_{7,6}^{(\infty)+}  & = &  21 m_{1}   \\[0.5mm]
   L_{7,7}^{(\infty)+}  & = &  1    
\end{eqnarray*}

\subsection{Lah symmetric functions of negative type}

\vspace*{-5mm}

\begin{eqnarray*}
   L_{1,1}^{(\infty)-}  & = &  1    \\[3mm]
   L_{2,1}^{(\infty)-}  & = &  m_{1}    \\[0.5mm]
   L_{2,2}^{(\infty)-}  & = &  1    \\[3mm]
   L_{3,1}^{(\infty)-}  & = &  3 m_{2}+4 m_{11}    \\[0.5mm]
   L_{3,2}^{(\infty)-}  & = &  3 m_{1}    \\[0.5mm]
   L_{3,3}^{(\infty)-}  & = &  1    \\[3mm]
   L_{4,1}^{(\infty)-}  & = &  15 m_{3}+25 m_{21}+36 m_{111}    \\[0.5mm]
   L_{4,2}^{(\infty)-}  & = &  15 m_{2}+22 m_{11}    \\[0.5mm]
   L_{4,3}^{(\infty)-}  & = &  6 m_{1}    \\[0.5mm]
   L_{4,4}^{(\infty)-}  & = &  1    \\[3mm]
   L_{5,1}^{(\infty)-}  & = &  105 m_{4}+250 m_{22}+210 m_{31}+380 m_{211}+576 m_{1111}    \\[0.5mm]
   L_{5,2}^{(\infty)-}  & = &  105 m_{3}+195 m_{21}+300 m_{111}    \\[0.5mm]
   L_{5,3}^{(\infty)-}  & = &  45 m_{2}+70 m_{11}    \\[0.5mm]
   L_{5,4}^{(\infty)-}  & = &  10 m_{1}    \\[0.5mm]
   L_{5,5}^{(\infty)-}  & = &  1    \\[3mm]
   L_{6,1}^{(\infty)-}  & = &  945 m_{5}+3010 m_{32}+2205 m_{41}+5810 m_{221}+4760 m_{311}+9156 m_{2111}+14400 m_{11111}    \\[0.5mm]
   L_{6,2}^{(\infty)-}  & = &  945 m_{4}+2590 m_{22}+2100 m_{31}+4130 m_{211}+6576 m_{1111}    \\[0.5mm]
   L_{6,3}^{(\infty)-}  & = &  420 m_{3}+840 m_{21}+1350 m_{111}    \\[0.5mm]
   L_{6,4}^{(\infty)-}  & = &  105 m_{2}+170 m_{11}    \\[0.5mm]
   L_{6,5}^{(\infty)-}  & = &  15 m_{1}    \\[0.5mm]
   L_{6,6}^{(\infty)-}  & = &  1    \\[3mm]
   L_{7,1}^{(\infty)-}  & = &  10395 m_{6}+48160 m_{33}+42525 m_{42}+27720 m_{51}+122010 m_{222}+97860 m_{321}  \\
      & & \quad +\, 69300 m_{411}+197904 m_{2211}+158928 m_{3111}+320544 m_{21111}+518400 m_{111111}    \\[0.5mm]
   L_{7,2}^{(\infty)-}  & = &  10395 m_{5}+38430 m_{32}+26775 m_{41}+79170 m_{221}+63000 m_{311}+129528 m_{2111} \\
      & & \quad +\, 211680 m_{11111}    \\[0.5mm]
   L_{7,3}^{(\infty)-}  & = &  4725 m_{4}+14350 m_{22}+11340 m_{31}+23660 m_{211}+38976 m_{1111}    \\[0.5mm]
   L_{7,4}^{(\infty)-}  & = &  1260 m_{3}+2660 m_{21}+4410 m_{111}    \\[0.5mm]
   L_{7,5}^{(\infty)-}  & = &  210 m_{2}+350 m_{11}    \\[0.5mm]
   L_{7,6}^{(\infty)-}  & = &  21 m_{1}    \\[0.5mm]
   L_{7,7}^{(\infty)-}  & = &  1    
\end{eqnarray*}

\end{document}